\documentclass{amsart}

\usepackage{amsmath,amssymb,amsbsy,amsthm,amsfonts}
\DeclareMathOperator*{\esslim}{ess\,lim}
\allowdisplaybreaks 

\usepackage{mathtools}
\usepackage{xfrac}
\usepackage{esint}
\usepackage{mathrsfs}
\usepackage{a4wide}
\usepackage[foot]{amsaddr}
\usepackage{verbatim}

\usepackage{graphicx}
	\DeclareGraphicsExtensions{.pdf,.png,.jpg}
\usepackage{pgfplots}
	\pgfplotsset{compat=1.3}
\usepackage{placeins}
\usepackage{tabularx}
\usepackage{enumitem}
\usepackage[font=footnotesize,margin=0.2cm]{caption}

\usepackage[left=3cm,right=3cm,top=3.5cm,bottom=2.5cm]{geometry}

\usepackage[obeyDraft,textsize=footnotesize,
bordercolor=gray]{todonotes}
\reversemarginpar
\usepackage{dsfont}
\definecolor{cobalt}{rgb}{0.0, 0.28, 0.67}
\definecolor{darkcerulean}{rgb}{0.03, 0.27, 0.49}
\usepackage[draft=false,colorlinks=true]{hyperref}
\hypersetup{citecolor=blue,linkcolor=darkcerulean}
\usepackage[all]{hypcap}

\makeatletter
\newcommand{\setword}[2]{%
  \phantomsection
  #1\def\@currentlabel{\unexpanded{#1}}\label{#2}%
}
\makeatother

\newtheorem{Theorem}{Theorem}[section]
\newtheorem{Example}[Theorem]{Example}
\newtheorem{Assumption}[Theorem]{Assumption}
\newtheorem{Lemma}[Theorem]{Lemma}

\newtheorem{Remark}[Theorem]{Remark}
\newtheorem{Definition}[Theorem]{Definition}
\newtheorem{Corollary}[Theorem]{Corollary}

\numberwithin{equation}{section}

\makeatletter
\def\subsubsection{\@startsection{subsubsection}{3}%
  \z@{.5\linespacing\@plus.7\linespacing}{-.5em}%
  {\normalfont\bfseries
  }}
 \makeatother
\def\abs#1{\left| #1 \right|}
\def\norm#1{\left|\!\left| #1 \right|\!\right|}
\def\ll{\left\langle}
\def\rr{\right\rangle}

\def\tens#1{\pmb{\mathsf{#1}}}
\def\vec#1{\boldsymbol{#1}}

\def\sym{\mathop{\mathrm{sym}}\nolimits}
\def\Rds{\mathbb{R}^{d \times d}_{\sym}}

\def\diam{\mathop{\mathrm{diam}}\nolimits} 		

\def\loc{\mathop{\mathrm{loc}}\nolimits}
\def\supp{\mathop{\mathrm{supp}}\nolimits}
\def\dom{\mathop{\mathrm{dom}}\nolimits}
\def\conv{\mathop{\mathrm{conv}}\nolimits}

\def\diver{\mathop{\mathrm{div}}\nolimits} 

\def\d{{\rm d}}

\def\wconv{\rightharpoonup}
\def\cpemb{\hookrightarrow \hookrightarrow}
\def\ctemb{\hookrightarrow}

\def\Cw{C_w}
\def\Leb{\mathrm{L}} 
\def\Sob{\mathrm{W}} 
\def\l2d{ \Leb^2_{\diver}(\Omega)^d}

\def\cS{\mathsf{S}}

\def\b0{\vec{0}}
\def\be{\vec{e}}
\def\bf{\vec{f}}
\def\bg{\vec{g}}
\def\bu{\vec{u}}
\def\bv{\vec{v}}
\def\bw{\vec{w}}
\def\bx{\vec{x}}
\def\bz{\vec{z}}
\def\bal{\vec{\alpha}}
\def\bxi{\vec{\xi}}
\def\bF{\vec{F}}
\def\bU{\vec{U}}
\def\bV{\vec{V}}
\def\bW{\vec{W}}

\def\B0{\tens{0}}
\def\Bv{\tens{v}}
\def\BA{\tens{A}}
\def\BB{\tens{B}}

\def\BD{\tens{D}}
\def\BF{\tens{F}}
\def\BG{\tens{G}}
\def\BH{\tens{H}}
\def\BS{\tens{S}}

\def\Vd{\mathbb{V}^n_{\diver}}
\def\Vn{\mathbb{V}^n}
\def\Qn{\mathbb{Q}^n}
\def\Pn{\Pi^n}

\def\T{\mathcal{T}}
\def\bad{\mathscr{B}}

\def\uk{\bU^{\kappa}} 
\def\oU{\overline{\bU}} 
\def\wU{\widetilde{\bU}} 
\def\L{\lambda} 
\def\Xod{\mathrm{X}^q_{\diver}(\Omega)}
\def\Xq{\mathrm{X}^q(Q)}
\def\Xqd{\mathrm{X}^q_{\diver}(Q)}
\def\wq{\hat{q}} 

\begin{document}

\title[FEM for implicitly constituted fluids]{Fully discrete finite element approximation of unsteady flows of implicitly constituted incompressible fluids
}

\author[E. S\"uli]{Endre S\"uli}
\address{Mathematical Institute,
University of Oxford, Woodstock Road, Oxford, OX2 6GG, UK}
\email{endre.suli@maths.ox.ac.uk}

\author[T. Tscherpel]{Tabea Tscherpel}
\thanks{This work was supported by Engineering and Physical Sciences Research Council [EP/L015811/1].
}
\email{tabea.tscherpel@maths.ox.ac.uk}
\date{\today}

\maketitle

\begin{abstract}
Implicit constitutive theory provides a very general framework for fluid flow models, including both Newtonian and generalized Newtonian fluids, where the Cauchy stress tensor and the rate of strain tensor are assumed to be related by
an implicit relation associated with a maximal monotone graph. For incompressible unsteady flows of such fluids, subject to a homogeneous Dirichlet boundary condition on a Lipschitz polytopal domain $\Omega \subset \mathbb{R}^d$, $d \in \{2,3\}$, we investigate a fully-discrete approximation scheme, using a spatial mixed finite element approximation combined with backward Euler time-stepping. We show convergence of a subsequence of approximate solutions, when the velocity field belongs to the space of solenoidal functions contained in $\Leb^\infty(0,T;\Leb^2(\Omega)^d)\cap \Leb^q(0,T;\Sob^{1,q}_0(\Omega)^d)$, provided that $q\in \big(\frac{2d}{d+2},\infty\big)$, which is the maximal range for $q$ with respect to existence of weak solutions. This is achieved by a technique based on splitting and regularizing, the use of a solenoidal parabolic Lipschitz truncation method, a local Minty-type monotonicity result, and
various weak compactness results.
\end{abstract}




\section{Introduction}

In the mechanics of viscous incompressible fluids typical constitutive relations relate the shear stress tensor to the rate of strain tensor through an explicit functional relationship.
In the case of a Newtonian fluid the relationship is linear, and in the
case of generalized Newtonian fluids it is usually a power-law-like nonlinear, but still explicit, functional relation.
Implicit constitutive theory was introduced in order to describe a wide range of non-Newtonian rheology, by admitting
implicit and discontinuous constitutive laws, see \cite{R.2003,R.2006}. The existence of weak solutions to mathematical models of this kind was explored, for example, in \cite{BGMS.2009, BGMS.2012} for steady and unsteady flows, respectively.
The aim of the present paper is to construct a fully discrete numerical approximation scheme, in the unsteady case, for a class of such implicitly constituted models, where the shear stress and rate of strain tensors are related through a (possibly discontinuous) maximal monotone graph. The scheme is based on a spatial mixed finite element approximation and a backward Euler discretization with respect to the temporal variable. We will show weak convergence (up to subsequences) to a weak solution of the model. The mathematical ideas contained in the paper are motivated by the existence theory formulated, in the unsteady case, in \cite{BGMS.2012}, and the convergence theory for finite element approximations of steady implicitly constituted fluid flow models developed in \cite{DKS.2013}.

\subsection{Implicit Constitutive Law \nopunct}

\subsubsection*{Statement of the Problem}
Let $\Omega \subset \mathbb{R}^d$ with $d \geq 2$ be a bounded Lipschitz domain and denote by $Q = (0,T)\times \Omega$ the parabolic cylinder for a given final time $T \in (0, \infty)$. Furthermore, let $\bf\colon Q \to \mathbb{R}^d$ be a given external force and let $\bu_0\colon \Omega \to \mathbb{R}^d$ be an initial velocity field.
We seek a velocity field $\bu: Q \to \mathbb{R}^d$, a pressure $\pi: Q \to \mathbb{R}$, and a stress tensor field $\BS: Q \to \Rds$ satisfying the balance law of linear momentum and the incompressibility condition:
\begin{alignat}{2} 
\begin{aligned}
\label{P-1}
\partial_t \bu + \diver(\bu \otimes \bu)
-\diver\BS &= -\nabla \pi + \bf, \qquad \quad &&\text{ on } (0,T) \times \Omega,\\
 \diver\bu &= 0, \qquad \quad &&\text{ on } (0,T) \times \Omega,
\end{aligned}
\end{alignat}
subject to the following initial condition and no-slip boundary condition:
\begin{alignat}{2} \label{P-2}
\bu(0,\cdot)&=\bu_0(\cdot)  \qquad &&\text{ in } \Omega,\\
\bu &= \b0  \qquad \quad &&\text{ on } (0,T) \times \partial \Omega. \label{P-3}
\end{alignat}
In order to close the system we need to impose a relation, the so-called \textit{constitutive law},
\begin{align} \label{P-4}
\BG(\cdot,\BD\bu, \BS) = \B0,
\end{align}
between the stress tensor $\BS$ and the symmetric gradient $\BD \bu = \frac{1}{2}(\nabla \bu + (\nabla \bu)^{\top})$, which represents the shear rate of the fluid. In the following we will refer to the problem consisting of \eqref{P-1}--\eqref{P-4} as \setword{\textbf{(P)}}{def:P}, and will use the notation $\bz:=(t,\bx)\in Q$.

The relation $\BG$ may be fully implicit and we assume that $\BG$ can be identified with a maximal monotone graph $\mathcal{A}(\bz) \subset \Rds \times \Rds$, for $\bz \in Q$, as
\begin{align*}
\BG(\bz,\BD \bu(\bz), \BS(\bz)) = \B0 \quad \Leftrightarrow \quad (\BD\bu(\bz), \BS(\bz)) \in \mathcal{A}(\bz),
\end{align*}
where $\mathcal{A}(\cdot)$ satisfies the following assumption, similarly as in \cite[p.~110]{BGMS.2009} and \cite[Sec.~1.2]{BGMS.2012}.

\begin{Assumption}[{Properties of $\mathcal{A}(\cdot)$}]\label{assump-A}
We assume that $\mathcal{A}(\cdot) \subset \Rds \times \Rds$ satisfies the following conditions for  a.e. $\bz \in Q $:
\begin{enumerate}[label={(A\arabic*)}]
\item \label{itm:A1} $(\tens{0},\tens{0})\in \mathcal{A}(\bz)$;
\item \label{itm:A2} $\mathcal{A}(\bz)$ is a monotone graph, i.e., for all $(\BD_1,\BS_1), (\BD_2,\BS_2) \in \mathcal{A}(\bz)$,
\begin{align*}
(\BD_1-\BD_2)\colon(\BS_1-\BS_2)\geq 0;
\end{align*}
\item \label{itm:A3}  $\mathcal{A}(\bz)$ is a maximal monotone graph, i.e., $(\BD, \BS) \in \Rds \times \Rds$ and
\begin{align*}
(\overline{\BD}-\BD)\colon(\overline{\BS}-\BS)\geq 0 \quad \text{ for all } (\overline{\BD}, \overline{\BS}) \in \mathcal{A}(\bz),
\end{align*}
implies that $(\BD, \BS)\in \mathcal{A}(\bz)$;
\item \label{itm:A4} There exists a constant $c_*>0$, a nonnegative function $g \in \Leb^1(Q)$ and $q \in (1, \infty)$ such that
\begin{align*}
\BD \colon \BS \geq -g(\bz)+ c_*(\abs{\BD}^{q}+ \abs{\BS}^{q'}) \quad \text{ for all }\, (\BD, \BS)\in \mathcal{A}(\bz),
\end{align*}
where $q'$ is the H\"older conjugate of $q$;
\item \label{itm:A5} $\mathcal{A}$ is $\mathcal{L}(Q)-(\mathcal{B}(\Rds)\otimes \mathcal{B}(\Rds))$ measurable, where $\mathcal{L}(Q)$ denotes the set of all Lebesgue measurable subsets of $Q$ and $\mathcal{B}(\Rds)$ denotes the set of all Borel subsets of $\Rds$.
\end{enumerate}
\end{Assumption}

\begin{Remark}[{Properties of $\mathcal{A}(\cdot)$}]\hfill
\begin{enumerate}[label={(\roman*)}]
\item
In \cite{BGMS.2012} the authors phrase the condition (A4) in the more general context of Orlicz--Sobolev spaces. Here we will restrict ourselves to the usual Sobolev setting. 
\item Conditions for \ref{itm:A5} to be satisfied are given in \cite[(A5)(i),(ii), p.~110]{BGMS.2009}.
\end{enumerate}
\end{Remark}

This framework covers explicit relations, including Newtonian fluids, where $q=2$, and $q$-fluids describing shear-thinning and shear-thickening behaviour, for $1<q<2$ and $q>2$, respectively. Also relations, where the the stress is a set-valued or discontinuous function of the symmetric gradient, as for Bingham and Herschel--Bulkley fluids, are included, which is shown in \cite[Lem.~1.1]{BGMS.2012}. Furthermore, fully implicitly constituted fluids are covered and the constitutive relation is allowed to depend on $(t,\bx)\in Q$.

\subsection{Overview of the Context \nopunct}\hfill\\
For explicit constitutive laws the list of existence results is long, see amongst others \cite{La.1969,L.1969,FMS.2003,DRW.2010} and the references therein. For those constitutive laws the main challenge is the lack of admissibility for small $q$, caused by the presence of the convective term. In the most recent results this difficulty was overcome by the use of the Lipschitz truncation method.

In the case of implicitly constituted fluids, the existence of weak solutions for $q>\frac{2d}{d+2}$ for steady and unsteady flows was proved in \cite{BGMS.2009, BGMS.2012}. The restriction on $q$ is required to ensure compactness of the embedding $\Sob^{1,q}(\Omega)\ctemb \Leb^2(\Omega)$, and implies that the convective term in the weak formulation is well-defined.
In \cite{BGMS.2012}, a Navier slip boundary condition and $C^{1,1}$ regularity of $\partial \Omega$ were assumed to avoid technicalities related to lack of regularity of the pressure in the unsteady case.
Due to the weak structural assumptions the existence of short-time strong solutions and uniqueness cannot be expected to hold in general. The proof in \cite{BGMS.2012} is constructive and is based on a three-level approximation using finite-dimensional Galerkin subspaces spanned by eigenfunctions of higher order elliptic operators. These Galerkin spaces are not available for practical computations and therefore we take an alternative route here in the construction of a numerical method for the
problem and for its convergence analysis.
Here we shall consider a mixed finite element approximation under minimal regularity hypotheses, so we can only hope for qualitative convergence results rather than qualitative error bounds in terms of the spatial and temporal discretization
parameters. The approximation scheme will be constructed for a regularized version of the equations, and after passing to the limit with the discretization parameters we shall pass to the limit with the regularization parameter.

Concerning the numerical analysis of implicitly constituted fluid flow models, to the best of our knowledge the only results available are those contained in \cite{DKS.2013, KS.2016}, which deal with the steady case.
By means of a discrete Lipschitz truncation method and various weak compactness results the authors of \cite{DKS.2013} prove the convergence of a large class of mixed finite element methods for $q>\frac{2d}{d+1}$ for discretely divergence-free finite element functions for the velocity, and for $q>\frac{2d}{d+2}$ for exactly (i.e., pointwise) divergence-free finite element functions for the velocity field. In the case of discretely  divergence-free mixed finite element approximations the more demanding requirement, that $q>\frac{2d}{d+1}$, arises from the (numerical) modification of the trilinear form associated with the convective term in the weak formulation of the problem,
whose purpose is to reinstate the skew-symmetry of the trilinear form, lost in the course of the spatial approximation. In \cite{KS.2016} an a~posteriori analysis is performed for implicitly constituted fluid flow models using discretely divergence-free finite element functions. In the unsteady case no convergence result is available for implicitly constituted fluid flow models, and even those contributions that are focussed on explicit constitutive laws, such as \cite{CHP.2010}, for example, are restricted to the case when $q>\frac{2(d+1)}{d+2}$.

\subsection{Aim and Main Result \nopunct}\hfill\\
Our objective here is to establish a convergence result for implicitly constituted fluids in the unsteady case for the whole range $q>\frac{2d}{d+2}$. More specifically, the aim is to show weak convergence (up to subsequences) of the sequence of approximate solutions to a weak solution of the problem, the main challenges being to deal with the implicit, possibly discontinuous, relation between stress and shear rate, and with small exponents $q$ arising in the coercivity and boundedness assumption on the implicit relation. To this end, we shall consider a three-level approximation, consisting of approximating the potentially multivalued graph $\mathcal{A}$ with a single-valued implicit functional relationship between the shear stress tensor and the
rate of strain tensor, performing a regularization of the resulting model, which will enable us to cover the entire range
of $q > \frac{2d}{d+2}$, and then constructing a fully discrete approximation of the regularized model.
%

%
The main contribution of the paper is therefore the following.
Let $\Omega$ be a Lipschitz polytopal domain, $q>\frac{2d}{d+2}$, and assume that we have a pair of inf-sup stable finite element spaces. Also we assume that a suitable approximation of the graph $\mathcal{A}$ is available, which will be constructed below. Then, a sequence of approximate solutions to the fully discrete problem exists and the corresponding sequence of finite element approximations converges weakly, up to subsequences, to a weak solution of problem \ref{def:P}, when first taking the graph approximation limit, then the spatial and temporal discretization limits, and finally the regularization limit.
The precise formulation of this result is contained in Theorem~\ref{Thm:main} and the notion of weak solution is given in Definition~\ref{def:w-sol}.
The main tools in the proof are a Minty type local monotonicity result proved in \cite{BGMS.2012}, and the solenoidal parabolic Lipschitz truncation constructed in \cite{BDS.2013} to overcome the admissibility problem for small $q$.

The paper is structured as follows. Section \ref{sec:2} provides the analytical setting, including the graph approximation and the Lipschitz truncation. Section \ref{Sec:FEM} describes the finite element approximation, the approximation of the convective term and the time stepping. In Section \ref{sec:4} we first introduce the approximation levels in detail before giving the convergence proof.

\section{Analytical Preliminaries}\label{sec:2}

By $\Rds$ we denote the set of all real-valued symmetric $d\times d$-matrices and we use $\colon$ for the Frobenius scalar product in $\mathbb{R}^{d\times d}$.
For $\omega \subset \mathbb{R}^d$ we denote by $\abs{\omega}$ the $d$-dimensional Lebesgue measure of $\omega$. By $\mathds{1}_{\omega}$ we denote the characteristic function  of the set $\omega$.
For the partial derivatives with respect to time we use the shorthand notation $\partial_t f \coloneqq \frac{\partial f}{\partial t}$.

For $\omega\subset \mathbb{R}^d$ open and $p\in [1, \infty)$ let $(\Leb^p(\omega), \norm{\cdot}_{\Leb^p(\omega)})$ be the standard Lebesgue space of $p$-integrable functions and the space of essentially bounded functions when $p=\infty$. For $s\in \mathbb{N}$ let $(\Sob^{s,p}(\omega), \norm{\cdot}_{\Sob^{s,p}(\omega)})$ be the respective Sobolev spaces.
For spaces of vector-valued and tensor-valued functions we use superscripts $d$ and $d \times d$, respectively (except for in norms). By $\Leb^p_0(\omega)$ we denote the set of functions in $\Leb^p(\omega)$ with zero mean integral.

For a general Banach space $(X, \norm{\cdot}_X)$, the dual space consisting of all continuous linear functionals on $X$ is denoted by $X'$ and the dual pairing is denoted by $\ll f, g \right\rangle_{X', X},$ if $f \in X'$ and $g \in X$.
 If $X$ is a space of functions defined on $\omega$ then we denote the dual pairing by
$\ll f,g\right\rangle_\omega \coloneqq \ll f,g\right\rangle_{X', X},$
in case the space $X$ is known from the context. We also use this notation for the integral of the scalar product $f \cdot g$ of two functions $f$ and $g$, provided that $f \cdot  g \in \Leb^1(\omega)$.
Furthermore, if $\omega \subset \mathbb{R}^d$ is measurable and $0<\abs{\omega}<\infty$, then we denote
$
\fint_{\omega} f(\bx) \, \d \bx \coloneqq \frac{1}{\abs{\omega}}\int_{\omega} f(\bx) \, \d \bx.
$

Now let $\Omega \subset \mathbb{R}^d$ be a bounded domain, let $T \in (0, \infty)$ and $Q=(0,T)\times \Omega$.
Denote by $C^\infty_0(\Omega)$ the set of all smooth and compactly supported functions on $\Omega$ and by $C^{\infty}_{0, \diver}(\Omega)^d$ the set of all functions in $C^\infty_0(\Omega)^d$ with vanishing divergence. Analogously, define $C^{\infty}_{0}(Q)^d$ and $C^{\infty}_{0,\diver}(Q)^d$. Define $\Sob^{1,p}_0(\Omega)\coloneqq \overline{C^{\infty}_0(\Omega)}^{\norm{\cdot}_{\Sob^{1,p}(\Omega)}}$,  for $p\in [1,\infty)$ and $\Sob^{1,\infty}_0(\Omega)\coloneqq \Sob^{1,1}_0(\Omega) \cap \Sob^{1,\infty}(\Omega)$.
For a given $p \in (1, \infty)$ we let the H\"{o}lder exponent $p'$ be defined by $\frac{1}{p}+ \frac{1}{p'}=1$. Then, if $p\in (1, \infty)$, $\Leb^{p'}(\Omega)$ is the dual space of $\Leb^p(\Omega)$ and $\Sob^{-1, p'}(\Omega)$ will denote the dual space of $\Sob^{1,p}_0(\Omega)$.
Further, we define the spaces of divergence-free functions
\begin{align*}
 \l2d  &\coloneqq \overline{C^{\infty}_{0,\diver}(\Omega)^d}^{\norm{\cdot}_{\Leb^{2}(\Omega)}},\\
\Sob^{1,p}_{0, \diver}(\Omega)^d&\coloneqq \overline{C^{\infty}_{0,\diver}(\Omega)^d}^{\norm{\cdot}_{\Sob^{1,p}(\Omega)}}, \quad \text{ for } p\in [1,\infty), \\
\Sob^{1,\infty}_{0,\diver}(\Omega)^d&\coloneqq \Sob^{1,1}_{0, \diver}(\Omega)^d \cap \Sob^{1,\infty}(\Omega)^d. 
\end{align*}

Let $C(\Omega)$ be the set of all continuous real-valued functions on $\Omega$ and $C^{0, \beta}(\Omega)$ the set of all H\"{o}lder continuous functions on $\Omega$ with exponent $\beta\in(0,1]$.
With $C([0,T];X)$ we denote the set of all functions defined on $[0,T]$ taking values in a Banach space $X$, which are continuous (with respect to the strong topology in $X$). Furthermore, we define the space of weakly continuous functions with values in $X$ by
\begin{align*}
\Cw([0,T];X)&\coloneqq \{\bv: [0,T] \to X:\, t \mapsto \ll \bw, \bv(t, \cdot) \right\rangle_{X',X} \in C([0,T]),\, \forall \bw \in X'\}.
\end{align*}
We denote by $\Leb^p(0,T; X)$ the standard Bochner space of $p$-integrable $X$-valued functions.
 We use the notation $\esslim_{t \to 0_+} f(t)$ to indicate that there exists a zero set $N(f)\subset [0,T]$ such that $t \in (0,T)\backslash N(f)$, when considering the limit of $f(t)$, as $t \to 0_+$.

 In the following, $c>0$ will denote a generic constant, which can change from line to line and depends only on the given data unless specified otherwise.

\subsection{Analytic Setting \nopunct}\hfill\\
For $q \in (1, \infty)$ denote
\begin{align}
\label{def:hq}
\wq \coloneqq \max\left( \left(\frac{q(d+2)}{2d} \right)', q \right) = \max\left(
\frac{q(d+2)}{q(d+2)-2d}, q\right) ,
\end{align}
and note that $\wq=q$, if $q \geq \frac{3d+2 }{d+2}$ and $\wq<\infty$ for any $q>\frac{2d}{d+2}$. Also denote
\begin{align}
\label{def:eta}
\eta &\coloneqq \max\left( 2q', \frac{q(d+2)}{d} \right) >2,\\
\label{def:tau}
\tau &\coloneqq \min\left( q', (2q')'\right)>1.
\end{align}
We denote the function spaces
\begin{align}
\label{def:Xq-s}
\Xod &\coloneqq \Sob^{1,q}_{0,\diver}(\Omega)^d \cap \Leb^{2q'}(\Omega)^d,\\
\label{def:Xq-u}
\Xq &\coloneqq \Leb^{q}(0,T;\Sob^{1,q}_0(\Omega)^d) \cap \Leb^{2q'}(Q)^d
,\\
\label{def:Xq-div-u}
\Xqd &\coloneqq \Leb^{q}(0,T;\Sob^{1,q}_{0,\diver}(\Omega)^d) \cap \Leb^{2q'}(Q)^d,
\end{align} 	
and the corresponding norms are given by
\begin{align*}
\norm{\cdot}_{X^q(\Omega)}\coloneqq \norm{\cdot}_{\Sob^{1,q}(\Omega)} + \norm{\cdot}_{\Leb^{2q'}(\Omega)}, \quad
\norm{\cdot}_{\Xq} \coloneqq \norm{\cdot}_{\Leb^{q}(0,T;\Sob^{1,q}(\Omega))} + \norm{\cdot}_{\Leb^{2q'}(Q)}.
\end{align*}

\subsubsection*{Weak Solutions}
In what follows, let $\Omega\subset \mathbb{R}^d$, with $d \in \{2,3\}$, be a bounded Lipschitz domain and for  $T \in (0,\infty)$ denote $Q =(0,T)\times \Omega$. Furthermore, assume that $q \in (1,\infty)$ is given and let $\mathcal{A}(\cdot) \subset \mathbb{R}^{d \times d}_{\sym} \times \mathbb{R}^{d \times d}_{\sym}$ be a monotone graph satisfying Assumption~\ref{assump-A} with respect to $q$.

\begin{Definition}[Weak Solution]\label{def:w-sol}
For a given $\bu_0 \in \l2d$ and $\bf \in \Leb^{q'}(0,T; \Sob^{-1,q'}(\Omega)^d)$
we call $(\bu, \BS)$ a weak solution to problem \ref{def:P}, if
\begin{align*}
\bu &\in \Leb^q(0,T; \Sob^{1,q}_{0, \diver}(\Omega)^d)
 \cap \Cw([0,T];\l2d), \text{ s.t. } \\
 \partial_t \bu &\in \Leb^{\wq'}(0,T;(\Sob^{1,\wq}_{0,\diver}(\Omega)^d)'),
\\
\BS &\in \Leb^{q'}(Q)^{d \times d},
\end{align*}
and
\begin{align}
\ll \partial_t \bu, \bv \rr_{\Omega} - \ll \bu \otimes \bu, \BD \bv \rr_{\Omega}
 +\ll \BS, \BD \bv\rr_{\Omega} 
&= \ll \bf, \bv \rr_{\Omega} \quad
&&\text{ for all } \,\bv \in C^{\infty}_{0,\diver}(\Omega)^d,\label{P-w-1}
\\
&&&
\text{ for a.e. } t \in (0,T),
 \notag
\\
\left(\BD\bu(\bz), \BS(\bz)\right) &\in \mathcal{A}(\bz) \quad &&\text{ for a.e. } \bz \in Q, \label{P-w-2} \\
\esslim_{t\to 0_+}\norm{\bu(t,\cdot)-\bu_0(%
\cdot)}_{\Leb^2(\Omega)} &= 0. \label{P-w-3}
\end{align}
\end{Definition}

We choose a pressure-free notion of weak solution because in the unsteady problem subject to homogeneous Dirichlet boundary conditions on Lipschitz domains one can only expect to establish a distributional (in time) pressure, see \cite[Ch.~III, \S~1, pp.~266]{T.1984} and also \cite{S.1999}.

\subsection{Implicit Constitutive Laws \nopunct}

\subsubsection*{Approximation of $\mathcal{A}$}

The implicit relation encoded by $\mathcal{A}$ can be viewed as a set-valued map. In order to perform the analysis we require a single-valued map and thus a measurable selection $\BS^{\star}$ of the graph $\mathcal{A}$ is chosen, which may have discontinuities.

\begin{Lemma}[{Measurable Selection, \cite[Rem.~1.1, Lem.~2.2]{BGMS.2012}}]\label{sel}
Let $\mathcal{A}(\cdot)\subset \Rds \times \Rds$ be the maximal monotone graph satisfying \ref{itm:A1}--\ref{itm:A5} in Assumption~\ref{assump-A}. 
Then, there exists a measurable selection $\BS^\star\colon Q \times \Rds  \rightarrow \Rds $, i.e.,
\begin{align}\label{sel-def}
(\BB, \BS^\star(\bz,\BB)) \in \mathcal{A}(\bz) \quad \text{ for all }\, \BB \in \Rds,  \text{ for a.e. } \bz \in Q,
\end{align}
and $\BS^{\star}$ is $(\mathcal{L}(Q) \otimes \mathcal{B}(\Rds))-\mathcal{B}(\Rds)$-measurable.
Furthermore, for a.e. $\bz \in Q$, one has that 
\begin{enumerate}[label={(a\arabic*)}]
\item \label{itm:a1}
$\dom \,\BS^{\star}(\bz, \cdot) = \mathbb{R}^{d  \times d}_{\sym}$;
\item \label{itm:a2}
$\BS^{\star}$ is monotone, i.e., for all $\BB_1, \BB_2 \in \mathbb{R}^{d\times d}_{\sym}$,
\begin{align*}
 	(\BS^{\star}(\bz, \BB_1)- \BS^{\star}(\bz, \BB_2))\colon (\BB_1-\BB_2) \geq 0;
\end{align*}
\item \label{itm:a3}
	for any $\BB \in \Rds$ one has that
\begin{align*}
\BB \colon \BS^{\star}(\bz, \BB) \geq -g(\bz)+ c_*(\abs{\BB}^{q}+ \abs{\BS^{\star}(\bz, \BB)}^{q'});
\end{align*}
\item \label{itm:a4} Let $U$ be a dense set in $\Rds$ and let $(\BD, \BS)\in \Rds \times \Rds$. The following are equivalent:
\begin{itemize}
\item[(i)]
	$(\BS-\BS^{\star}(\bz, \BB))\colon(\BD-\BB)\geq 0 \quad \text{ for all } \, \BB \in U$;
\item[(ii)]
	$(\BD, \BS) \in \mathcal{A}(\bz)$;
\end{itemize}
\item \label{itm:a5}
	$\BS^{\star}$ is locally bounded, i.e., for a given $r>0$ there exists a constant $c=c(r)$ such that
\begin{align*}
\abs{\BS^{\star}(\bz, \BA)} \leq c \qquad \text{ for all } \bz\in Q \text{ and for all } \BA \in B_r(\B0) \subset \Rds.
\end{align*}
\end{enumerate}
\end{Lemma}

To show the existence of solutions to any of the approximate problems considered below, continuity of the (approximate) stress tensor is required. Hence, we introduce the following assumptions on a sequence of approximations of the selection $\BS^{\star}$.

\begin{Assumption}[Properties of $\BS^k$, $k \in \mathbb{N}$]\label{Assump-Sk}
Given the selection $\BS^{\star}: Q \times \Rds \to \Rds$ in Lemma~\ref{sel}, assume that there is a sequence $\{\BS^k\}_{k \in \mathbb{N}}$ of Carath\'{e}odory functions $\BS^{k}: Q \times \Rds \to \Rds$ such that:
\begin{enumerate}[label={($\alpha$\arabic*)}]
\item \label{itm:al-2} $\BS^{k}(\bz, \cdot)$ is monotone, i.e., for all $\BA_1, \BA_2 \in \mathbb{R}^{d\times d}_{\sym}$ and for a.e. $\bz \in Q$, we have
\begin{align*}
 (\BS^{k}(\bz, \BA_1)- \BS^{k}(\bz, \BA_2))\colon (\BA_1-\BA_2) \geq 0;
\end{align*}
\item \label{itm:al-3} There exists a constant $\widetilde{c}_*>0$ and a nonnegative function $\widetilde{g} \in \Leb^1(Q)$ such that, for all $k \in \mathbb{N}$, for any $\,\BA \in \Rds$ and for a.e. $\bz\in Q$, one has that
\begin{align*}
\BA \colon \BS^{k}(\bz, \BA) \geq -\widetilde{g}(\bz)+ \widetilde{c}_*\left(\abs{\BA}^{q}+ \abs{\BS^{k}(\bz, \BA)}^{q'}\right);
\end{align*}
\item \label{itm:al-4} Let $U \subset \Rds$ be a dense set. For any sequence $\{\BD^k\}_{k \in \mathbb{N}}$ bounded in $\Leb^{\infty}(Q)^{d \times d}$, for any $\BB \in U$ and all $\varphi \in C^{\infty}_0(Q)$ such that $\varphi \geq 0$, we have
\begin{align*}
\liminf_{k \to \infty} \int_Q \left( \BS^k(\cdot, \BD^k)- \BS^{\star}(\cdot, \BB) \right)\colon \left(\BD^k - \BB \right)\varphi \, \d \bz \geq 0.
\end{align*}
\end{enumerate}
 \end{Assumption}

In the existence proofs in \cite{BGMS.2009, BGMS.2012} and also in \cite{DKS.2013} the approximating sequence $\BS^k$ is chosen as the convolution of the selection $\BS^{\star}$ in the second argument with a mollification kernel.

\begin{Example}[Approximation by Mollification]
Let $\rho\in C^{\infty}_0(\Rds)$ be a mollification kernel, i.e., a nonnegative, radially symmetric function, the support of which is contained in the unit ball $B_1(\B0)\subset \Rds$, and which satisfies
$
\int_{\Rds} \rho(\BA)\, \d \BA = 1.
$
For $k \in \mathbb{N}$ set $\rho^k(\BB)\coloneqq k^{d^2} \rho(k \BB)$ and define the mollification of $\BS^{\star}$ with respect to the last argument by
\begin{align}\label{def:Sk}
\BS^k(\bz, \BB) \coloneqq (\BS^{\star} \ast \rho^k)(\bz, \BB) = \int_{\Rds} \BS^{\star}(\bz, \BA)\rho^k(\BB-\BA)\, \d\BA, \quad \bz \in Q,\; \BB \in \Rds.
\end{align}
Lemma~\ref{Lem:prop-Sk} in the Appendix shows that $\BS^k$ satisfies Assumption~\ref{Assump-Sk}.
\end{Example}

A possibly more practicable approximation based on a piecewise affine interpolant can be used in the case of a radially symmetric selection function $\BS^{\star}$ under additional regularity assumptions. 

\begin{Example}[Approximation by Affine Interpolation]
Assume that $\cS^{\star}: Q \times \mathbb{R}_{\geq 0} \to \mathbb{R}_{\geq 0}$ is a measurable function with $\cS^{\star}(\bz, 0)=0$, for any $\bz \in Q$, such that $\BS^{\star}: Q \times \Rds \to \Rds$, defined by
\begin{align*}
\BS^{\star}(\bz,\BB) = \begin{cases}
\B0 &\quad \text{ if } \BB = \B0, \\
\cS^{\star}(\bz,\abs{\BB}) \frac{\BB}{\abs{\BB}}
 &\quad \text{ else, }
 \end{cases}
\end{align*}
 is a measurable selection of a graph $\mathcal{A}$ satisfying Assumption~\ref{assump-A}. Furthermore, we assume that
\begin{enumerate}[label={(a\arabic*')}]
\item \label{ref:a1-v} $\cS^{\star}(\bz, \cdot): \mathbb{R}_{\geq 0} \to  \mathbb{R}_{\geq 0} $ is monotone  for a.e. $\bz \in Q$;
\item \label{ref:a2-v} Let $I \in \mathbb{N}$ and $\{a_i\}_{i \in \{1, \ldots, I\}}\subset \mathbb{R}>0$, s.t. $a_1 < \cdots < a_I$. Let $a_0=0$ and denote by $A \coloneqq \bigcup_{i \in \{0, \ldots, I\}} a_i$ the finite set of possible discontinuities.
Then, assume that, for a.e. $\bz \in Q$, we have
\begin{alignat*}{2}
\cS^{\star}(\bz, \cdot)|_{(a_{i-1,a_i})} &\in \Sob^{1, \infty}((a_{i-1}, a_i)) \quad &&\text{ for all } i \in \{1, \ldots, I\},\nonumber\\
\cS^{\star}(\bz, \cdot)|_{(a_{I},b)} &\in \Sob^{1, \infty}((a_I, b)) \quad &&\text{ for all } b \in (a_I, \infty).\nonumber
\end{alignat*}
\end{enumerate}
We construct the approximation as follows. Since the set $A$ is finite, there exists a $k_0 \in \mathbb{N}$ such that $\frac{2}{k_0} < \min_{i \in \{1, \ldots, I\}} (a_i - a_{i-1})$. Let $k \in \mathbb{N}$, with $k \geq k_0$ be arbitrary but fixed.
Denote for $i \in \{0, \ldots, I\}$
\begin{align*}
a^k_{i,-} \coloneqq a_i - \frac{1}{k},
\quad a^k_{i,+} \coloneqq a_i + \frac{1}{k},
\quad A^k_i \coloneqq [a^k_{i,-},a^k_{i,+}]
\quad \text{ and } \quad A^k\coloneqq \bigcup_{i \in \{0, \ldots, I\}} A^k_i.
\end{align*}
Let $\bz \in Q$ be arbitrary but fixed. First we extend $\cS^{\star}(\bz, \cdot)$ as an odd function to $\mathbb{R}$ and we still denote the extension by $\cS^{\star}(\bz, \cdot)$. Since the point evaluations $\cS^{\star}(\bz, a^k_{i,\pm})$, for $ i \in \{0, \ldots, I\}$, are well-defined, we can define
\begin{align}\begin{split}
\label{def:ap-Sk}
\overline{\cS}^{k}_{i}(\bz, B)
&\coloneqq
\frac{k}{2}\left(\cS^{\star}(\bz, a^k_{i,-}) \frac{a^k_{i,+}}{a^k_{i,-}} - \cS^{\star}(\bz, a^k_{i,+}) \right)(a^k_{i,-}- B) + \cS^{\star}(\bz, a^k_{i,-}) \frac{B}{a^k_{i,-}},\\
\cS^{k}(\bz, B) &\coloneqq \begin{cases}
\cS^{\star}(\bz, B)\quad &\text{ if } B \notin A^k,\\
\overline{\cS}^{k}_{i}(\bz, B)
 \quad &\text{ if } B \in A^k_i, \quad i \in \{ 0, \ldots, I\}.
\end{cases}
\end{split}
\end{align}
On $A^k_i$ the approximation $\cS^{\star}(\bz, \cdot)$ is the affine interpolant between $\cS^{\star}(\bz, a^k_{i,-})$ and $\cS^{\star}(\bz, a^k_{i,+})$ and otherwise $\cS^{\star}$ is unchanged.

In Lemma~\ref{Lem:prop-Sk-2} and Corollary~\ref{Cor:prop-Sk-2} in the Appendix it is proved that the resulting approximating sequence $\BS^k(\bz, \cdot)$ satisfies Assumption~\ref{Assump-Sk}.
\end{Example}

\subsubsection*{Minty's Trick} The following lemma is one of the crucial tools for the identification of the implicit constitutive law upon passage to the limit.

\begin{Lemma}[{Localized Minty's Trick, \cite[Lem.~2.4]{BGMS.2012} and \cite[Lem.~3.1]{BM.2016}}]\label{minty}
Let $\mathcal{A}(\cdot)$ be a maximal monotone graph satisfying \ref{itm:A2}, \ref{itm:A3} in Assumption~\ref{assump-A} and assume that there are sequences $\{\BS^l\}_{n \in \mathbb{N}}$ and $\{\BD^l\}_{n \in \mathbb{N}}$ and there is a measurable set $\widetilde{Q}\subset Q$ and a $p\in (1,\infty)$ such that
\begin{alignat*}{2}
(\BD^l(\bz), \BS^l(\bz))&\in \mathcal{A}(\bz)\qquad &&\text{ for a.e. } \bz \in \widetilde{Q},\\
\BD^l &\rightharpoonup \BD \qquad &&\text{ weakly in } \Leb^p(\widetilde{Q})^{d \times d}, \\
\BS^l &\rightharpoonup \BS \qquad &&\text{ weakly in } \Leb^{p'}(\widetilde{Q})^{d \times d}, \\
\limsup_{n \to \infty} \ll \BS^l , \BD^l \rr_{\widetilde{Q}} &\leq  \ll \BS , \BD \rr_{\widetilde{Q}}.&&
\end{alignat*}
Then, we have that $
(\BD(\bz), \BS(\bz))\in \mathcal{A}(\bz)$ for a.e. $\bz\in \widetilde{Q}$.
\end{Lemma}

\subsection{Lipschitz Approximation \nopunct}\hfill\\
For small $q\in (1,\infty)$ a weak solution according to Definition~\ref{def:w-sol} is not an admissible test function because of the presence of the convective term. The Lipschitz truncation method helps to identify the implicit relation despite the lack of admissibility.
It first appeared in \cite{AF.1988} and since then the method was further developed and refined in a series of papers, see, e.g., \cite{KL.2002, FMS.2003, DMS.2008,DRW.2010, BDF.2012, BDS.2013, DKS.2013}, to mention just a few.

For a sequence of solutions to a sequence of divergence-form evolution equations a solenoidal parabolic Lipschitz truncation was developed in \cite{BDS.2013}. 
Note that the sets $\bad_{l,j}$ in the following Lemma satisfy $\bad_{l,j} = \mathcal{O}_{l,j} \cap Q_0$, where $\mathcal{O}_{l,j}$ are the ``bad sets'' in the construction in \cite{BDS.2013}.

\begin{Lemma}[{Parabolic Solenoidal Lipschitz Approximation, \cite[Thm.~2.2, Cor.~2.4]{BDS.2013}}\nopunct]\label{Lem:LA-par-div} \hfill\\
Let $p \in (1, \infty)$, $\sigma \in (1, \min(p,p'))$ and let $Q_0 = I_0 \times B_0 \subset \mathbb{R} \times \mathbb{R}^d$ 
  be a parabolic cylinder,  for $d=3$, for an open interval $I_0$ and an open ball $B_0$. Let $\{\bv^l\}_{l \in \mathbb{N}}$ be a sequence
  of weakly divergence-free functions,
   which is converging to zero weakly in $\Leb^{p}(I_0;\Sob^{1,p}(B_0)^d)$, strongly in $\Leb^{\sigma}(Q_0)^d$ and is uniformly bounded in $\Leb^{\infty}(I_0,\Leb^{\sigma}(B_0)^d)$. Consider a sequence $\{\BG^l_1\}_{l \in \mathbb{N}}$, converging to zero weakly in $\Leb^{p'}(Q_0)^{d \times d}$ and a second sequence, $\{\BG^l_2\}_{l \in \mathbb{N}}$, converging to zero strongly in $\Leb^{\sigma}(Q_0)^{d \times d}$. Furthermore, denoting $\BG^l \coloneqq \BG^l_1 + \BG^l_2$, assume that, for any $l \in \mathbb{N}$, the equation
\begin{align}
\ll \partial_t \bv^l, \bxi  \rr_{Q_0} = \ll \BG^l , \nabla \bxi \rr_{Q_0}\quad \text{ for all } \bxi \in C^{\infty}_{0,\diver}(Q_0)^d
\end{align}
is satisfied.
Then, there exists a $j_0 \in \mathbb{N}$,
\begin{itemize}
\item a double sequence $\{\lambda_{l,j}\}_{l,j \in \mathbb{N}} \subset \mathbb{R}$ with $\lambda_{l,j} \in \left[ 2^{2^j}, 2^{2^{j+1}-1}\right]$, for any $l,j \in \mathbb{N}$,
\item a double sequence of open sets $\bad_{l,j} \subset Q_0 $, $l,j \in \mathbb{N}$,
\item a double sequence of functions $\{\bv^{l,j}\}_{l,j \in \mathbb{N}}\subset \Leb^{1}(Q_0)^d$ and
\item a nonnegative function $\zeta \in C^{\infty}_0(\tfrac{1}{6}Q_0)$ such that $\mathds{1}_{\tfrac{1}{8}Q_0} \leq \zeta \leq \mathds{1}_{\tfrac{1}{6}Q_0}$, 
\end{itemize}
such that
\begin{enumerate}[label={(\roman*)}]
\item \label{itm:LA-u-sp}
$\bv^{l,j} \in \Leb^{s}(\tfrac{1}{4}I_0;\Sob^{1,s}_{0,\diver}(\tfrac{1}{6}B_0)^d)$ for all $s \in [1, \infty)$, and $\supp (\bv^{l,j}) \subset \tfrac{1}{6}Q_0$, for any  $j \geq j_0$ and any $l \in \mathbb{N}$;
\item \label{itm:LA-u-eq}
$\bv^{l,j}=\bv^l$ on $\tfrac{1}{8}Q_0 \backslash \bad_{l,j}$, i.e., $\{\bv^{l,j}\neq\bv^l \}\cap \tfrac{1}{8}Q_0 \subset \bad_{l,j}$, for any  $j \geq j_0$ and any $l \in \mathbb{N}$;
\item \label{itm:LA-u-bad}
there exists a constant $c>0$ such that
\begin{align*}
\limsup_{l \to \infty} \lambda_{l,j}^p \abs{\bad_{l,j}} \leq c 2^{-j}\quad \text{ for all } j \geq j_0;
\end{align*}
\item \label{itm:LA-u-est}
there exists a constant $c>0$ such that
\begin{align*}
\norm{\nabla \bv^{l,j}}_{\Leb^{\infty}(\tfrac{1}{4}Q_0)} \leq c \lambda^{l,j}\quad \text{ for all } j \geq j_0 \text{ and all } l \in \mathbb{N};
\end{align*}
\item \label{itm:LA-u-conv}
for any fixed $j \geq j_0$ we have
\begin{align*}
\bv^{l,j} &\to \b0 \quad &&\text{ strongly in }\Leb^{\infty}(\tfrac{1}{4}Q_0)^d,\\ 
\nabla \bv^{l,j} &\rightharpoonup \B0 \quad &&\text{ weakly in }\Leb^{s}(\tfrac{1}{4}Q_0)^{d\times d} \quad \text{ for all } s \in [1,\infty),
\end{align*}
as $l \to \infty$;
\item \label{itm:LA-u-G}
 there exists a constant $c>0$ such that
\begin{align*}
\limsup_{l \to \infty} \abs{ \ll \BG^l, \nabla \bv^{l,j} \rr }\leq c 2^{-j} \quad \text{ for all } j \geq j_0;
\end{align*}
\item \label{itm:LA-u-GH}
there exists a constant $c>0$ such that, for any $\BH \in \Leb^{p'}(\tfrac{1}{6}Q_0)^{d \times d}$, we have that
\begin{align*}
\limsup_{l \to \infty} \abs{ \ll (\BG^l_1 + \BH), \nabla \bv^{l} \zeta \mathds{1}_{\bad_{l,j}^c}\rr }\leq c 2^{-\sfrac{j}{p}} \quad \text{ for all } j \geq j_0.
\end{align*}
\end{enumerate}
\end{Lemma}

The lemma is stated for $d=3$, but according to \cite[Rem.~2.1, p.~2692]{BDS.2013} the result holds for all $d \geq 2$, in fact, with minor modifications of the proof.

\subsection{Continuity and Compactness in Time \nopunct}\hfill

\begin{Lemma}[Parabolic Interpolation]\label{par-interpol}
Let $d\geq 2$, let $\Omega \subset \mathbb{R}^d$ be a Lipschitz domain, let $T \in (0,\infty)$,  $Q=(0,T)\times \Omega$ and let $p >1$. Then, there exists a constant $c>0$ such that
\begin{align}\label{est:par-interpol}
\int_{Q} \abs{v}^{\frac{p(d+2)}{d}} \d \bz \leq c \norm{v}_{\Leb^{p}(0,T;\Sob^{1,p}(\Omega))}^p
\left(  \norm{v}_{\Leb^{\infty}(0,T;\Leb^2(\Omega))}
\right)^{\frac{2p}{d}}
\end{align}
for all $v \in \Leb^p(0,T; \Sob^{1,p}(\Omega))\cap \Leb^{\infty}(0,T; \Leb^2(\Omega))$.
\begin{proof}
For $p\geq d$, the reader is referred to \cite[Prop.~3.1, p.~8]{D.1993}. If $p<d$, then by H\"{o}lder's inequality with $s=\frac{d}{d-p}>1$ and $s'=\frac{d}{p}$ and with Sobolev exponent  $p^* = \frac{dp}{d-p}$ we have that
\begin{align*}
\int_{\Omega} \abs{v}^{\frac{p(d+2)}{d}} \d \bx
= \int_{\Omega} \abs{v}^p \abs{v}^{\frac{2p}{d}} \d \bx
\leq \left( \int_{\Omega} \abs{v}^{\frac{dp}{d-p}} \, \d \bx \right) ^{\frac{d-p}{d}} \left( \int_{\Omega} \abs{v}^{2} \, \d \bx \right) ^{\frac{p}{d}}d
= \norm{v}_{\Leb^{p^*}(\Omega)}^{p} \norm{v}_{\Leb^2(\Omega)}^{\frac{2p}{d}}.
\end{align*}
Then, by integrating over $(0,T)$, further estimating the second factor and using the continuous embedding of $\Sob^{1,p}(\Omega) \ctemb \Leb^{p^*}(\Omega)$, we obtain
\begin{align*}
\int_{Q} \abs{v}^{\frac{p(d+2)}{d}} \d \bx \, \d t
&
\leq c \norm{v}_{\Leb^{p}(0,T;\Leb^{p^*}(\Omega))}^p
\left(  \norm{v}_{\Leb^{\infty}(0,T;\Leb^2(\Omega))}
\right)^{\frac{2p}{d}}
\\
&
\leq \norm{v}_{\Leb^{p}(0,T;\Sob^{1,p}(\Omega))}^p
\left(  \norm{v}_{\Leb^{\infty}(0,T;\Leb^2(\Omega))}
\right)^{\frac{2p}{d}},
\end{align*}
which completes the proof.
\end{proof}
\end{Lemma}

 \begin{Lemma}[{Continuity I}]\label{Lem:cont-1}
Let $Z$ be a reflexive Banach space. Then, any function $v \in \Leb^1(0,T;Z')$, with distributional derivative $\partial_t v \in \Leb^1(0,T;Z')$ is contained in
$\Cw ([0,T];Z')$.
\begin{proof} The proof follows from  \cite[Lem.~1.1, Ch.~III, \S~1, p.~250]{T.1984}. 
\end{proof}
\end{Lemma}

\begin{Lemma}[{Continuity II, \cite[Lem.~1.4, Ch.~III, \S~1, p.~263]{T.1984}
}]\label{Lem:cont-2}
Let $X$ and $Z$ be Banach spaces and let $X$ be reflexive such that the embedding $X \hookrightarrow Z$ is continuous. Then one has that
\begin{align*}
\Leb^{\infty}(0,T;X) \cap \Cw([0,T];Z) \subset \Cw([0,T];X).
\end{align*}
\end{Lemma}


\begin{Lemma}[{Simon, \cite[Thm.~3, p.~80]{S.1987}}] \label{Lem:cp-2}
Let $X,B$ be Banach spaces such that the embedding $X \cpemb B$ is compact. Let $\mathcal{F} \subset \Leb^p(0,T;B)$ for some $p \in [1, \infty]$ and let
\begin{itemize}
\item[(i)] $\mathcal{F}$ be bounded in $\Leb^1_{\loc}(0,T;X)$,
\item[(ii)] $
\int_0^{T-\varepsilon} \norm{f(s+\varepsilon,\cdot)-f(s,\cdot)}_B^p \, \d s \to 0, \quad \text{ as } \varepsilon \to 0, \quad \text{ uniformly for } f \in \mathcal{F}.
$
\end{itemize}
Then, $\mathcal{F}$ is relatively compact in $\Leb^p(0,T;B)$, and in $C([0,T];B)$ if $p=\infty$.
\end{Lemma}

\section{Finite Element Approximation}\label{Sec:FEM}
\subsection{Finite Element Spaces and Assumptions \nopunct}\hfill\\
The setting here is slightly more general than the one in \cite{DKS.2013}.

\begin{Assumption}[Triangulations $\{\mathcal{T}_n\}_{n \in \mathbb{N}}$]\label{assump-triang} Let us assume that $d \geq 2$ and that $ \Omega$ is a bounded  Lipschitz polytopal
 domain. Furthermore, assume that $\{\mathcal{T}_n\}_{n\in\mathbb{N}}$ is a family of simplicial partitions of $\overline{\Omega}$ (in the sense of \cite[Sec.~2.1, p.~38]{C.2002}) 
such that the following conditions hold:
\begin{itemize}
\item[(i)]
Each element $K\in \mathcal{T}_n$ is affine-equivalent to the closed standard reference simplex $\widehat{K}\coloneqq \conv \{ \b0, \be_1, \ldots, \be_d \}\subset\mathbb{R}^d$, i.e., there exists an affine invertible function $\BF_K\colon K \to \widehat{K}$;
\item[(ii)]  
$\{\T_n\}_{n \in \mathbb{N}}$ is shape-regular, i.e., there exists a constant $c_r$
(independent of $n \in \mathbb{N}$) such that
\begin{align*}
\tfrac{h_K}{\rho_K} \leq c_r \quad \text{ for all } K \in \mathcal{T}_n \text{ and all } n \in \mathbb{N},
\end{align*}
where $h_K \coloneqq \diam (K)$ and $\rho_K \coloneqq \sup \{\diam(B)\colon \, B \text{ is a ball contained in } K \}$;
\item[(iii)] The grid size $
h_n\coloneqq \max \{h_K\colon K \in \mathcal{T}_n\}
$ vanishes, as $n \to \infty$.
\end{itemize}
\end{Assumption}

\subsubsection*{Finite Element Spaces} Let $\widehat{\mathbb{P}}_{\mathbb{V}} \subset \Sob^{1,\infty}(\widehat{K})^d$ and let $\widehat{\mathbb{P}}_{\mathbb{Q}} \subset \Leb^\infty(\widehat{K})$ be finite-dimensional function spaces on the reference simplex $\widehat{K}$ (with a slight abuse of notation) as in \cite{DKS.2013}. Further, let $\mathbb{V} \subset C(\overline{\Omega})^d$ and let $\mathbb{Q}\subset \Leb^{\infty}(\Omega)$.
Then we define the conforming finite element spaces $\mathbb{V}^n$ and $\mathbb{Q}^n$ with respect to the partition $\mathcal{T}_n$ by
\begin{align}\label{Vn}
\mathbb{V}^n &\coloneqq \{ \bV \in \mathbb{V}\colon \quad \bV|_K \circ \BF_K^{-1}\in \widehat{\mathbb{P}}_{\mathbb{V}}, \, K \in \mathcal{T}_n \text{ and } \bV|_{\partial \Omega} =\b0 \},\\ \label{Qn}
\mathbb{Q}^n &\coloneqq \{ Q \in \mathbb{Q}\colon \quad Q|_K \circ \BF_K^{-1}\in \widehat{\mathbb{P}}_{\mathbb{Q}}, \, K \in \mathcal{T}_n \}.
\end{align}
Let us also introduce the subspace of discretely divergence-free functions of $\mathbb{V}^n$ and the subspace of zero integral mean functions of $\mathbb{Q}^n$ by
\begin{align}\label{Vd}
\mathbb{V}^n_{\diver} &\coloneqq \{ \bV \in \mathbb{V}^n: \, \ll \diver \bV, Q \right\rangle_{\Omega} =0 \, \text{ for all } Q \in \mathbb{Q}^n\},\\
\label{Qn0}
\mathbb{Q}^n_0 &\coloneqq	\{Q \in \mathbb{Q}^n: \, \int_\Omega Q\, \d \bx =0 \}.
\end{align}
Note that the functions in $\Vd$ are in general not divergence-free, so in general $\Vd \not\subset \Sob^{1, \infty}_{0, \diver}(\Omega)^d$.

\begin{Assumption}[{Approximability, \cite[Assump.~5]{DKS.2013}}]\label{assump-approx}
Assume that for all $p \in [1, \infty)$, we have that
\begin{align}\label{Vn-approx}
\inf_{\bV \in \Vn} \norm{\bv - \bV}_{\Sob^{1,p}(\Omega)} &\to 0, \quad \text{ as } n \to \infty \quad \text{ for all } \bv \in \Sob^{1,p}_0(\Omega)^d,\\ \label{Qn-approx}
\inf_{Q \in \Qn} \norm{h - Q}_{\Leb^{p}(\Omega)} &\to 0, \quad \text{ as } n \to \infty \quad \text{ for all } h  \in \Leb^{p}_0(\Omega).
\end{align}
\end{Assumption}

\subsubsection*{Projectors}
For the convergence analysis we require certain projectors to the respective finite element spaces and suitable assumptions on them. Since we do not need local stability of the projector $\Pn$, we assume less than in \cite{DKS.2013}.
\begin{Assumption}[Projector $\Pn$] \label{assump-Pin}\hfill\\
Assume that for each $n \in \mathbb{N}$ there exists a linear  projector $\Pn \colon \Sob^{1,1}_0(\Omega)^d \to \mathbb{V}^n$ such that:
{\begin{enumerate}[label=(\roman*)]
\item \label{itm:Pin-div} (preservation of the divergence in $(\Qn)^{'}$)
for any $\bv \in \Sob^{1,1}_0(\Omega)^d$ one has that
\begin{align*}
\ll \diver \bv, Q \rr_\Omega = \ll \diver (\Pn \bv), Q \rr_\Omega \qquad \text{ for all } Q \in \Qn;
\end{align*}
\item \label{itm:Pin-stab}  ($\Sob^{1,p}$-stability)
 for any $p \in (1,\infty)$
  there exists a constant $c(p)>0$ (independent of $n$) such that
\begin{align*}
\norm{\Pn \bv}_{\Sob^{1,p}(\Omega)}  \leq c
\norm{ \bv}_{\Sob^{1,p}(\Omega)}  \quad \text{ for all } \bv \in \Sob^{1,p}_0(\Omega)^d \; \text{ and all } n \in \mathbb{N}.
\end{align*}
\end{enumerate}
}
\end{Assumption}

\begin{Assumption}[Projector $\Pn_{\mathbb{Q}}$]\label{assump-PQ}\hfill\\
Assume that for each $n \in \mathbb{N}$ there exists a linear projector $\Pn_{\mathbb{Q}} \colon \Leb^{1}(\Omega) \to \Qn$ such that, for any $p \in (1, \infty)$, there exists a constant $c(p)>0$ such that
\begin{align}\label{Pq-stab}
\norm{\Pn_{\mathbb{Q}} h}_{\Leb^p(\Omega)} \leq c \norm{h}_{\Leb^p(\Omega)} \quad \text{ for all } h \in \Leb^{p}(\Omega) \text{ and all } n \in \mathbb{N}.
\end{align}
\end{Assumption}

\begin{Remark}[Properties of $\Pn$ and $\Pn_{\mathbb{Q}}$] \label{Rem:proj} \hfill
\begin{enumerate}[label=(\roman*)]
\item \label{itm:approx-Pin} The global stability in Assumption~\ref{assump-Pin} \ref{itm:Pin-stab} and the approximability in \eqref{Vn-approx} yield that
\begin{align*}
\norm{\bv - \Pn \bv}_{\Sob^{1, p}(\Omega)} \to 0, \quad \text{ as } n \to \infty
\end{align*}
for all $\bv \in \Sob^{1,p}_0(\Omega)^d$ with $p \in [1, \infty)$.
\item \label{itm:approx-PQ} The stability in
\eqref{Pq-stab} and the approximability in \eqref{Qn-approx} imply that
\begin{align*}
\norm{h - \Pn_{\mathbb{Q}} h}_{\Leb^p(\Omega)} \to 0, \quad \text{ as } n \to \infty
\end{align*}
for all $h \in \Leb^p(\Omega)$ with $p \in [1, \infty)$.
\item \label{itm:inf-sup} The existence of the Bogovski\u{\i} operator, see \cite{B.1979} and also \cite[p.~223]{DMS.2008} implies that the continuous inf-sup condition holds for any $p \in (1, \infty)$. With this and the Assumption~\ref{assump-Pin} the corresponding discrete inf-sup condition holds uniformly in $n \in \mathbb{N}$; see Fortin's Lemma for Banach spaces in \cite[Lem.~4.19]{EG.2004}. This means that the framework results in an inf-sup stable pair $(\Vn, \Qn)$.
\end{enumerate}
\end{Remark}

\begin{Example}[Finite Element Spaces]\hfill\\
The following elements satisfy Assumptions~\ref{assump-approx}--\ref{assump-PQ}:
\begin{enumerate}[label={(\roman*)}]
\item the $\mathbb{P}_2-\mathbb{P}_0$ element for $d=2$, see \cite[Sec.~8.4.3]{BBF.2013}, where the projector $\Pn$ is given and Assumption~\ref{assump-Pin} \ref{itm:Pin-div} is shown; The stability in \ref{itm:Pin-stab} can be proved similarly as for the MINI element, see \cite[App. A.1]{BBDR.2012} and \cite[pp.~990]{DKS.2013}.
\item the conforming Crouzeix--Raviart element, for $d=2$, see \cite[Ex.~8.6.1]{BBF.2013} and \cite{CR.1973}; the projector $\Pn$ satisfying Assumption~\ref{assump-Pin} \ref{itm:Pin-div} is given in \cite[pp.~49]{CR.1973} and it can be shown to satisfy \ref{itm:Pin-stab}, see for example \cite[Thm.~3.3]{GS.2003}.
\item the MINI element for $d \in \{2,3\}$ ($k=1$), see \cite[Sec.~8.4.2, 8.7.1]{BBF.2013}; the proof that Assumption~\ref{assump-Pin} is satisfied is given in \cite[App.~A.1]{BBDR.2012}, see also \cite[Lem.~4.5]{GL.2001b} and \cite[pp.~990]{DKS.2013}.
\item the Taylor--Hood element and its generalizations for $d \in \{2,3\}$ and $k \geq d$, see \cite[Sec.~8.8.2]{BBF.2013}; the proof of Assumption~\ref{assump-Pin} can be found in \cite[Thm.~3.1, 3.2]{GS.2003}.
\end{enumerate}
The following element satisfies Assumption~\ref{assump-approx}--\ref{assump-PQ} and additionally, that $\Vd \subset \Sob^{1,\infty}_{0,\diver}(\Omega)^d$:
\begin{enumerate}[label=(\roman*),start=5]
\item the high-dimensional family of Guzm\'an--Neilan elements for $k \geq 1$, $d \in \{2,3\}$, see \cite{GN.2014, GN.2014b}. Note that the stability in Assumption~\ref{assump-Pin} \ref{itm:Pin-stab} is only shown for $p=2$, if $d=2$. 
\end{enumerate}
\end{Example}

\subsubsection*{$\Leb^2$-Projector to $\Vd$}

 Let us introduce the projector onto $\Vd$, given by
\begin{align} \label{def-Pn} \begin{split}
P^n_{\diver}\colon \Leb^2(\Omega)^d &\to \Vd, \quad  \text{ and for }  \bv\in \Leb^2(\Omega)^d,\\
\ll P^n_{\diver} \bv, \bV\right\rangle_{\Omega}&=(\bv, \bV)_{\Omega}\quad \text{ for all } \bV \in \Vd.
\end{split}
\end{align}
Directly from the definition we have $\Leb^2$-stability and optimality of the approximation in $\Leb^2(\Omega)^d$, i.e., for $\bv \in \Leb^2(\Omega)^d$ we have
\begin{align}
\norm{P^n_{\diver}\bv}_{\Leb^2(\Omega)} &\leq \norm{\bv}_{\Leb^2(\Omega)}, \label{est:L2-stab} \\
\label{opt-L2}
\norm{\bv-P^n_{\diver} \bv}_{\Leb^2(\Omega)} &\leq \norm{\bv -\bV}_{\Leb^2(\Omega)} \quad  \text{ for all } \bV \in \Vd.
\end{align}
By these properties and an approximation argument using the properties of $\Pn$ (see Remark~\ref{Rem:proj} \ref{itm:approx-Pin}) one can show that
\begin{align}\label{conv-Pnd}
P^n_{\diver} \bw \to \bw, \quad \text{ strongly in } \Leb^2(\Omega)^d, \text{ as } n \to \infty,
\end{align}
for any $\bw \in \l2d$.

\subsection{Convective Term and its Numerical Approximation \nopunct}\hfill\\
Motivated by the form of the convective term in the conservation of momentum equation, we consider the trilinear form $b$ defined by
\begin{align}\label{def:b}
b(\bu, \bv, \bw) \coloneqq - \ll \bu \otimes \bv, \nabla \bw \right\rangle_{\Omega} =
\ll \bu \otimes \bw, \nabla \bv \right\rangle_{\Omega} - \ll \diver \bu, \bv \cdot \bw \right\rangle_{\Omega},
\end{align}
for $\bu, \bv, \bw \in C^{\infty}_0(\Omega)^d$, where the second equality follows by integration by parts.
Hence for divergence-free functions $\bu$ the last term vanishes and $b(\bu, \cdot, \cdot)$ is skew-symmetric, i.e., $b(\bu, \bv, \bv)=0$ for $\bu \in C^{\infty}_{0, \diver}(\Omega)^d$ and $\bv \in C^{\infty}_0(\Omega)^d$.
This property can be extended to Sobolev functions.

As in general $\Vd \not\subset \Sob^{1,\infty}_{0, \diver}(\Omega)^d$, the second term in \eqref{def:b} need not vanish. To preserve the skew-symmetry of the trilinear form associated with the convective term the usual approach in the numerical analysis literature (see, e.g., \cite{T.1984}) is therefore to consider instead the skew-symmetric trilinear form
\begin{align}\begin{split}\label{def:wb}
\widetilde{b}(\bu, \bv, \bw) &\coloneqq \tfrac{1}{2}\left( \ll\bu \otimes \bw, \nabla \bv \right\rangle_{\Omega} - \ll \bu \otimes \bv, \nabla \bw \right\rangle_{\Omega}  \right)
= - \ll \bu \otimes \bv, \nabla \bw \right\rangle_{\Omega}  +\tfrac{1}{2} \ll \diver \bu, \bv \cdot \bw \right\rangle_{\Omega},   \end{split}
\end{align}
for $\bu, \bv, \bw \in C^{\infty}_0(\Omega)^d$.  Thus we have that $\widetilde{b}(\bu, \bv, \bv)=0$ regardless of the solenoidality of $\bu$.
Note that we have $b(\bu, \cdot, \cdot)=\widetilde{b}(\bu, \cdot, \cdot)$ for divergence-free functions $\bu$.

In the equations the terms appear in the form $b(\bu, \bu, \bv)$ and $\widetilde{b}(\bu, \bu, \bv)$, for the velocity $\bu$ and a test function $\bv$.
The natural function space for weak solutions of the problem is $\Leb^{\infty}(0,T;\Leb^2(\Omega)^d)\cap \Leb^{q}(0,T;\Sob^{1,q}_0(\Omega)^d)$, which embeds by Lemma~\ref{par-interpol} continuously into $\Leb^{\frac{q(d+2)}{d}}(Q)^d$. Also, provided that $q\geq \frac{2d}{d+2}$, we have that the embedding $\Leb^{\frac{q(d+2)}{d}}(Q)^d \ctemb \Leb^2(Q)^d$ is continuous, which means that the expression $b(\bu(t,\cdot),\bu(t,\cdot),\bv)$ is integrable on $(0,T)$, for any $\bv \in \Sob^{1,\infty}(\Omega)^d$. More specifically we can show that
for $\wq$ as defined in \eqref{def:hq} we have
\begin{align}\label{est:b-3}
\abs{\ll \bu(t,\cdot) \otimes \bu(t,\cdot), \nabla \bv \rr_{\Omega}}
\leq
 c \norm{\bu(t,\cdot)}_{\Leb^{\frac{q(d+2)}{d}}(\Omega)}^2 \norm{\bv}_{\Sob^{1,\wq}(\Omega)}, \quad \text{ provided that } q \geq \tfrac{2d}{d+2}.
\end{align}

On the other hand, for the modification (cf. the first term in \eqref{def:wb}) of the trilinear form $b$ associated with the convective term one obtains
\begin{align}\label{est:wb-3}
\abs{\ll \bu(t,\cdot) \otimes \bv, \nabla \bu(t,\cdot) \rr_{\Omega}}
\leq
 c \norm{\bu(t,\cdot)}_{\Leb^{\frac{q(d+2)}{d}}(\Omega)} \norm{ \bv}_{\Sob^{1,\wq}(\Omega)} \norm{\nabla \bu(t,\cdot)}_{\Leb^{q}(\Omega)}, \;\; \text{ if } q\geq\tfrac{2(d+1)}{d+2}.
\end{align}
Evidently, the source of this more restrictive requirement on $q$ is the modification of the trilinear form $b$, introduced in order to reinstate the skew symmetry of $b$, lost in the course of approximating the pointwise divergence-free
solution by discretely divergence-free finite element functions. We note in passing that the restriction $q\geq\tfrac{2(d+1)}{d+2}$ in the unsteady case corresponds to the restriction $q \geq \frac{2d}{d+1}$ in the steady case in \cite{DKS.2013}. 

To motivate the choice of the form of the additional regularization term that we shall add to the weak form to relax the excessive restriction $q\geq\tfrac{2(d+1)}{d+2}$  to the natural restriction on $q > \frac{2d}{d+2}$, we note that by H\"older's inequality we have that
\begin{align}\begin{split}\label{est:wb-r}
\norm{\widetilde{b}(\bu,\bu,\bv) }_{\Leb^1((0,T))}
&\leq
\norm{\bu}_{\Leb^{2q'}(Q)}^2 \norm{\nabla \bv}_{\Leb^{q}(Q)} + \norm{\bu}_{\Leb^{2q'}(Q)} \norm{\bv}_{\Leb^{2q'}(Q)} \norm{\nabla \bu}_{\Leb^{q}(Q)}\\
 &\leq c \norm{\bu}_{\Xq}^2 \norm{\bv}_{\Xq},
 \end{split}
\end{align}
for $\bu,\bv \in \Xq:=\Leb^{q}(0,T;\Sob^{1,q}_0(\Omega)^d)\cap \Leb^{2q'}(Q)^d$, see \eqref{def:Xq-u}, without any restrictions on the range of $q$, other than $q \in (1,\infty)$. This justifies the use of a regularizing term guaranteeing additional $\Leb^{2q'}$-integrability.

 \subsection{Time Discretization \nopunct}\label{FEM:time} \hfill \\
For the purpose of time discretization, let $l \in \mathbb{N}$ and define the time step by $\delta_l = {T}/{l}\to 0$, as $l \to \infty$. For $l \in \mathbb{N}$, we shall use the equidistant temporal grid on $[0,T]$ defined by $\{t^l_i\}_{i \in \{0,\ldots, l\}}$, where $t^l_i \coloneqq i \delta_l$, for $i \in \{ 0, \ldots, l\}$. In the following we will suppress the superscript $l$ and write $t_i$, $i \in \{0, \ldots, l\}$.

For a Banach space $X$ of functions, $l \in \mathbb{N}$ and a sequence $\{\varphi_{i}\}_{i \in \{1, \ldots, l\}}\subset X$ we consider the temporal difference quotient
\begin{align}\label{def:ap-t}
d_t \varphi_i \coloneqq \frac{1}{\delta_l}\left(\varphi_i- \varphi_{i-1}\right).
\end{align}
Furthermore, for $l \in \mathbb{N}$ we denote by $\mathbb{P}^l_0(0,T;X)$ the linear space of left-continuous piecewise constant mappings from $(0,T]$ into $X$ with respect to the equidistant temporal grid $\{t^l_0,\ldots, t_l^l \}\subset[0,T]$,
and by $\mathbb{P}^l_1(0,T;X)$ the space of continuous, piecewise affine functions from $[0,T]$ into $X$ with respect to the
same temporal grid. Let the piecewise constant and the piecewise affine interpolants $\overline{\varphi}$ and $\widetilde{\varphi}$ of $\{\varphi_i\}_{i \in \{0, \ldots, l\}}$ be defined by
\begin{alignat}{2}\label{def:interpol-c}
\overline{\varphi}(t)&\coloneqq\varphi_i,\qquad && \mbox{for $t \in (t_{i-1},  t_i]$, $i\in \{1,\ldots,l\}$,}\\ \label{def:interpol-l}
\widetilde{\varphi}(t)&\coloneqq\varphi_i\frac{t-t_{i-1}}{\delta_l} +\varphi_{i-1}\frac{t_i-t}{\delta_l},
\qquad && \mbox{for $t \in [t_{i-1},  t_i]$, $i\in \{1,\ldots,l\}$,}
\end{alignat}
so that $\overline{\varphi}, \widetilde{\varphi}, \partial_t \widetilde{\varphi} \in \Leb^{\infty}(0,T;X)$. Choosing the representative $\partial_t \widetilde{\varphi} \in \mathbb{P}^l_0(0,T;X)$, for $t \in (t_{i-1}, t_i]$ we have $\partial_t \widetilde{\varphi}(t) = d_t \varphi_i$
and
\begin{align}\label{ip-id-1}
\overline{\varphi}(t)- \widetilde{\varphi}(t) = (t_i-t)\partial_t \widetilde{\varphi}(t).
\end{align}
Furthermore, note that one has
\begin{align}\label{est:ip-est-1}
\norm{\overline{\varphi}}_{\Leb^{\infty}(0,T;X)}
&=  \max_{i \in \{1, \ldots, l\}} \norm{\varphi_i}_X, \quad&
\norm{\overline{\varphi}}_{\Leb^{p}(0,T;X)}^p
&=\delta_l \sum_{i=1}^l \norm{\varphi_i}_X^p,\quad \mbox{for $p \in [1, \infty)$}, \\
\label{est:ip-est-2}
\norm{\widetilde{\varphi}}_{\Leb^{\infty}(0,T;X)}
&= \max_{i \in \{0, \ldots, l\}} \norm{\varphi_i}_X,
\quad&
\norm{\widetilde{\varphi}}_{\Leb^{p}(0,T;X)}^p
&\leq c(p)\delta_l
 \sum_{i=0}^l \norm{\varphi_i}_X^p,\quad \mbox{for $p \in [1, \infty)$},
\end{align}
where $0<c(p) \leq 1$ by the Riesz--Thorin interpolation theorem (cf. \cite[Thm. 1.1.1, p.2]{BL.1976}). 

For a Bochner function $\psi \in \Leb^{p}(0,T;X)$, $p \in [1, \infty)$, we define the time averages with respect to the time grid $\{t_0, \ldots, t_l\}$, for $l \in \mathbb{N}$, by
\begin{align}\label{def:time-av}
\psi_i \coloneqq \fint_{t_{i-1}}^{t_i} \psi(t, \cdot)\, \d t \in X, \quad i \in \{1, \ldots, l\}.
\end{align}
Considering the piecewise constant interpolant $\overline{\psi}$ of the set of values $\{\psi_i\}_{i \in \{1, \ldots, l\}}$,
with $\psi_i$ defined by \eqref{def:time-av}, it follows by Jensen's and by H\"older's inequality that
\begin{align}
\label{t-av-stab}
\norm{\overline{\psi}}_{\Leb^{p}(0,T;X)}& \leq \norm{\psi}_{\Leb^{p}(0,T;X)} \quad \text{ for all } p \in [1,\infty],
\end{align}
and, for any $p \in [1,\infty)$,
\begin{align}
\label{t-av-conv}
\overline{\psi} \to \psi \quad &\text{ strongly in } \Leb^{p}(0,T;X), \quad \text{ as } l \to \infty,
\end{align}
thanks to the inequality $\|\psi - \overline{\psi}\|_{\Leb^p(0,T;X)} \leq T^{\frac{1}{p}} \delta_l \|\psi\|_{C^{0,1}([0,T];X)}$
for all $\psi \in C^{0,1}([0,T];X)$ and $p \in [1,\infty]$, the density of $C^{0,1}([0,T];X)$ in $\Leb^p(0,T;X)$ for $p \in [1,\infty)$, which, together with \eqref{t-av-stab} implies \eqref{t-av-conv} for $p \in [1,\infty)$.

To simplify the notation we will denote $Q_s^t \coloneqq (s,t)\times \Omega,$ for $0\leq s<t\leq T$, and $Q_s \coloneqq Q_0^s$, for $s \in (0,T]$. Furthermore, let us introduce the notation $Q_{i-1}^i \coloneqq Q_{t_{i-1}}^{t_i}$ and $Q_i \coloneqq Q_{t_i}$, for $i \in \{1, \ldots, l\}$.

\section{Convergence Proof}\label{sec:4}

\noindent Motivated by the approach in \cite[Sec.~3.1]{BGMS.2012} we consider the following levels of approximation.
\begin{enumerate}
\item[$k \in \mathbb{N}$:] The selection $\BS^{\star}$ given in Lemma~\ref{sel} is approximated by a family of Carath\'{e}odory functions $\{\BS^k\}_{k \in \mathbb{N}}$, which satisfy Assumption~\ref{Assump-Sk}. The approximation of the stress is then explicit and continuous in $\BD \bu$.
\item[$l \in \mathbb{N}$:] A time stepping based on the implicit Euler method is introduced similarly as, e.g., in \cite{CHP.2010, T.1984}, see Subsection~\ref{FEM:time}.
\item[$n \in \mathbb{N}$:] The velocity $\bu$ is approximated by a Galerkin approximation in finite element spaces in the spatial variable, see Section~\ref{Sec:FEM}.
\item[$m \in \mathbb{N}$:] The regularizing term $\frac{1}{m}\abs{\bu}^{2q'-2}\bu$ is added to the equation to gain admissibility of the approximate solutions if $q \leq \frac{3d+2}{d+2}$ and to enable us to use the bound on $\widetilde{b}(\cdot, \cdot, \cdot)$ in \eqref{est:wb-r}, without imposing the restriction $q > \frac{3d+2}{d+2}$.
\end{enumerate}
This results in a fully discrete approximation.
The limits are taken in the order $k \to \infty$, $l,n \to \infty$, and then $m \to \infty$, and we can take the limits in $l,n \to \infty$ simultaneously.
To simplify the notation we shall write
\begin{align}\label{not:3conv}
\bv^{k,l,n,m} \underset{k}{\to} \underset{(l,n)}{\to} \underset{m}{\to} \bv \quad \text{ in } X, \quad \text{ as } k \to \infty,\, l,n \to \infty,\, m \to \infty,
\end{align}
to denote the fact that the limits $k,(l,n),m$ are taken successively in the order of indexing (from left-to right) and the space $X$ describes the weakest topology of the three limits. We will use analogous notation for weak and weak* convergence. In each step one has to identify the equation and the implicit relation, which is the most challenging part.
The most significant difference compared to \cite{BGMS.2012} lies in the passage to the limits $l, n \to \infty$ and the identification of the implicit law.

As both the external force $\bf$ and the approximate stress $\BS^k$ will be allowed to be time-dependent, and the
time-dependence will not be assumed to be continuous, we shall consider an integral-averaged version in the approximate problem. Let us recall the notation in Subsection~\ref{FEM:time} and introduce for $\bf \in \Leb^{q'}(0,T;\Sob^{-1,q'}(\Omega)^d)$ and $\BS^k: Q \times \Rds \to \Rds$ as in Assumption~\ref{Assump-Sk} and $l \in \mathbb{N}$ the averages with respect to the time grid $\{t_i\}_{i \in \{0, \ldots, l\}}$ defined, for $i\in \{1,\ldots,l\}$, by
\begin{align} \label{def:f-i}
\bf_i(\bx) &\coloneqq \fint_{t_{i-1}}^{t_i} \bf(t, \bx)\, \d t, 
\quad \bx \in \Omega,\\
\label{def:Sk-i}
\BS^{k}_i(\bx, \BB) &\coloneqq \fint_{t_{i-1}}^{t_i} \BS^k(t, \bx,\BB)\, \d t, \quad \bx \in \Omega, \,\BB \in \Rds,
\end{align}
and let the corresponding piecewise constant interpolants $\overline{\bf}$ and $\overline{\BS}^k$ be defined as in Subsection~\ref{FEM:time}  \eqref{def:interpol-c}.
Recall that by \eqref{t-av-stab} and \eqref{t-av-conv} we have that
\begin{align}
\label{of-bd}
\norm{\overline{\bf}}_{\Leb^{q'}(0,T;\Sob^{-1,q'}(\Omega))}& \leq \norm{\bf}_{\Leb^{q'}(0,T;\Sob^{-1,q'}(\Omega))} \quad \text{ for all } l \in \mathbb{N},\\
\label{of-conv}
\overline{\bf} &\to \bf \quad \text{ strongly in } \Leb^{q'}(0,T;\Sob^{-1,q'}(\Omega)^d), \quad \text{ as } l \to \infty.
\end{align}

For $\bu, \bv \in \Vd$ we introduce
\begin{align}\begin{split}\label{Lklnm-i}
\mathfrak{L}^{k,l,n,m}_i[\bu; \bv]\coloneqq&
 - \widetilde{b}(\bu,\bu,\bv)
- \ll \BS^k_i(\cdot,\BD\bu), \BD\bv\rr_{\Omega}
- \frac{1}{m}\ll \abs{\bu}^{2q'-2} \bu, \bv \rr_{\Omega}
+ \ll \bf_i, \bv\rr_{\Omega},
\end{split}
\end{align}
for $k,l,n, m \in \mathbb{N}$, $i \in \{1, \ldots, l\}$ and $\widetilde{b}(\cdot,\cdot,\cdot)$ as defined in \eqref{def:wb}.

\subsubsection*{Approximate Problem:} For $k,l,n, m \in \mathbb{N}$ find a sequence $\{\bU^{k,l,n,m}_{i}\}_{i\in \{0, \ldots, l\}}\subset \Vd$ such that
\begin{align}\label{ap-k-i-1}
\bU^{k,l,n,m}_0=P^n_{\diver} \bu_0,
\end{align}
and for a given $\bU^{k,l,n,m}_{i-1}\in \Vd$ the approximate solution on the next time level, $\bU^{k,l,n, m}_i\in \Vd$, is defined, for $i\in \{1,\ldots,l\}$, by
\begin{alignat}{2}\label{ap-k-i-2}
\ll d_t \bU^{k,l,n,m}_i, \bW\rr_{\Omega}
&= \mathfrak{L}^{k,l,n,m}_i[\bU^{k,l,n,m}_i;\bW]  \qquad &&\text{ for all } \bW \in \Vd,
\end{alignat}
 where $P^n_{\diver}$ is the $\Leb^2$-projector onto $\Vd$, defined in \eqref{def-Pn}.
\smallskip

For each $i \in \{1, \ldots, l\}$ a fully implicit problem has to be solved, since the numerical solution from the previous time level only appears in the term involving $d_t \bU^{k,l,n,m}_i$, as defined in \eqref{def:ap-t}.\smallskip

\begin{Theorem}[Main Result]\label{Thm:main}
In addition to the assumptions of Definition~\ref{def:w-sol} let $\BS^k$ satisfy Assumption~\ref{Assump-Sk}.
 For the finite element approximation let Assumption~\ref{assump-triang} on the domain and on the family of simplicial partitions be satisfied. Let $\Vn$ and let $\Vd$ be as introduced in \eqref{Vn} and \eqref{Vd}, respectively, and assume that Assumptions~\ref{assump-approx}, \ref{assump-Pin} and \ref{assump-PQ} hold. Then, for all $k,l,n,m \in \mathbb{N}$ there exists a sequence $\{\bU^{k,l,n,m}_i\}_{i \in \{0, \ldots, l\}} \subset \Vd$ solving \eqref{ap-k-i-1}, \eqref{ap-k-i-2}.
Moreover, if $q \in \left( \frac{2d}{d+2}, \infty \right)$, then there exists a weak solution $(\bu, \BS)$ of \ref{def:P} according to Definition~\ref{def:w-sol} and for the piecewise constant interpolant $\oU^{k,l,n,m}$ and the continuous, piecewise affine interpolant  $\wU^{k,l,n,m}$ of $\{\bU^{k,l,n,m}_i\}_{i \in \{0, \ldots, l\}}$ as defined in \eqref{def:interpol-c} and \eqref{def:interpol-l}, respectively, (up to not relabelled subsequences) one has that
\begin{alignat*}{2}
 \oU^{k,l,n,m},  \wU^{k,l,n,m} &
 \underset{k}{\to}
 \underset{(l,n)}{\to}
 \underset{m}{\to} \bu \quad && \text{ strongly in  } \Leb^q(0,T;\Leb^{2}(\Omega)^d),\\ 
     \oU^{k,l,n,m},  \wU^{k,l,n,m}&
\overset{*}{\underset{k}{\rightharpoonup}}
\overset{*}{\underset{(l,n)}{\rightharpoonup}}
\overset{*}{\underset{m}{\rightharpoonup}}
\bu \quad && \text{ weakly in  } \Leb^{\infty}(0,T;\Leb^2(\Omega)^d),\\
  \oU^{k,l,n,m} &
  \underset{k}{\rightharpoonup}
  \underset{(l,n)}{\rightharpoonup}
  \underset{m}{\rightharpoonup}
  \bu \quad && \text{ weakly in  } \Leb^q(0,T;\Sob^{1,q}_{0}(\Omega)^d),\\
  \BS^k(\cdot, \cdot, \BD \oU^{k,l,n,m}) &
  \underset{k}{\rightharpoonup}
  \underset{(l,n)}{\rightharpoonup}
  \underset{m}{\rightharpoonup} \BS \quad && \text{ weakly in } \Leb^{q'}(Q)^{d\times d},
\end{alignat*}
as $k \to \infty$, $(l,n) \to \infty$ (combined) and $m \to \infty$, when taking the limits successively, without restrictions on the relation between the discretization parameters $\delta_l$ and $h_n$.
\end{Theorem}

\begin{Remark}\hfill
\begin{enumerate}[label=(\roman*)]
\item 
In the proof of Theorem~\ref{Thm:main} it is essential that the limits are taken in the indicated order.
\item If $\BS^{\star}$ is a Carath\'{e}odory function, then the approximation level corresponding to $k \in \mathbb{N}$ can be skipped.
\item 
For Lipschitz polytopal domains Theorem~\ref{Thm:main} is also a new existence result, since in \cite{BGMS.2012} a Navier slip boundary condition and $\partial \Omega \in C^{1,1}$ are assumed.
\item 
The convergence proof is presented for discretely divergence-free velocity functions. If $\Vd \subset \Sob^{1,\infty}_{0,\diver}(\Omega)^d$, then no modification of the convective term is required and the proof that $\bu^m$ is divergence-free is also simpler.
\end{enumerate}
\end{Remark}

The rest of this section consists of the proof of Theorem~\ref{Thm:main}, which relies on Lemmas~\ref{Lem:k-ex}--\ref{Lem:k-id} dealing with the existence of the discrete solution, and the limit $k \to \infty$, Lemmas~\ref{Lem:ln-conv} and \ref{Lem:ln-id} covering the combined limit $l,n\to \infty$, and Lemmas~\ref{Lem:m-conv} and \ref{Lem:m-id} the limit $m \to \infty$.

\subsection*{Limit $k \to \infty$}
The existence and convergence in Lemmas~\ref{Lem:k-ex} and \ref{Lem:k-conv} follow by a standard approach presented, e.g., in \cite{T.1984}, with minor modifications required to deal with the time-dependence of $\BS^k$.
Taking $k \to \infty$ we remain in the finite-dimensional setting and hence strong convergence of the sequence of symmetric gradients follows. Consequently, the identification of the limiting equation is based on the properties of the sequence $\{\BS^k\}_{k \in \mathbb{N}}$ according to Assumption~\ref{Assump-Sk}; c.f. \cite{BGMS.2012}.

\begin{Lemma}[Existence of Approximate Solutions]\label{Lem:k-ex}
For each $\kappa \coloneqq (k,l,n,m) \in \mathbb{N}^4$, there exists a sequence $\{\uk_i\}_{i \in \{0, \ldots, l\}} \subset \Vd$, which satisfies \eqref{ap-k-i-1},
\eqref{ap-k-i-2}. Furthermore, there exists a constant $c>0$ such that for all $\kappa=(k,l,n,m) \in \mathbb{N}^4$ one has that
  \begin{align}\begin{split}\label{ap-est:k}
\max_{j \in \{1, \ldots, l \}} \norm{\uk_j}_{\Leb^2(\Omega)}^2
&+ \sum_{j=1}^l \norm{\uk_j - \uk_{j-1}}_{\Leb^2(\Omega)}^2
+ \delta_l \sum_{j=1}^l \norm{\uk_j}_{\Sob^{1,q}(\Omega)}^q \\
&
+  \sum_{j=1}^l \norm{\BS^k(\cdot,\cdot,\BD \uk_j)}_{\Leb^{q'}(Q_{j-1}^j)}^{q'}
+  \frac{\delta_l}{m}\sum_{j=1}^l \norm{\uk_j}_{\Leb^{2q'}(\Omega)}^{2q'} \leq c.
\end{split}
 \end{align}
 \begin{proof}\hspace*{\fill}
\smallskip

\noindent
\textit{Step 1: A~priori estimates}. 
The a~priori estimates follow from standard arguments, see \cite{T.1984}, in combination with the estimates available for $\BS^k$ by Assumption~\ref{Assump-Sk}: testing \eqref{ap-k-i-2} with $\bW = \uk_i \in \Vd$ one obtains
\begin{align}\label{est:enid-k}
\ll d_t \uk_i, \uk_i \rr_{\Omega} + \ll \BS^k_i(\cdot, \BD \uk_i), \BD \uk_i \rr_{\Omega} + \frac{1}{m} \norm{\uk_i}_{\Leb^{2q'}(\Omega)}^{2q'} = \ll \bf_i, \uk_i \rr_{\Omega},
\end{align}
since the term involving $\widetilde{b}$ vanishes by skew-symmetry.
By the fact that $2a(a-b)=a^2-b^2+(a-b)^2$, for $a,b \in \mathbb{R}$ and by the definition of $\d_t \uk_i$ in \eqref{def:ap-t}, the first term in \eqref{est:enid-k} can be rewritten as
\begin{align}\begin{split}\label{ap-est:k-1}
\ll d_t \uk_i, \uk_i \rr_{\Omega}
&= \frac{1}{\delta_l} \ll \uk_i - \uk_{i-1}, \uk_i \rr_{\Omega} \\
&= \frac{1}{2\delta_l} \left(\norm{\uk_i}_{\Leb^2(\Omega)}^2 - \norm{\uk_{i-1}}_{\Leb^2(\Omega)}^2 + \norm{\uk_i-\uk_{i-1}}_{\Leb^2(\Omega)}^2 \right).
\end{split}
\end{align}
Using the definition of $\BS^k_i$ in \eqref{def:Sk-i} and Assumption~\ref{Assump-Sk} \ref{itm:al-3} one has that
\begin{align}
\ll \BS^k_i(\cdot, \BD \uk_i), \BD \uk_i \rr_{\Omega}
&\overset{\eqref{def:Sk-i}}{=}
\ll \fint_{t_{i-1}}^{t_i} \BS^k(t,\cdot,\BD \uk_i)\, \d t, \BD \uk_i \rr_{\Omega}
= \frac{1}{\delta_l} \ll \BS^k(\cdot,\cdot, \BD \uk_i), \BD \uk_i \rr_{Q_{i-1}^i} \notag \\
&\geq \frac{1}{\delta_l} \int_{Q^i_{i-1}} - \abs{\widetilde{g}(\cdot)} + \widetilde{c}_{*} \left( \abs{\BD \uk_i}^q + \abs{\BS^k(\cdot,\cdot,\BD \uk_i)}^{q'} \right)
 \d \bz  \label{ap-est:k-2} \\
 &\geq - \frac{1}{\delta_l} \norm{\widetilde{g}}_{\Leb^1(Q_{i-1}^i)}
 + c \norm{\uk_i}^q_{\Sob^{1,q}(\Omega)}
 + \frac{ \widetilde{c}_*}{\delta_l} \norm{\BS^k (\cdot,\cdot,\BD \uk_i)}_{\Leb^{q'}(Q_{i-1}^i)}^{q'},\notag
\end{align}
where the last inequality follows by Korn's and Poincar\'e's inequality.
On the term on the right-hand side of \eqref{est:enid-k} by duality of norms and by Young's inequality with $\varepsilon>0$ we obtain that
\begin{align}\begin{split}\label{ap-est:k-4}
\ll \bf_i, \uk_i \right\rangle_{\Omega}
&
\leq \norm{\bf_i}_{\Sob^{-1,q'}(\Omega)} \norm{\uk_i}_{\Sob^{1,q}(\Omega)}
\leq c(\varepsilon) \norm{\bf_i}_{\Sob^{-1,q'}(\Omega)}^{q'} +  \varepsilon \norm{ \uk_i}_{\Sob^{1,q}(\Omega)}^q\\
&\leq \frac{c(\varepsilon)}{\delta_l} \norm{\bf}_{\Leb^{q'}(t_{i-1}, t_i;\Sob^{-1,q'}(\Omega))}^{q'} +  \varepsilon \norm{ \uk_i}_{\Sob^{1,q}(\Omega)}^q, \end{split}
\end{align}
where the last inequality follows by \eqref{of-bd}.
Applying the estimates \eqref{ap-est:k-1}--\eqref{ap-est:k-4} in \eqref{est:enid-k}, after rearranging, choosing $\varepsilon>0$ small enough and multiplying by $\delta_l$, we arrive at
\begin{align}\begin{split}\label{ap-est:k-5}
\norm{\uk_i}_{\Leb^2(\Omega)}^2 &- \norm{\uk_{i-1}}_{\Leb^2(\Omega)}^2 + \norm{\uk_i-\uk_{i-1}}_{\Leb^2(\Omega)}^2 \\
& + \delta_l \norm{\uk_i}^q_{\Sob^{1,q}(\Omega)}
 + \norm{\BS^k (\cdot,\cdot,\BD \uk_i)}_{\Leb^{q'}(Q_{i-1}^i)}^{q'}
+ \frac{\delta_l}{m} \norm{\uk_i}_{\Leb^{2q'}(\Omega)}^{2q'}\\
&\quad  \leq c(\varepsilon)\norm{\bf}_{\Leb^{q'}(t_{i-1}, t_i;\Sob^{-1,q'}(\Omega))}^{q'}
 + c \norm{\widetilde{g}}_{\Leb^1(Q_{i-1}^i)}.
 \end{split}
\end{align}
 For arbitrary $j \in \{1, \ldots, l\}$, summing over $i \in \{1, \ldots, j\}$, yields
\begin{align}\begin{split}\label{ap-est:k-6}
\norm{\uk_j}_{\Leb^2(\Omega)}^2 &- \norm{\uk_{0}}_{\Leb^2(\Omega)}^2 + \sum_{i=1}^j\norm{\uk_i-\uk_{i-1}}_{\Leb^2(\Omega)}^2 \\
&+ \delta_l \sum_{i=1}^j \norm{\uk_i}^q_{\Sob^{1,q}(\Omega)}
+ \sum_{i=1}^j \norm{\BS^k(\cdot,\cdot,\BD \uk_i)}_{\Leb^{q'}(Q_{i-1}^i)}^{q'}
+ \frac{\delta_l}{m} \sum_{i=1}^j \norm{\uk_i}_{\Leb^{2q'}(\Omega)}^{2q'}\\
& \qquad\quad  \leq c (\norm{\bf}_{\Leb^{q'}(0, T;\Sob^{-1,q'}(\Omega))}^{q'}
 + \norm{\widetilde{g}}_{\Leb^1(Q)}),
 \end{split}
\end{align}
because of cancellation in the first term. Applying the estimate
\begin{align}\label{ap-est:k-7}
\norm{\uk_{0}}_{\Leb^2(\Omega)}^2 \overset{\eqref{ap-k-i-1}}{=} \norm{P^n_{\diver}\bu_0}_{\Leb^2(\Omega)}^2 \overset{\eqref{est:L2-stab}}{\leq}  \norm{\bu_0}_{\Leb^2(\Omega)}^2,
\end{align}
and taking the supremum over all $j \in \{1, \ldots, l\}$ in \eqref{ap-est:k-6} finishes the proof of \eqref{ap-est:k}.
\smallskip

\noindent 
\textit{Step 2: Existence of $\{\uk_i\}_{i \in \{0, \ldots, l\}}$}. Let us fix $\kappa \in \mathbb{N}^4$ and since $\uk_0 = P^n_{\diver} \bu_0$ by \eqref{ap-k-i-1}, we only have to show that for a given $\uk_{i-1} \in \Vd$, there exists a $\uk_{i} \in \Vd$ such that \eqref{ap-k-i-2} is satisfied. So let $i \in \{1, \ldots, l\}$ be fixed and let $\uk_{i-1}\in \Vd$ be given. We wish to find $\bU = \sum_{j=1}^{d_n} \alpha_j \bW_j \in \Vd$, where $\{\bW_1, \ldots, \bW_{d_n}\}$ is a basis of $\Vd$ and $\bal \in \mathbb{R}^{d_n}$ such that
\begin{align}\label{ex:k-1}
\frac{1}{\delta_l} \ll \bU - \uk_{i-1} , \bW \rr_{\Omega} = \mathfrak{L}^{\kappa}_i[\bU;\bW] \quad \text{ for all } \bW \in \Vd,
\end{align}
where $\mathfrak{L}^{\kappa}_i[\cdot,\cdot]$ is defined in \eqref{Lklnm-i}.
Defining $\bF:\mathbb{R}^{d_n} \to \mathbb{R}^{d_n}$ by
\begin{align}\label{def:F-t}
\bF(\bal) \coloneqq \left(-  \mathfrak{L}^{\kappa}_i\left[\sum_{j=1}^{d_n} \bal_j \bW_j; \bW_{r}\right] + \frac{1}{\delta_l} \ll \sum_{j=1}^{d_n} \bal_j \bW_j - \uk_{i-1}, \bW_r \rr_{\Omega} \right)_{r \in \{1, \ldots, d_n\}},
\end{align}
this amounts to finding $\bal \in \mathbb{R}^{d_n}$ such that $\bF(\bal)= \b0$.
Considering term by term one can see that $\bF$ is continuous in $\bal$.
By use of an $\Leb^2$-orthonormal basis of $\Vd$ one can show that
\begin{align}\label{ap-est:k-9-v}
\abs{\bal}^2 \leq c(n) \norm{\bU}^2_{\Leb^2(\Omega)}.
\end{align}
Furthermore, we find that
\begin{align}\begin{split}
\bF(\bal)\cdot \bal &= \widetilde{b}(\bU,\bU,\bU) + \ll \BS^k_i (\cdot, \BD \bU), \BD \bU \rr_{\Omega} + \frac{1}{m}\norm{\bU}_{\Leb^{2q'}(\Omega)}^{2q'} \\
&\quad - \ll \bf_i, \bU \rr_{\Omega} + \frac{1}{\delta_l}\ll \bU - \uk_{i-1}, \bU \rr_{\Omega}.
\end{split}
\end{align}
The first term vanishes thanks to the skew-symmetry and the other terms can be estimated similarly as in \eqref{ap-est:k-1}--\eqref{ap-est:k-5} to obtain
\begin{align*}
\bF(\bal)\cdot \bal &\geq -c (\widetilde{g}, \delta_l, \bf, \uk_i) + \frac{1}{2 \delta_l}\norm{\bU}^2_{\Leb^2(\Omega)} \overset{\eqref{ap-est:k-9-v}}{\geq} - c + c(n,l) \abs{\bal}^2.
\end{align*}
Thus, there exists an $R>0$ such that $\bF(\bal)\cdot \bal >0$, for any $\bal \in \mathbb{R}^{d_n}$ with $\abs{\bal}=R$. So $\bF$ is outward normal on $\partial B_R(\b0)$, and hence, as a consequence of Brouwer's fixed point theorem, $\bF$ has a zero in $B_R(\b0)$, see \cite[\S~5.7, (G.7), p.~104]{GD.2003}. This means that $\bU$ satisfying \eqref{ex:k-1} exists and the claim is proved.
 \end{proof}
\end{Lemma}

For $\kappa =(k,l,n,m) \in \mathbb{N}^4$ let the sequence of coefficients $\{\bal^{\kappa}_i\}_{i \in \{0,\ldots, l\}}\subset \mathbb{R}^{d_n}$ be such that $\uk_i=\sum_{j=1}^{d_n} (\alpha^{\kappa}_i)_j \bW_j$, where $\{\bW_1, \ldots, \bW_{d_n}\}$ is a basis of $\Vd$. Uniqueness is in general not guaranteed, so we choose one such sequence for each $\kappa \in \mathbb{N}^4$. 
Let $\overline{\bal}^{\kappa} \in \mathbb{P}^l_0(0,T;\mathbb{R}^{d_n}) \subset \Leb^{\infty}(0,T)^{d_n}$ and $\widetilde{\bal}^{\kappa}\in \mathbb{P}^l_1(0,T;\mathbb{R}^{d_n}) \subset \Sob^{1,\infty}(0,T)^{d_n}$
be the piecewise constant and piecewise affine interpolants as in
 \eqref{def:interpol-c} and \eqref{def:interpol-l}. We denote
\begin{align}\begin{split}\label{def:Uk-interpol}
\oU^{\kappa} (t,\bx) &\coloneqq \sum_{j=1}^{d_n} \overline{\alpha}^{\kappa}_j(t) \bW_j(\bx) \in \mathbb{P}^l_0(0,T;\Vd)
, \\
\wU^{\kappa} (t,\bx) &\coloneqq \sum_{j=1}^{d_n} \widetilde{\alpha}^{\kappa}_j(t) \bW_j(\bx) \in \mathbb{P}^l_1(0,T;\Vd) .
\end{split}
\end{align}
The so defined functions coincide with the respective interpolants of $\{\uk_i\}_{i \in \{0, \ldots, l\}}$,  defined in \eqref{def:interpol-c} and \eqref{def:interpol-l}.
Furthermore, for $t \in (0,T]$, $\bu \in \mathbb{P}^l_0(0,T;\Vd)$ and $\bv \in \Vd$ we introduce
\begin{align}\begin{split}\label{Lklnm}
\mathfrak{L}^{\kappa}[\bu; \bv](t)&\coloneqq
 - \widetilde{b}(\bu(t, \cdot),\bu(t, \cdot),\bv)
- \ll \overline{\BS}^k(t,\cdot,\BD\bu(t, \cdot)), \BD\bv\rr_{\Omega}\\
&\quad\,\, - \frac{1}{m}\ll \abs{\bu(t, \cdot)}^{2q'-2} \bu(t, \cdot), \bv \rr_{\Omega}
+ \ll \overline{\bf}(t,\cdot), \bv\rr_{\Omega},
\end{split}
\end{align}
for $\kappa=(k,l,n, m) \in \mathbb{N}^4$ and $\widetilde{b}(\cdot,\cdot,\cdot)$ as defined in
 \eqref{def:wb}.
Recall that $\overline{\bf}\in
\mathbb{P}^l_0(0,T;\Sob^{-1,q'}(\Omega)^d)
$ is the piecewise constant interpolant  of $\{\bf_i\}_{i \in \{1, \ldots, l\}}$, as defined in  \eqref{def:interpol-c} in Subsection~\ref{FEM:time}, and similarly,
$\overline{\BS}^k(t, \cdot,\cdot) = \BS^k_i(\cdot,\cdot)$, for $t \in (t_{i-1},t_i]$, which is piecewise constant with respect to the variable $t \in (0,T]$.

\begin{Lemma}[{Equation for $t \in (0,T]$ and Convergence $k \to \infty$}]\label{Lem:k-conv}
The functions $\oU^{\kappa}$ and $\wU^{\kappa}$ defined as in \eqref{def:Uk-interpol} satisfy
\begin{align} \label{ap-k-1}
\ll \partial_t \wU^{\kappa}(t, \cdot), \bW \rr_{\Omega}
&= \mathfrak{L}^{\kappa}[\oU^{\kappa};\bW](t)\quad
&&\text{ for  all } \bW \in \Vd, \text{ for all } t \in (0,T],\\ \label{ap-k-2}
\wU^{\kappa}(0,\cdot)
&= P^n_{\diver} \bu_0 \quad
&&\text{ in } \Omega,
\end{align}
for any $\kappa=(k,l,n,m) \in \mathbb{N}^4$.
For each $\L \coloneqq (l,n, m) \in \mathbb{N}^3$, i.e., $\kappa=(k,\L)$, there exists a sequence $\{\bU^{{\L}}_i\}_{i \in \{0, \ldots, l\}} \subset \Vd$, and subsequences such that the piecewise constant and continuous piecewise affine interpolants $\oU^{\L}\in \mathbb{P}^l_0(0,T;\Vd) $ and $\wU^{\L}\in \mathbb{P}^l_1(0,T;\Vd) $ of $\{\bU^{{\L}}_i\}_{i \in \{0, \ldots, l\}}$, as defined in \eqref{def:interpol-c} and \eqref{def:interpol-l}, satisfy
 \begin{alignat}{2}\label{ck-3}
 \oU^{k,{\L}}&\to \oU^{\L} \quad &&\text{ strongly in } \Leb^{\infty}(0,T;C(\Omega)^d),\\\label{ck-4}
\BD \oU^{k,{\L}}&\to \BD \oU^{\L} \quad &&\text{ strongly in } \Leb^{\infty}(0,T;\Leb^{\infty}(\Omega)^{d \times d}),\\\label{ck-5}
\nabla \oU^{k,{\L}}&\to \nabla \oU^{\L} \quad &&\text{ strongly in } \Leb^{\infty}(0,T;\Leb^{\infty}(\Omega)^{d\times d}),\\
 \label{ck-6v}
 \wU^{k,{\L}}&\to \wU^{\L} \quad &&\text{ strongly in } C([0,T];C(\Omega)^d),\\
\label{ck-6}
 \partial_t \wU^{k,{\L}} &\to  \partial_t \wU^{{\L}} \quad &&\text{ strongly in } \Leb^{\infty}(0,T;C(\Omega)^d),
 \end{alignat}
 as $k \to \infty$.  Furthermore, for each ${\L} \in \mathbb{N}^3$ there exist  $\BS^{\L} \in \Leb^{q'}(Q)^{d \times d}$ and $\overline{\BS}^{\L} \in \mathbb{P}^l_0(0,T;\Leb^{q'}(\Omega)^{d \times d})$ and subsequences such that
 \begin{alignat}{2}\label{ck-7}
 \BS^k(\cdot,\cdot, \BD \oU^{k,{\L}}) &\wconv \BS^{\L}  \quad &&\text{ weakly in } \Leb^{q'}(Q)^{d \times d}, \\\label{ck-8}
 \overline{\BS}^k(\cdot,\cdot, \BD \oU^{k,{\L}}) &\wconv \overline{\BS}^{\L}  \quad &&\text{ weakly in } \Leb^{q'}(Q)^{d \times d},
 \end{alignat}
 as $k \to \infty$, where, up to a representative, we have
 \begin{align}\label{SL-i}
 \overline{\BS}^{\L}(t,\cdot) = \BS^{\L}_i(\cdot) \coloneqq \fint_{t_{i-1}}^{t_i} \BS^{\L}(t, \cdot)\, \d t \quad \text{ for all } t \in (t_{i-1},t_i], \, i \in \{1, \ldots, l\}.
 \end{align}
 \begin{proof}\hspace*{\fill}
\smallskip

\noindent
\textit{Step 1: Identification of the equation}.
 \noindent We have that $\wU^{\kappa}(0,\cdot)=\uk_0 = P^n_{\diver} \bu_0$ 
by definition of $\wU^{\kappa}$ and
 by \eqref{ap-k-i-1}, which shows \eqref{ap-k-2}. The equation \eqref{ap-k-1} follows from
  \eqref{ap-k-i-2} and the fact that for $t \in (t_{i-1},t_i]$ we have that
 \begin{align*}
 \oU^{\kappa}(t,\cdot) = \uk_i,  \quad
 \partial_t \wU^{\kappa}(t,\cdot) = d_t \uk_i,  \quad
 \overline{\bf}(t,\cdot) = \bf_i  \;\;\text{ and }\;\;
 \overline{\BS}^k(t,\cdot,\cdot) = \BS^k_i(\cdot,\cdot),\qquad i \in \{1,\ldots,l\}.
 \end{align*}

\noindent
\textit{Step 2: Estimates}.
\noindent Let $\L=(l,n,m) \in \mathbb{N}^3$ be arbitrary, but fixed. When taking $k \to \infty$ we stay in the finite-dimensional setting, hence it suffices to focus on estimates for the coefficient functions $\overline{\bal}^{\kappa} = \overline{\bal}^{k,\L}$ and $\widetilde{\bal}^{\kappa} = \widetilde{\bal}^{k,\L}$, uniformly in $k \in \mathbb{N}$. Using the definition of $\oU^{k,\L}$ in  \eqref{def:Uk-interpol} with the estimate \eqref{ap-est:k-9-v}, the fact that the interpolants are piecewise constant and also the a~priori estimate in \eqref{ap-est:k}, we obtain
\begin{align*}
\abs{\overline{\bal}^{k,\L}(t)}^2
\overset{\eqref{ap-est:k-9-v}}{\leq} c(n) \norm{\oU^{k,\L}(t,\cdot)}_{\Leb^2(\Omega)}^2
\leq c(n) \max_{i \in \{1, \ldots, l \}}\norm{\bU^{k,\L}_i}_{\Leb^2(\Omega)}^2 \overset{\eqref{ap-est:k}}{\leq} c(n),
\end{align*}
for any $t \in (0,T)$ uniformly in $k \in \mathbb{N}$. 
The corresponding estimate holds for $|\widetilde{\bal}^{k,\L}(t)|^2$, using also that $\wU^{k,\L}(0,\cdot)=P^n_{\diver} \bu_0$ by \eqref{ap-k-2} and the $\Leb^2$-stability of $P^n_{\diver}$ in \eqref{est:L2-stab}. Thus, we have that
\begin{align}\label{est:k-1}
\norm{\overline{\bal}^{k,\L}}_{\Leb^{\infty}(0,T)} +  \norm{\widetilde{\bal}^{k,\L}}_{\Leb^{\infty}(0,T)} \leq c \quad \text{ for all } k \in \mathbb{N}.
\end{align}
Since the space of continuous piecewise affine functions $\mathbb{P}^l_1(0,T;\mathbb{R}^{d_n})\subset \Sob^{1,\infty}(0,T)^{d_n}$ with respect to the time grid $\{t_0, \ldots, t_l\}\subset [0,T]$ is finite-dimensional, all norms on it are equivalent, with the norm-equivalence constants depending on the (here fixed) dimension. So we also have that
\begin{align}\label{est:k-2}
\norm{\widetilde{\bal}^{k,\L}}_{\Sob^{1,\infty}(0,T)} \leq c(l,n) \quad \text{ for all } k \in \mathbb{N}.
\end{align}
From the a~priori estimate \eqref{ap-est:k} it follows directly that
\begin{align}\label{est:k-3}
\norm{\BS^k\left(\cdot,\cdot, \BD \oU^{k,\L}\right)}_{\Leb^{q'}(Q)}^{q'} \leq c \quad \text{ for all } k \in \mathbb{N}.
\end{align}
By the definition of $\BS^k_i$ in \eqref{def:Sk-i}, we have that
\begin{align}\begin{split}\label{est:k-4}
\norm{\overline{\BS}^k(\cdot,\cdot,\BD \oU^{k,\L})}_{\Leb^{q'}(Q)}^{q'}
&= \sum_{i=1}^l \norm{\BS^k_i(\cdot,\BD \bU
^{k,\L}_i)}_{\Leb^{q'}(Q_{i-1}^i)}^{q'} = \sum_{i=1}^l \delta_l \norm{\BS^k_i(\cdot,\BD \bU
^{k,\L}_i)}_{\Leb^{q'}(\Omega)}^{q'}\\
& \overset{\eqref{t-av-stab}}{\leq}
\sum_{i=1}^l
 \norm{\BS^k(\cdot,\cdot,\BD \bU
^{k,\L}_i)}_{\Leb^{q'}(Q_{i-1}^i)}^{q'} =  \norm{\BS^k(\cdot,\cdot,\BD \oU^{k,\L})}_{\Leb^{q'}(Q)}^{q'}\overset{\eqref{est:k-3}}{\leq} c
\end{split}
\end{align}
for all $k \in \mathbb{N}$. This also shows that
\begin{align}\label{est:k-5}
 \norm{\BS^k_i(\cdot,\BD \bU
^{k,\L}_i)}_{\Leb^{q'}(\Omega)}^q \leq \frac{c}{ \delta_l} \leq c(l),
\end{align}
for any $k \in \mathbb{N}$ and any $i \in \{1, \ldots, l\}$.

\smallskip

\noindent
\textit{Step 3: Convergence as $k \to \infty$.} Since $\{\overline{\bal}^{k,\L}\}_{k \in \mathbb{N}} \subset \mathbb{P}^l_0(0,T;\mathbb{R}^{d_n})$ and the space $\mathbb{P}^l_0(0,T;\mathbb{R}^{d_n})$ is finite-dimensional, \eqref{est:k-1} implies strong convergence of a subsequence, i.e., there exists an $\overline{\bal}^\L \in \mathbb{P}^l_0(0,T;\mathbb{R}^{d_n})$ such that
\begin{align}\label{conv-k-1}
\overline{\bal}^{k,\L} \to \overline{\bal}^{\L} \quad \text{ strongly in } \Leb^{\infty}(0,T)^{d_n}, \quad \text{ as } k \to \infty.
\end{align}
Similarly, we obtain from \eqref{est:k-2} that there exists a subsequence and an $\widetilde{\bal}^\L \in \mathbb{P}^l_1(0,T;\mathbb{R}^{d_n})$ such that
\begin{align}\label{conv-k-2}
\widetilde{\bal}^{k,\L} \to \widetilde{\bal}^{\L} \quad \text{ strongly in } \Sob^{1,\infty}(0,T)^{d_n}, \quad \text{ as } k \to \infty.
\end{align}
Note that the convergence holds pointwise everywhere in $(0,T]$, and hence,
\begin{align*}
\overline{\bal}^\L(t_{i}) \leftarrow \overline{\bal}^{k,\L}(t_{i})  = \bal^{k,\L}_i &= \widetilde{\bal}^{k,\L}(t_i) \rightarrow \widetilde{\bal}^\L(t_i), \quad \text{ as } k \to \infty,
\end{align*}
for any $i \in \{1, \ldots, l\}$, so the limits coincide and we can set $\bal^\L_i = \widetilde{\bal}^\L(t_i)=\overline{\bal}^{\L}(t_{i})$, for $i \in \{1,\ldots, l\}$. Since we also have that $\widetilde{\bal}^{k,\L}(0,\cdot)\to \widetilde{\bal}^{\L}$, we can set $\bal^{\L}_0 = \widetilde{\bal}^{\L}(0,\cdot)$. 
Then $\overline{\bal}^\L \in \mathbb{P}^l_0(0,T;\mathbb{R}^{d_n})$ and $\widetilde{\bal}^\L \in \mathbb{P}^l_1(0,T;\mathbb{R}^{d_n})$ are the piecewise constant and continuous piecewise affine interpolants of $\{\bal^\L_i\}_{i\in \{0, \ldots, l\}}$, respectively.
Let us set
 \begin{align}\label{def:UL}
 \oU^\L(t, \bx) \coloneqq \sum_{j = 1}^{d_n} \overline{\alpha}^\L_j(t) \bW_j(\bx),\quad
  \wU^\L(t, \bx) \coloneqq \sum_{j = 1}^{d_n} \widetilde{\alpha}^\L_j(t) \bW_j(\bx),
 \end{align}
 and let $\bU^\L_i = \sum_{j=1}^{d_n} (\alpha^\L_i)_j \bW_j \in \Vd$ for $i \in \{0, \ldots, l\}$. Note that by the above considerations concerning the coefficients we have that $\oU^\L$ and $ \wU^\L$ coincide with the respective interpolants of $\{\bU^{\L}_i\}_{i \in \{0, \ldots, l\}}$.
 By the convergence in \eqref{conv-k-1}, \eqref{conv-k-2} one obtains for the so-defined functions the convergence results \eqref{ck-3}--\eqref{ck-6}, as $k \to \infty$.
 By the (sequential) Banach--Alaoglu theorem, \eqref{est:k-3}--\eqref{est:k-5} imply that there exist $\BS^{\L}, \overline{\BS}^{\L} \in \Leb^{q'}(Q)^{d \times d}$ and $\BS^{\L}_i \in \Leb^{q'}(\Omega)^{d \times d}$ for $i \in \{1, \ldots, l\}$, and subsequences such that
 \begin{alignat}{2}\label{conv-k-7}
 \BS^k(\cdot,\cdot, \BD \oU^{k,{\L}}) &\wconv \BS^{\L}  \quad &&\text{ weakly in } \Leb^{q'}(Q)^{d \times d}, \\\label{conv-k-8}
 \overline{\BS}^k(\cdot,\cdot, \BD \oU^{k,{\L}}) &\wconv \overline{\BS}^{\L}  \quad &&\text{ weakly in } \Leb^{q'}(Q)^{d \times d},\\ \label{conv-k-9}
 \BS^{k}_i(\cdot, \BD \bU^{k,{\L}}_i) &\wconv \BS^{\L}_i \quad &&\text{ weakly in } \Leb^{q'}(\Omega)^{d \times d}, \;\; \text{ for } i \in \{1, \ldots, l\},
 \end{alignat}
 as $k \to \infty$.

It remains to show the identification of $\BS^{\L}$, $\overline{\BS}^{\L}$ and $\{\BS^{\L}_i\}_{i \in \{1, \ldots, l\}}$. Let $i \in \{1, \ldots, l\}$ be arbitrary, but fixed.
First let  $\varphi \in C^{\infty}_0((t_{i-1}, t_i))$ and $\Bv \in C^{\infty}_0(\Omega)^{d \times d}$. On the one hand, by \eqref{conv-k-8} we have
\begin{align}\label{conv-k-10}
\ll  \overline{\BS}^k(\cdot,\cdot, \BD \oU^{k,{\L}}), \varphi \Bv \rr_{Q_{i-1}^i} \to \ll  \overline{\BS}^{\L}, \varphi \Bv \rr_{Q_{i-1}^i}, \quad \text{ as } k \to \infty.
\end{align}
On the other hand by the definition of $\overline{\BS}^k(\cdot,\cdot,\BD \oU^{k,{\L}})$ as piecewise constant interpolant of the sequence $\{\BS^{k}_i(\cdot,\BD \bU^{k,\L}_i)\}_{i \in \{1, \ldots, l\}} $ and by \eqref{conv-k-9} we have
\begin{align}\begin{split}\label{conv-k-11}
\ll  \overline{\BS}^k(\cdot,\cdot, \BD \oU^{k,{\L}}), \varphi \Bv \rr_{Q_{i-1}^i}
&= \ll  \ll \BS^k_i(\cdot, \BD \bU^{k,{\L}}_i),\Bv \rr_{\Omega} \varphi \rr_{(t_{i-1},t_i)}\\
&= \ll 1, \varphi \rr_{(t_{i-1},t_i)} \ll \BS^k_i(\cdot, \BD \bU^{k,{\L}}_i),\Bv \rr_{\Omega} \\
&\to
\ll 1, \varphi \rr_{(t_{i-1},t_i)} \ll \BS^{\L}_i,\Bv \rr_{\Omega}
=
\ll \BS^{\L}_i,\Bv \varphi \rr_{Q_{i-1}^i},\quad \text{ as } k \to \infty.
\end{split}
\end{align}
Now, \eqref{conv-k-10} and  \eqref{conv-k-11} imply, by the uniqueness of the limit, that $\overline{\BS}^{\L}(t,\bx) = \BS^{\L}_i(\bx)$ for a.e. $(t,\bx) \in Q_{i-1}^i$, i.e., $\overline{\BS}^{\L}$ is piecewise constant in $t$ and we can choose the representative in $\mathbb{P}^l_0(0,T;\Leb^{q'}(\Omega)^{d \times d})$.
Again for $\Bv \in C^{\infty}_0(\Omega)^{d \times d}$ we have by \eqref{conv-k-9} that
\begin{align}\label{conv-k-12}
\ll \BS^{k}_i(\cdot, \BD \bU^{k,{\L}}_i), \Bv \rr_{\Omega} \to \ll \BS^{{\L}}_i, \Bv \rr_{\Omega}, \qquad \text{ as } k \to \infty.
\end{align}
On the other hand by the definition of $\BS^{k}_i(\cdot,\BD \bU^{k,{\L}}_i)$ in \eqref{def:Sk-i} and by \eqref{conv-k-7} we obtain that
\begin{align}\begin{split}\label{conv-k-13}
\ll \BS^{k}_i(\cdot, \BD \bU^{k,{\L}}_i), \Bv \rr_{\Omega}  &= \ll \fint_{t_{i-1}}^{t_i} \BS^k(t,\cdot,\BD \bU^{k,{\L}}_i) \, \d t , \Bv \rr_{\Omega} \\
&= \frac{1}{\delta_l} \ll \BS^k(\cdot,\cdot, \BD \oU^{k,{\L}}), \mathds{1}_{(t_{i-1}, t_i)} \Bv \rr_{Q} \\
&\to \frac{1}{\delta_l }\ll \BS^{\L}, \mathds{1}_{(t_{i-1}, t_i)} \Bv \rr_{Q} = \ll \fint_{t_{i-1}}^{t_i} \BS^{\L}(t, \cdot) \, \d t , \Bv \rr_{\Omega}, \qquad \text{ as } k \to \infty,
\end{split}
\end{align}
so by the uniqueness of limits, we conclude from  \eqref{conv-k-12}, \eqref{conv-k-13}, that $\BS^{\L}_i(\bx) = \fint_{t_{i-1}}^{t_i} \BS^{\L}(t, \bx) \, \d t$ for a.e. $\bx \in \Omega$, which completes the proof.
 \end{proof}
\end{Lemma}

For ${\L}=(l,n,m) \in \mathbb{N}^3$, $t \in (0,T]$,  $\bu \in \mathbb{P}^l_0(0,T;\Vd)$ and $\bv \in \Vd$, let us introduce
\begin{align}\begin{split}\label{Llnm-t}
\mathfrak{L}^{{\L}}[\bu; \bv](t)
\coloneqq&
- \widetilde{b}(\bu(t, \cdot),\bu(t, \cdot), \bv) - \ll \overline{\BS}^{{\L}}(t,\cdot), \BD\bv\rr_{\Omega}\\
&- \frac{1}{m}\ll \abs{\bu(t,\cdot)}^{2q'-2} \bu(t,\cdot), \bv \rr_{\Omega}
+ \ll \overline{\bf}(t, \cdot), \bv\rr_{\Omega},
\end{split}
\end{align}
where $\overline{\BS}^{\L} \in \mathbb{P}^l_0(0,T;\Leb^{q'}(\Omega)^{d \times d})$ is given in Lemma~\ref{Lem:k-conv}. 

Furthermore, for ${\L}=(l,n,m) \in \mathbb{N}^3$, $i \in \{1, \ldots, l\}$ and $\bu, \bv \in \Vd$ let us denote
\begin{align}\label{Llnm-i-t}
\mathfrak{L}^{{\L}}_i[\bu; \bv]
\coloneqq
- \widetilde{b}(\bu,\bu, \bv) - \ll \BS^{{\L}}_i, \BD\bv\rr_{\Omega}- \frac{1}{m}\ll \abs{\bu}^{2q'-2} \bu, \bv \rr_{\Omega}
+ \ll \bf_i, \bv\rr_{\Omega},
\end{align}
where $\BS^{\L}_i \in \Leb^{q'}(\Omega)^{d \times d}$ is given in Lemma~\ref{Lem:k-conv} \eqref{SL-i} and $\bf_i \in \Sob^{-1,q'}(\Omega)^d$ is defined in \eqref{def:f-i}.

\begin{Lemma}[Identification of the PDE as $k \to \infty$]\label{Lem:k-id}
The functions  $\oU^{\L}\in \mathbb{P}^l_0(0,T;\Vd) $, $\wU^{\L}\in \mathbb{P}^l_1(0,T;\Vd) $ and $\BS^{\L} \in \Leb^{q'}(Q)^{d \times d}$
given in Lemma~\ref{Lem:k-conv} satisfy
\begin{alignat}{2}
\ll \partial_t \wU^{{\L}}(t,\cdot), \bW\rr_{\Omega}
&= \mathfrak{L}^{{\L}}[\oU^{{\L}};\bW](t)  \qquad &&\text{ for all } \bW \in \Vd, \,\text{ for all } t \in (0,T], \label{ap-l-1}\\
\wU^{{\L}}(0,\cdot) &= P^n_{\diver} \bu_0(\cdot) \qquad &&\text{ in } \Omega,\label{ap-l-2}\\
(\BD \oU^{{\L}}(\bz), \BS^{{\L}}(\bz)) &\in \mathcal{A}(\bz)   \qquad &&\text{ for a.e. } \bz \in Q, \label{impl-l}
\end{alignat}
for all ${\L}=(l,n,m) \in \mathbb{N}^3$, where $\mathfrak{L}^{{\L}}[\cdot; \cdot](\cdot)$ is defined by \eqref{Llnm-t}, 
using $\overline{\BS}^{\L} \in \mathbb{P}^l_0(0,T;\Leb^{q'}(\Omega)^{d \times d})$ given by \eqref{SL-i} in Lemma~\ref{Lem:k-conv}.
Furthermore, the sequence $\{\bU^{{\L}}_i\}_{i \in \{0, \ldots, l\}} \subset \Vd$ given in Lemma~\ref{Lem:k-conv} satisfies
\begin{align}\label{ap-l-i-2}
\bU^{{\L}}_0&=P^n_{\diver} \bu_0,\\
 \label{ap-l-i-1}
\ll d_t \bU^{{\L}}_i, \bW\rr_{\Omega}
&= \mathfrak{L}^{{\L}}_i[\bU^{{\L}}_i;\bW]  \; &&\text{ for all } \bW \in \Vd, \quad \text{ for all } i \in \{1, \ldots, l\}.
\end{align}

\begin{proof}
Let ${\L}=(l,n,m) \in \mathbb{N}^3$ be arbitrary but fixed. \hspace*{\fill}
\smallskip

\noindent
\textit{Step 1: Identification of the initial condition}. By \eqref{ap-k-2} the family of continuous functions $\wU^{k,{\L}}$ satisfies the initial condition $\wU^{k,{\L}}(0, \cdot) = P^n_{\diver} \bu_0(\cdot)$ in $\Omega$ for all $k \in \mathbb{N}$. Thus, the sequence $\{\wU^{k,{\L}}(0, \cdot)\}_{k \in \mathbb{N}}$  is constant and the strong convergence in \eqref{ck-6v} implies that $P^n_{\diver} \bu_0(\cdot) = \wU^{k,{\L}}(0, \cdot) = \wU^{{\L}}(0, \cdot) = \bU^{\L}_0$ in $\Omega$, so \eqref{ap-l-2} and \eqref{ap-l-i-2} are satisfied.
\smallskip

\noindent
\textit{Step 2: Identification of the limiting equation}. Let $\bW \in \Vd$ and let $\varphi \in C^{\infty}_0((0,T))$ be arbitrary but fixed.
With the convergence of $\partial_t \wU^{k,{\L}}$ in  \eqref{ck-6} it follows that
\begin{align}\label{lk-1}
\ll \partial_t \wU^{k,{\L}}, \varphi \bW \rr_{Q} \to \ll \partial_t \wU^{{\L}}, \varphi \bW \rr_{Q}, \quad \text{ as } k \to \infty.
\end{align}
Further, by the strong convergence in \eqref{ck-3}, \eqref{ck-5} and \eqref{ck-8} it is straightforward to show that
\begin{align}\label{lk-2}
\ll \mathfrak{L}^{k,{\L}} [\oU^{k,{\L}}, \bW](\cdot),\varphi \rr_{(0,T)} \to \ll \mathfrak{L}^{{\L}} [\oU^{{\L}}, \bW](\cdot),\varphi \rr_{(0,T)}, \quad \text{ as } k \to \infty.
\end{align}
In particular, the strong convergence in \eqref{ck-3} and in \eqref{ck-5} allows us to take the limit in the numerical convective term without any restriction.
Finally, \eqref{lk-1} and \eqref{lk-2} applied in
 \eqref{ap-k-1} imply that \eqref{ap-l-1} holds for a.e. $t \in (0,T)$.
 Recall that $\overline{\BS}^{\L}(t,\cdot)= \BS^{\L}_i$, for any $t \in (t_{i-1},t_i]$ and $i \in \{1, \ldots, l\}$ by \eqref{SL-i}.
 Since now the terms in \eqref{ap-l-1} are constant on each interval $(t_{i-1},t_i]$, for $i \in \{1, \ldots, l\}$, the equation also holds for all $t \in (0,T]$.
  Also, \eqref{ap-l-i-1} follows from \eqref{ap-l-1} since $\partial_t \wU^{\L}(t,\cdot)= d_t \bU^{\L}_i$, $\oU^{\L}(t,\cdot)= \bU^{\L}_i$ and the corresponding holds for $\BS^{\L}_i$ and $\bf_i$, for any $t \in (t_{i-1},t_i]$ and $i \in \{1, \ldots, l\}$.
\smallskip

\noindent
\textit{Step 3: Identification of the implicit relation}. The proof of the implicit relation \eqref{impl-l} relies on the strong convergence of $\{\BD \oU^{k,\L}\}_{k \in \mathbb{N}}$ and the properties of $\BS^k$ stated in Assumption~\ref{Assump-Sk}. By the property \ref{itm:al-4} in Assumption~\ref{Assump-Sk} on $\BS^k$ and the boundedness of $\{\BD \oU^{k,{\L}}\}_{k \in \mathbb{N}}$ in $\Leb^{\infty}(Q)^{d \times d}$ resulting from  \eqref{ck-4}, we have
\begin{align}\label{impl-k-1}
0 \leq \liminf_{k \to \infty} \ll \BS^k(\cdot, \BD \oU^{k,{\L}}) - \BS^{\star}(\cdot, \BB), (\BD \oU^{k,{\L}}- \BB)\varphi \rr_{Q}
\end{align}
for all $\varphi \in C^{\infty}_0(Q)$ such that $\varphi \geq 0$ and for all matrices $\BB \in U$, for the dense set $U\subset\Rds$ given in the assumption.
Then, by the strong convergence of $\BD \oU^{k,{\L}}$ in \eqref{ck-4} and the weak convergence of $\BS^k(\cdot,\cdot, \BD \oU^{k,{\L}})$ in \eqref{ck-7} we obtain
\begin{align}\begin{split}\label{impl-k-2}
0 &\leq \liminf_{k \to \infty} \ll \BS^k(\cdot, \BD \oU^{k,{\L}}) - \BS^{\star}(\cdot, \BB), (\BD \oU^{k,{\L}}- \BB)\varphi \rr_{Q}\\
&=  \ll \BS^{{\L}} - \BS^{\star}(\cdot, \BB), (\BD \oU^{{\L}}- \BB)\varphi \rr_{Q}
\end{split}
\end{align}
for all $\varphi \in C^{\infty}_0(Q)$ such that $\varphi \geq 0$ and for all matrices $\BB \in U$.
By Lemma~\ref{sel} \ref{itm:a4} this allows us to conclude that
\begin{align*}
(\BD \oU^{{\L}}(\bz), \BS^{{\L}}(\bz)) &\in \mathcal{A}(\bz) \quad \text{ for a.e. }\bz \in Q,
\end{align*}
so \eqref{impl-l} is shown.
%
%
\end{proof}
\end{Lemma}	

\subsection*{Limit $l,n \to \infty$}  	
 We are taking the limits $l,n \to \infty$ simultaneously without imposing any condition on $\delta_l$ and $h_n$. The condition $q>\frac{2d}{d+2}$ is required to gain compactness. Two additional difficulties, compared to \cite{BGMS.2012}, arise from the discretization.
The first is that in order to prove a uniform bound on the sequence of approximations to the time derivative one would require the stability of the $\Leb^2$-projector onto $\Vd$ in Sobolev norms, which would impose stronger requirements on the finite element partition of $\Omega$. To avoid this, instead of the Aubin--Lions lemma we shall employ an alternative compactness result due to Simon (cf. Lemma \ref{Lem:cp-2}), which requires convergence properties of time-increments. 
The second difficulty is that, in the identification of the implicit relation we have to deal with the discrepancy between $\overline{\BS}^{\L}$ and $\BS^{\L}$, since $\overline{\BS}^{\L}$ appears in the equation  \eqref{ap-l-1} and $\BS^{\L}$ satisfies the implicit relation in \eqref{impl-l}.

\begin{Lemma}[Convergence as $l,n \to \infty$]\label{Lem:ln-conv} Let $\oU^{l,n,m} \in \mathbb{P}^l_0(0,T;\Vd)$, $\wU^{l,n,m} \in \mathbb{P}^l_1(0,T;\Vd)$, $\BS^{l,n,m} \in \Leb^{q'}(Q)^{d \times d}$ and
$\overline{\BS}^{l,n,m}\in \mathbb{P}^l_0(0,T;\Leb^{q'}(\Omega)^{d \times d})$ satisfy \eqref{ap-l-1}--\eqref{impl-l}, for any $l,n,m \in \mathbb{N}$, by Lemma~\ref{Lem:k-id}. Further, let $\eta$ be as defined in \eqref{def:eta}. For any $0 \leq s_0 < s \leq T$ and all ${\L}=(l,n,m) \in \mathbb{N}^3$ one has that
\begin{align}\begin{split} \label{en-ineq-ln-t}
\frac{1}{2}\norm{\wU^{{\L}}(s,\cdot)}_{\Leb^2(\Omega)}^2
+ \ll \overline{\BS}^{\L}, \BD \oU^{\L} \rr_{Q_{s_0}^s}
+& \frac{1}{m} \norm{\oU^{\L}}_{\Leb^{2q'}(Q_{s_0}^s)}^{2q'} \\
&\leq \ll \overline{\bf}, \oU^{\L}\rr_{Q_{s_0}^s} +
\frac{1}{2}\norm{\wU^{{\L}}(s_0,\cdot)}_{\Leb^2(\Omega)}^2.
\end{split}
\end{align}
Furthermore, for each $m \in \mathbb{N}$ there exists a $\bu^m \in \Leb^{\infty}(0,T;\l2d) \cap \Xqd$, $\BS^m \in \Leb^{q'}(Q)^{d \times d}$ and subsequences such that, as $l,n \to \infty$,
\begin{alignat}{2}
\label{cla-0}
 \wU^{l,n,m} &\to \bu^m  \quad &&  \text{ strongly in } \Leb^p(0,T;\Leb^{2}(\Omega)^d) \text{ for all } p \in [1, \infty),\\
\label{cla-1}
 \wU^{l,n,m}(s,\cdot) &\to \bu^m(s,\cdot)  \quad &&  \text{ strongly in } \Leb^{2}(\Omega)^d \text{ for a.e. } s\in(0,T),
 \\
 \label{cla-2}
  \wU^{l,n,m} (0,\cdot) &\to \bu_0 \quad && \text{ strongly in } \Leb^{2}(\Omega)^{d},\\
 \label{cla-3}
 \oU^{l,n,m} &\to \bu^m  \quad &&  \text{ strongly in } \Leb^p(0,T;\Leb^{2}(\Omega)^d) \cap \Leb^{r}(Q)^d \\
 &&& \quad \text{ for all } p \in [1, \infty) \text{ and all } r \in [1, \eta), \notag\\
 \label{cla-9}
 \oU^{l,n,m}(s,\cdot) &\to \bu^m(s,\cdot)  \quad &&  \text{ strongly in } \Leb^{2}(\Omega)^d \text{ for a.e. } s\in(0,T),
 \\
 \label{cla-4-v}
 \oU^{l,n,m} &\rightharpoonup \bu^m
 \quad && \text{ weakly in  }\Leb^q(0,T;\Sob^{1,q}_{0}(\Omega)^d)
 \cap \Leb^{\eta}(Q)^d, \\
   \label{cla-5}
 \oU^{l,n,m} & \overset{*}{\rightharpoonup} \bu^m
 \quad && \text{ weakly* in  } \Leb^{\infty}(0,T;\Leb^2(\Omega)^d),\\
 \label{cla-6}
\abs{\oU^{l,n,m}}^{2q'-2}\oU^{l,n,m} & \rightharpoonup \abs{\bu^{m}}^{2q'-2}\bu^{m} \quad  &&  \text{ weakly in } \Leb^{(2q')'}(Q)^{d},\\
 \label{cla-7}
  \overline{\BS}^{l,n,m} &\rightharpoonup \BS^m \quad && \text{ weakly in } \Leb^{q'}(Q)^{d\times d},\\
   \label{cla-8}
  \BS^{l,n,m} &\rightharpoonup \BS^m \quad && \text{ weakly in } \Leb^{q'}(Q)^{d\times d}.
 \end{alignat}
 \begin{proof} \hspace*{\fill}
\smallskip

\noindent
\textit{Step 1: Energy inequality}. Let $\L=(l,n,m) \in \mathbb{N}^3$, $i \in \{1,\ldots, l\}$ and let $t \in (t_{i-1},t_i]$. In \eqref{ap-l-1} we test with $\bW = \oU^{\L}(t,\cdot) \in \Vd$. For the first term adding and subtracting $\wU^{\L}(t,\cdot)$ one obtains with \eqref{ip-id-1} that
\begin{align}\begin{split}\label{en-ineq-l-1}
 \ll \partial_t \wU^{\L}(t,\cdot), \oU^{\L}(t, \cdot) \rr_{\Omega}
&=
 \ll \partial_t \wU^{\L}(t,\cdot), \wU^{\L}(t, \cdot) \rr_{\Omega}
 +
\ll \partial_t \wU^{\L}(t,\cdot), \oU^{\L}(t, \cdot)-\wU^{\L}(t, \cdot) \rr_{\Omega}\\
&= \frac{1}{2}\frac{\d}{\d t} \norm{\wU^{\L}(t, \cdot)}_{\Leb^2(\Omega)}^2 + (t_i-t)  \norm{ \partial_t \wU^{\L}(t, \cdot)}_{\Leb^2(\Omega)}^2\\
&\geq \frac{1}{2}\frac{\d}{\d t} \norm{\wU^{\L}(t, \cdot)}_{\Leb^2(\Omega)}^2,
\end{split}
\end{align}
since $t\leq t_i$. By the continuity of $\wU^{\L}$, upon integration over $(s_0,s)$, 
 for $0 \leq s_0 < s \leq T$, this yields
\begin{align}\label{en-ineq-l-2}
\int_{s_0}^s  \ll \partial_t \wU^{\L}(t,\cdot), \oU^{\L}(t, \cdot) \rr_{\Omega}
\geq  \frac{1}{2}  \norm{\wU^{\L}(s, \cdot)}_{\Leb^2(\Omega)}^2 - \frac{1}{2}  \norm{\wU^{\L}(s_0, \cdot)}_{\Leb^2(\Omega)}^2.
\end{align}
The other terms follow immediately and \eqref{en-ineq-ln-t} is proved.
\smallskip

\noindent
\textit{Step 2: Estimates on the discrete level}.
Similarly as in the proof of Lemma~\ref{Lem:k-ex}, testing \eqref{ap-l-i-1} with $\bW = \bU^{\L}_i\in \Vd$ for $i \in \{1, \ldots, l\}$, yields
\begin{align}\label{est:enid-l}
\ll d_t \bU^{\L}_i, \bU^{\L}_i \rr_{\Omega} + \ll \BS^{\L}_i, \BD \bU^{\L}_i \rr_{\Omega} + \frac{1}{m} \norm{\bU^{\L}_i}_{\Leb^{2q'}(\Omega)}^{2q'} = \ll \bf_i, \bU^{\L}_i \rr_{\Omega}.
\end{align}
The first term on the left-hand side is bounded as in \eqref{ap-est:k-1} and the term on the right-hand side
is bounded as in \eqref{ap-est:k-4}.
The only difference arises in bounding the term involving the stress tensor; cf. \eqref{ap-est:k-2}: we use that $(\BD \oU^{\L}, \BS^{\L}) \in \mathcal{A}$ a.e. in $Q$ by \eqref{impl-l} and Assumption~\ref{assump-A} \ref{itm:A4} to obtain
 \begin{align} 
\ll \BS^{\L}_i, \BD \bU^{\L}_i \rr_{\Omega}
&\overset{\eqref{SL-i}}{=}
\ll \fint_{t_{i-1}}^{t_i} \BS^{\L}(t,\cdot)\, \d t, \BD \bU^{\L}_i \rr_{\Omega}
= \frac{1}{\delta_l} \ll \BS^{\L}, \BD \bU^{\L}_i \rr_{Q_{i-1}^i} \notag \\
&\geq \frac{1}{\delta_l} \int_{Q^i_{i-1}} - \abs{g(\cdot)} + c_{*} \left( \abs{\BD \bU^{\L}_i}^q + \abs{\BS^{\L}}^{q'} \right)
 \d \bz \label{est:l-1} \\
 &\geq - \frac{1}{\delta_l} \norm{g}_{\Leb^1(Q_{i-1}^i)}
 + c \norm{\bU^{\L}_i}^q_{\Sob^{1,q}(\Omega)}
 + \frac{ c_*}{\delta_l} \norm{\BS^{\L}}_{\Leb^{q'}(Q_{i-1}^i)}^{q'},\notag
\end{align}
where again Poincar\'e's and Korn's inequalities were used.
Following the same procedure as in \eqref{ap-est:k-5}--\eqref{ap-est:k-7} we arrive at
 \begin{align}\begin{split}\label{est:l-2}
\max_{j \in \{1, \ldots, l \}} \norm{\bU^{\L}_j}_{\Leb^2(\Omega)}^2
&+ \sum_{j=1}^l \norm{\bU^{\L}_j - \bU^{\L}_{j-1}}_{\Leb^2(\Omega)}^2
+ \delta_l \sum_{j=1}^l \norm{\bU^{\L}_j}_{\Sob^{1,q}(\Omega)}^q \\
&
\,+  \norm{\BS^{\L}}_{\Leb^{q'}(Q)}^{q'}
+  \frac{\delta_l}{m}\sum_{j=1}^l \norm{\bU^{\L}_j}_{\Leb^{2q'}(\Omega)}^{2q'} \leq c \quad \text{ for all } \L=(l,n,m) \in \mathbb{N}^3.
\end{split}
 \end{align}
It follows by the relation between $\overline{\BS}^{\L}$ and $\BS^{\L}$ in \eqref{SL-i}, by \eqref{t-av-stab} and the estimate \eqref{est:l-2} that
 \begin{align}\label{est:l-3}
 \norm{\overline{\BS}^{\L}}_{\Leb^{q'}(Q)}^q
\leq  \norm{\BS^{\L}}_{\Leb^{q'}(Q)}^{q'}\leq c \quad \text{ for all } \L \in\mathbb{N}^3.
 \end{align}
\smallskip

\noindent
\textit{Step 3: Estimates on the continuous level}. By the definition of the piecewise constant interpolant  according to \eqref{def:interpol-c} it follows from the discrete
estimates that
\begin{align}\begin{split}\label{est:l-4}
\norm{\oU^{\L}}_{\Leb^{\infty}(0,T;\Leb^2(\Omega))}
&+
\norm{\oU^{\L}}_{\Leb^{2q'}(Q)}^{2q'}
+
\norm{\oU^{\L}}_{\Leb^{q}(0,T;\Sob^{1,q}(\Omega))}^{q}\\
&
\overset{\eqref{est:ip-est-1}}{=}
\max_{j \in \{1, \ldots, l\}} \norm{\bU^{\L}_j}_{\Leb^2(\Omega)}
+
\delta_l \sum_{j=1}^l \norm{\bU^{\L}_j}_{\Leb^{2q'}(\Omega)}^{2q'}
+
\delta_l \sum_{j=1}^l \norm{\bU^{\L}_j}_{\Sob^{1,q}(\Omega)}^{q}\\
& \overset{\eqref{est:l-2}}{\leq} c(m) \quad \text{ for all } {\L}=(l,n,m) \in \mathbb{N}^3.
\end{split}
\end{align}
With this and the parabolic interpolation from Lemma~\ref{par-interpol} we have that
\begin{align}\label{est:l-4v}
\norm{\oU^{\L}}_{\Leb^{\frac{q(d+2)}{d}}(Q)}\leq c(m)
\quad \text{ for all } {\L} =(l,n,m)\in \mathbb{N}^3.
\end{align}

For the estimates of the continuous, piecewise affine interpolant $\wU^{\L}$ according to
 \eqref{def:interpol-l} one also has to estimate the corresponding norms of $\bU^{\L}_0$; by \eqref{ap-l-i-2} and the stability of the $\Leb^2$-projector in
  \eqref{est:L2-stab} we have that
 \begin{align}\label{est:l-5}
 \norm{\bU^{\L}_0}_{\Leb^2(\Omega)} \overset{\eqref{ap-l-i-2}}{=} \norm{P^n_{\diver} \bu_0}_{\Leb^2(\Omega)} \overset{\eqref{est:L2-stab}}{\leq} \norm{\bu_0}_{\Leb^2(\Omega)} \quad \text{ for all } \L = (l,n,m) \in \mathbb{N}^3.
 \end{align}
Together with the discrete estimate in \eqref{est:l-2} this yields that
 \begin{align}\begin{split}\label{est:l-6}
 \norm{\wU^{\L}}_{\Leb^{\infty}(0,T;\Leb^2(\Omega)}
\overset{\eqref{est:ip-est-2}}{=}
 \max_{j \in \{0, \ldots, l\}} \norm{\bU^{\L}_j}_{\Leb^2(\Omega)}
\overset{\eqref{est:l-5},\eqref{est:l-2}}{\leq} c
\quad \text{ for all } \L\in \mathbb{N}^3.
\end{split}
 \end{align}
For the compactness argument instead of $\wU^{\L}$
we consider $\widehat{\bU}^{\L} \in C([0,T];\Vd)$ defined by
\begin{align}\label{def:u-hat}
\widehat{\bU}^{\L}(t,\cdot) \coloneqq \begin{cases}
\wU^{\L}(t,\cdot) \quad &\text{ if } t \in (\delta_l, T],\\
\oU^{\L}(t,\cdot)=\bU^{\L}_1(\cdot) &\text{ if } t \in [0, \delta_l],
\end{cases}
\end{align}
which is constant on $[0,\delta_l]$.
Note that by definition and by the discrete estimate we have that
\begin{align}\label{est:l-7}
\norm{\widehat{\bU}^{\L}}_{\Leb^{\infty}(0,T;\Leb^2(\Omega))} &\leq \norm{\widetilde {\bU}^{\L}}_{\Leb^{\infty}(0,T;\Leb^2(\Omega))} \overset{\eqref{est:l-6}}{\leq} c,\\
\label{est:l-12}
\norm{\widehat{\bU}^{\L}}_{\Leb^{2q'}(Q)}^{2q'}
&= \norm{\bU^{\L}_1}_{\Leb^{2q'}(Q^{\delta_l}_0)}^{2q'} + \norm{\wU^{\L}}_{\Leb^{2q'}(Q_{\delta_l}^T)}^{2q'} \leq c \delta_l \sum_{i=1}^l \norm{\bU^{\L}_i}_{\Leb^{2q'}(\Omega)}^{2q'} \overset{\eqref{est:l-2}}{\leq} c(m),\\
\label{est:l-8}
\norm{\widehat{\bU}^{\L}}_{\Leb^{q}(0,T;\Sob^{1,q}(\Omega))}^q
&\leq c \delta_l \sum_{i=1}^l \norm{\bU^{\L}_i}_{\Sob^{1,q}(\Omega)}^q \overset{\eqref{est:l-2}}{\leq} c
\end{align}
for all ${\L}=(l,n,m) \in \mathbb{N}^3$.

By the fact that $\partial_t \wU^l(t,\cdot) = d_t \bU^{\L}_i$, for $t \in (t_{i-1},t_i]$, $i \in \{1, \ldots, l\}$, and by the discrete estimate
 \eqref{est:l-2} we obtain
\begin{align}\label{est:l-9} 
\delta_l \norm{\partial_t \wU^{\L}}_{\Leb^{2}(Q)}^2 &=
\delta_l  \sum_{i=1}^l \norm{\frac{1}{\delta_l } ({\bU}^{\L}_i-{\bU}^{\L}_{i-1})}_{\Leb^{2}(Q_{i-1}^i)}^2
 =  \sum_{i=1}^l\norm{ {\bU}^{\L}_i-{\bU}^{\L}_{i-1}}_{\Leb^{2}(\Omega)}^2 \overset{\eqref{est:l-2}}{\leq} c
\end{align}
for all ${\L} \in \mathbb{N}^3$.

  Finally, we also estimate $\mathfrak{L}^{\L}[\bu;\bv](t)$, as defined in \eqref{Llnm-t}: by
   \eqref{est:wb-r}, duality of norms, and H\"older's and Poincar\'e's inequality we obtain
\begin{align}\begin{split}\label{est:l-10}
\int_{a}^{b} \mathfrak{L}^{\L}[\bu;\bv](t) \, \d t
&=
 \ll \widetilde{b}(\bu(t, \cdot),\bu(t, \cdot), \bv)\rr_{(a,b)} - \ll \overline{\BS}^{{\L}}(t,\cdot), \BD\bv\rr_{Q_a^b}\\
&\quad - \frac{1}{m}\ll \abs{\bu(t,\cdot)}^{2q'-2} \bu(t,\cdot), \bv \rr_{Q_a^b}
+ \ll \overline{\bf}(t, \cdot), \bv\rr_{Q_a^b}\\
&\leq
\norm{\bu}_{\Leb^{2q'}(Q_a^b)}^{2} \norm{\nabla \bv}_{\Leb^{q}(Q_a^b)}
 + \norm{\bu}_{\Leb^{2q'}(Q_a^b)}  \norm{\nabla \bu}_{\Leb^{q}(Q_a^b)} \norm{\bv}_{\Leb^{2q'}(Q_a^b)}\\
&\quad +
\norm{\overline{\BS}^{{\L}}}_{\Leb^{q'}(Q_a^b)} \norm{\BD \bv}_{\Leb^{q}(Q_a^b)}
+ \frac{1}{m} \norm{\bu}_{\Leb^{2q'}(Q_a^b)}^{2q'-1} \norm{\bv}_{\Leb^{2q'}(Q_a^b)} \\
&\quad + \norm{\overline{\bf}}_{\Leb^{q'}(a,b;\Sob^{-1,q'}(\Omega))}
\norm{\bv}_{\Leb^{q}(a,b;\Sob^{1,q}(\Omega))}\\
&\leq
c \left(1+
\norm{\bu}_{\Leb^{2q'}(Q_a^b)}^{2}
 \right) \norm{ \nabla \bv}_{\Leb^{q}(Q_a^b)}
 \\
 &\quad
 + c \left(
   \norm{\bu}_{\Leb^{2q'}(Q_a^b)}  \norm{\nabla \bu}_{\Leb^{q}(Q_a^b)}
   +
   \frac{1}{m} \norm{\bu}_{\Leb^{2q'}(Q_a^b)}^{2q'-1}
 \right)
 \norm{\bv}_{\Leb^{2q'}(Q_a^b)},
\end{split}
\end{align}
for $0\leq a < b \leq T$, for any ${\L}=(l,n,m)\in \mathbb{N}^3$, where we have used the estimate
 \eqref{est:l-3} on $\overline{\BS}^{\L}$ and
  \eqref{of-bd} on $\overline{\bf}$.
With the estimates on $\oU^{\L}$ in \eqref{est:l-4} this yields
\begin{align}\begin{split}\label{est:l-11}
\int_{a}^{b} \mathfrak{L}^{\L}[\oU^{\L};\bv](t) \, \d t
&\overset{\eqref{est:l-10}}{\leq}
c \left( 1 +
\norm{\oU^{\L}}_{\Leb^{2q'}(Q_a^b)}^{2}
 \right) \norm{ \nabla \bv}_{\Leb^{q}(Q_a^b)}
 \\
 &\quad
 + c \left(
   \norm{\oU^{\L}}_{\Leb^{2q'}(Q_a^b)}  \norm{\nabla \oU^{\L}}_{\Leb^{q}(Q_a^b)}
   +
   \frac{1}{m} \norm{\oU^{\L}}_{\Leb^{2q'}(Q_a^b)}^{2q'-1}
 \right)
 \norm{\bv}_{\Leb^{2q'}(Q_a^b)}\\
 &\overset{\eqref{est:l-4}}{\leq} c(m) \left(  \norm{ \nabla \bv}_{\Leb^{q}(Q_a^b)} + \norm{\bv}_{\Leb^{2q'}(Q_a^b)}\right),
\end{split}
\end{align}
for $0\leq a < b \leq T$ and any ${\L}=(l,n,m)\in \mathbb{N}^3$.
 \smallskip

\noindent
\textit{Step 4: Convergence of the time increments} (compare \cite[pp.~174]{CHP.2010}). Instead of applying the Aubin--Lions lemma, as in \cite{BGMS.2012}, here we apply the compactness result due to Simon, stated in Lemma~\ref{Lem:cp-2}. This means that we do not need uniform bounds on the time derivatives but only convergence properties for time increments, which avoids the use of stability results in Sobolev norms for the $\Leb^2$-projector onto $\Vd$.
We wish to apply Lemma~\ref{Lem:cp-2} to the sequence $\{\widehat{\bU}^{l,n,m}\}_{l,n \in \mathbb{N}}$, for fixed $m \in \mathbb{N}$, with $X=\Sob^{1,q}(\Omega)^d$, $B=\Leb^2(\Omega)^d$ and $p=2$. Let us show that
\begin{align}\label{est:cp-1}
\int_0^{T-\varepsilon} \norm{\widehat{\bU}^{\L}(s+\varepsilon,\cdot)-\widehat{\bU}^{\L}(s,\cdot)}_{\Leb^2(\Omega)}^2 \, \d s \to 0, \quad \text{ as } \varepsilon \to 0, \quad \text{ uniformly for } l,n \in \mathbb{N}.
\end{align}
Consider the term $
\ll \widehat{\bU}^{\L}(s+\varepsilon,\cdot)-\widehat{\bU}^{\L}(s,\cdot), \bW \rr_{\Omega}$,
for $\bW \in \Vd$, $s \in (0,T)$ and $\varepsilon>0$ such that $s+\varepsilon <T$.
If $s+\varepsilon \leq \delta_l$, then we have $\widehat{\bU}^{\L}(s+\varepsilon) = \widehat{\bU}^{\L}(s) = \bU^{\L}_1$, so the term vanishes. Now let $s+\varepsilon > \delta_l$.
By the definition of $\widehat{\bU}^{\L}$ in
 \eqref{def:u-hat} we have that $\widehat{\bU}^{\L}(s,\cdot)= \widehat{\bU}^{\L}(\max(s,\delta_l),\cdot)$.
By the continuity of $\widehat{\bU}^{\L}$ and since $\partial_t \widehat{\bU}^{\L}$ is integrable, we obtain
\begin{align}\begin{split} \label{est:cp-2}
\ll \widehat{\bU}^{\L}(s+\varepsilon,\cdot)-\widehat{\bU}^{\L}(s,\cdot), \bW \rr_{\Omega}
&=
\int_{\max(s,\delta_l)}^{s+\varepsilon} \ll \partial_t \widehat{\bU}^{\L}(t, \cdot), \bW \rr_{\Omega} \, \d t \\
&=
\int_{\max(s,\delta_l)}^{s+\varepsilon} \ll \partial_t \wU^{\L}(t, \cdot), \bW \rr_{\Omega} \, \d t,
\end{split}
\end{align}
where in the last line we have used that $\widehat{\bU}^{\L}(t,\cdot)$ and $\wU^{\L}(t,\cdot)$ coincide on $(\max(s,\delta_l), s+\varepsilon) \subset [\delta_l, T]$. Applying the equation \eqref{ap-l-1} for a.e. $t \in (\max(s,\delta_l), s+\varepsilon)$, integrating and applying the bounds in
 \eqref{est:l-11} yields
\begin{align}\begin{split} \label{est:cp-3}
\int_{\max(s,\delta_l)}^{s+\varepsilon} \ll \partial_t \wU^{\L}(t, \cdot), \bW \rr_{\Omega} \, \d t &=
\int_{\max(s,\delta_l)}^{s+\varepsilon} \mathfrak{L}^{\L}[\oU^{\L}; \bW](t)\, \d t\\
&\overset{\eqref{est:l-11}}{\leq}
c(m) \left( \norm{ \nabla \bW}_{\Leb^{q}\left(Q_{\max(s,\delta_l)}^{s+\varepsilon}\right)} + \norm{\bW}_{\Leb^{2q'}\left(Q_{\max(s,\delta_l)}^{s+\varepsilon}\right)}\right)\\
&= c(m) \left(\varepsilon^{\sfrac{1}{q}} + \varepsilon^{\sfrac{1}{2q'}}\right)
\norm{  \bW}_{X^{q}(\Omega)},
\end{split}
\end{align}
since $\bW$ is constant in time and the length of the time interval is
\begin{align*}
s+\varepsilon- \max(s, \delta_l) = \min(\varepsilon, s+\varepsilon-\delta_l) \leq \varepsilon.
\end{align*}

For all $s \in (0,T)$ and $\varepsilon>0$ such that $s+\varepsilon<T$ we have that $ \widehat{\bU}^{\L}(s+\varepsilon,\cdot)$, $\widehat{\bU}^{\L}(s,\cdot) \in \Vd$; so, applying \eqref{est:cp-2} and
 \eqref{est:cp-3} with $\bW =  \widehat{\bU}^{\L}(s+\varepsilon,\cdot)-\widehat{\bU}^{\L}(s,\cdot)$, which is piecewise constant in time, shows that
\begin{align}\begin{split} \label{est:cp-4}
\norm{\widehat{\bU}^{\L}(s+\varepsilon,\cdot)-\widehat{\bU}^{\L}(s,\cdot)}_{\Leb^2(\Omega)}^2
&\leq c(m) \left( \varepsilon^{\sfrac{1}{q}}
+ \varepsilon^{\sfrac{1}{2q'}} \right)
 \norm{  \widehat{\bU}^{\L}(s+\varepsilon,\cdot)
- \widehat{\bU}^{\L}(s,\cdot)}_{X^{q}(\Omega)}.
\end{split}
\end{align}
Integrating over $(0,T-\varepsilon)$, using the triangle inequality, H\"older's inequality and the estimates in
 \eqref{est:l-12} and \eqref{est:l-8} yields
\begin{align}\begin{split} \label{est:cp-5}
& \int_0^{T-\varepsilon} \norm{\widehat{\bU}^{\L}(s+\varepsilon,\cdot)-\widehat{\bU}^{\L}(s,\cdot)}_{\Leb^2(\Omega)}^2 \, \d s
 \\
&\quad\quad \overset{\eqref{est:cp-4}}{\leq} c(m)  \left( \varepsilon^{\sfrac{1}{q}} +
\varepsilon^{\sfrac{1}{2q'}}
\right)
 \int_{0}^{T-\varepsilon}\left(
 \norm{  \widehat{\bU}^{\L}(s+\varepsilon,\cdot)}_{X^{q}(\Omega)} +
 \norm{ \widehat{\bU}^{\L}(s,\cdot)}_{X^{q}(\Omega)} \right) \, \d s \\
 &\quad\quad \leq c(m)
  \left(
  \varepsilon^{\sfrac{1}{q}} + \varepsilon^{\sfrac{1}{2q'}} \right)
   \norm{ \widehat{\bU}^{\L}}_{X^q(Q)}  \overset{\eqref{est:l-12},\eqref{est:l-8}}{\leq} c(m) (\varepsilon^{\sfrac{1}{q}} +  \varepsilon^{\sfrac{1}{2q'}}  ) \to 0,
\end{split}
\end{align}
as $\varepsilon \to 0$ uniformly in $l,n \in \mathbb{N}$, where $\L=(l,n,m) \in\mathbb{N}^3$. This proves \eqref{est:cp-1}.
 \smallskip

\noindent
\textit{Step 5: Convergence as $l,n \to \infty$}. Recall that we have $\L=(l,n,m) \in \mathbb{N}^3$ and let $m \in \mathbb{N}$ be fixed. By the estimates \eqref{est:l-7} and \eqref{est:l-8} we have that $\{\widehat{\bU}^{l,n,m}\}_{l,n \in \mathbb{N}}$  is bounded in particular in $\Leb^2(Q)^d$ and $\Leb^1(0,T;\Sob^{1,q}(\Omega)^d)$. By the condition that $q>\frac{2d}{d+2}$, the embedding $\Sob^{1,q}(\Omega)\cpemb \Leb^2(\Omega)$ is compact and with \eqref{est:cp-1} all the assumptions in Lemma~\ref{Lem:cp-2} are satisfied for $X=\Sob^{1,q}(\Omega)^d$, $B=\Leb^2(\Omega)^d$ and $p=2$. Hence, there exists $\bu^m\in \Leb^2(Q)^d$ and a subsequence such that
\begin{align}\label{conv-l-1}
\widehat{\bU}^{l,n,m} \to \bu^m \quad \text{ strongly in } \Leb^2(Q)^d,\quad \text{ as } l,n \to \infty.
\end{align}
By the definition of $\widehat{\bU}^{l,n,m}$ in \eqref{def:u-hat} and the property \eqref{ip-id-1} of the interpolants defined in \eqref{def:interpol-c} and \eqref{def:interpol-l} we have that
\begin{align}
\norm{\widehat{\bU}^{l,n,m}- \wU^{l,n,m}}_{\Leb^2(Q)}^2
&=
\norm{\oU^{l,n,m}- \wU^{l,n,m}}_{\Leb^2(0,\delta_l; \Leb^2(\Omega))}^2 \notag
\\ \label{conv-l-2}
&
\overset{\eqref{ip-id-1}}{\leq}
\norm{(\delta_l - t)\partial_t \wU^{l,n,m}}_{\Leb^2(0,\delta_l; \Leb^2(\Omega))}^2
 \leq \delta_l^2  \norm{\partial_t \wU^{l,n,m}}_{\Leb^2(0,\delta_l; \Leb^2(\Omega))}^2\\  \notag
 & \overset{\eqref{est:l-9}}{\leq} c \delta_l
\to 0, \quad \text{ as } l \to \infty.
\end{align}
With \eqref{conv-l-1} it follows that $\wU^{l,n,m} \to \bu^m$ strongly in $\Leb^{2}(Q)^d$, as $l,n \to \infty$.
By the boundedness in $\Leb^{\infty}(0,T;\Leb^2(\Omega)^d)$ in \eqref{est:l-6} and interpolation, this implies that
\begin{align}\label{conv-l-3}
\wU^{l,n,m} \to \bu^m \quad \text{ strongly in } \Leb^{p}(0,T;\Leb^2(\Omega)^d),\quad \text{ as } l,n \to \infty,
\end{align}
 for any $p \in [1, \infty)$.
Similarly, by \eqref{ip-id-1} we have that
\begin{align}
\begin{split}\label{conv-l-4}
\norm{\oU^{l,n,m}- \wU^{l,n,m}}_{\Leb^2(Q)}^2
&\overset{\eqref{ip-id-1}}{=} \sum_{i=1}^l
 \norm{(t_i-t) \partial_t \wU^{l,n,m}}_{\Leb^2(Q_{i-1}^i)}^2  \\
& \leq \delta_l^2  \norm{\partial_t \wU^{l,n,m}}_{\Leb^2(Q)}^2 \overset{\eqref{est:l-9}}{\leq} c \delta_l
\to 0, \quad \text{ as } l \to \infty.
\end{split}
\end{align}
Consequently, with \eqref{conv-l-3} it follows that
$\oU^{l,n,m} \to \bu^m$ strongly in $\Leb^{2}(Q)^d$, as $l,n \to \infty$.
In particular, $t \mapsto \|\oU^{l,n,m}(t,\cdot)- \bu^m(t,\cdot)\|_{\Leb^2(\Omega)}$ converges to zero strongly in $\Leb^2(0,T)$, as $l,n \to \infty$. Thus, there exists a subsequence such that $t \mapsto \|\oU^{l,n,m}(t,\cdot)- \bu^m(t,\cdot)\|_{\Leb^2(\Omega)}$ converges to zero a.e. in $(0,T)$, as $l,n \to \infty$, which implies \eqref{cla-9}. Analogously, \eqref{cla-2} follows from the strong convergence of $\wU^{l,n,m}$ in \eqref{conv-l-3}.

The uniform bounds in $\Leb^{\infty}(0,T;\Leb^2(\Omega)^d)$ and $\Leb^{\eta}(Q)^d$, with $\eta = \max(2q',\frac{q(d+2)}{d})$, by \eqref{est:l-4} and \eqref{est:l-4v}, and the strong convergence in $\Leb^2(Q)^d$, yield by interpolation, that
\begin{align}\label{conv-l-5}
\oU^{l,n,m} \to \bu^m \quad \text{ strongly in } \Leb^{p}(0,T;\Leb^2(\Omega)^d)\cap \Leb^r(Q)^d, \text{ as } l,n \to \infty,
\end{align}
for any $p \in [1, \infty)$ and any $r \in [1, \eta)$.
By the uniform bounds in \eqref{est:l-4} and \eqref{est:l-4v} and the (sequential) Banach--Alaoglu theorem, up to subsequences, we have that
\begin{alignat}{2}\label{conv-l-6}
\oU^{l,n,m} &\overset{*}{\wconv} \bu^m \quad &&\text{ weakly* in } \Leb^{\infty}(0,T;\Leb^2(\Omega)^d),\\\label{conv-l-7}
\oU^{l,n,m} &\wconv \bu^m \quad &&\text{ weakly in } \Leb^{q}(0,T;\Sob^{1,q}_0(\Omega)^d)\cap \Leb^{\eta}(Q)^d,
\end{alignat}
as $l,n \to \infty$, where the identification of the limiting functions with $\bu^m$ follows by the strong convergence in \eqref{conv-l-5}.

The argument that $\bu^m$ is divergence-free follow as in \cite[p.~1001]{DKS.2013}: Let $h \in \Leb^{q'}(\Omega)$ and note that by the Assumption~\ref{assump-PQ} on the projector $\Pn_{\mathbb{Q}}$ we have that $\Pn_{\mathbb{Q}} h \to h$ in particular in $\Leb^{q'}(\Omega)$, as $n \to \infty$, compare Remark~\ref{Rem:proj} \ref{itm:approx-PQ}. Also, let $\varphi \in C^{\infty}_0(0,T)$.
 By \eqref{conv-l-7} we have that $\diver \oU^{l,n,m} \rightharpoonup \diver \bu^m$ weakly in $\Leb^{q}(\Omega)$, and hence
\begin{align}
\ll \diver \oU^{l,n,m},  \varphi \Pn_{\mathbb{Q}}h  \rr_{Q}
\to \ll \diver \bu^m, \varphi h \rr_{Q}, \quad \text{ as } l,n \to \infty.
\end{align}
Since $\oU^{l,n,m} \in \mathbb{P}^l_0(0,T;\Vd)$ the left-hand side vanishes for all $l,n \in \mathbb{N}$, and hence we have that $\ll \diver \bu^m, h \varphi \rr_{Q} =0$ for all $h \in \Leb^{q'}(\Omega)$ and all $\varphi \in C^{\infty}(0,T)$, so by density $\bu^m$ is (weakly) divergence-free.

By \eqref{est:l-4} it follows that $\{|\oU^{l,n,m}|^{2q'-2}\,\oU^{l,n,m}\}_{l,n \in \mathbb{N}}$ is bounded in $\Leb^{(2q')'}(Q)^d$ and thus, by the (sequential) Banach--Alaoglu theorem there exists a subsequence and $\psi^m \in \Leb^{(2q')'}(Q)^d$ such that
\begin{align}\label{conv-l-8}
\abs{\oU^{l,n,m}}^{2q'-2}\oU^{l,n,m} \wconv \psi^m \quad \text{ weakly in } \Leb^{(2q')'}(Q)^d, \quad \text{ as } l,n \to \infty.
\end{align}
By the strong convergence in \eqref{conv-l-5}, there exists a subsequence, which converges a.e. in $Q$, and hence we can identify $\psi^m= \abs{\bu^m}^{2q'-2}\bu^m$, which shows \eqref{cla-6}.

Because $\wU^{l,n,m}(0,\cdot)= P^n_{\diver} \bu_0$ by \eqref{ap-l-2}, with \eqref{conv-Pnd} it follows that
\begin{align}\label{conv-l-9}
\wU^{l,n,m}(0,\cdot)=P^n_{\diver} \bu_0 \to \bu_0, \quad \text{ strongly in } \Leb^2(\Omega)^d, \quad \text{ as } n \to \infty,
\end{align}
so \eqref{cla-2} is proven.

The uniform estimates in \eqref{est:l-3} and the (sequential) Banach--Alaoglu theorem imply that there exist $\overline{\BS}^m, \BS^m \in \Leb^{q'}(Q)^{d \times d}$ such that
\begin{align}\label{conv-l-10}
\overline{\BS}^{l,n,m} &\wconv \overline{\BS}^{m}, \quad &&\text{ weakly in } \Leb^{q'}(Q)^{d\times d},\\\label{conv-l-11}
\BS^{l,n,m} &\wconv \BS^{m}, \quad &&\text{ weakly in } \Leb^{q'}(Q)^{d\times d},
\end{align}
as $l,n \to \infty$. It remains to show that $\overline{\BS}^m=\BS^m$; to this end, let $\BB \in C^{\infty}_0(Q)^{d \times d}$ be arbitrary but fixed. On the one hand the weak convergence in \eqref{conv-l-10} shows that
\begin{align}\label{conv-l-12}
\ll \overline{\BS}^{l,n, m}, \BB \rr_{Q} \to \ll \overline{\BS}^{ m}, \BB \rr_{Q}, \quad \text{ as } l,n \to \infty.
\end{align}
On the other hand, by the relation between $\overline{\BS}^{l,n, m}$ and $\BS ^{l,n,m}$ according to \eqref{SL-i}, we have that
\begin{align}\begin{split}\label{conv-l-13} 
\ll \overline{\BS}^{l,n, m}, \BB \rr_{Q}
&= \sum_{i=1}^l \int_{t_{i-1}}^{t_i} \int_{\Omega} \fint_{t_{i-1}}^{t_i} \BS^{l,n, m}(s, \bx) \, \d s \, \BB(t, \bx) \, \d \bx \, \d t \\
&= \sum_{i=1}^l \frac{1}{\delta_l} \int_{\Omega} \left( \int_{t_{i-1}}^{t_i} \BS^{l,n, m}(s, \bx) \, \d s \right) \left( \int_{t_{i-1}}^{t_i} \BB(t, \bx) \, \d t \right) \, \d \bx \\
&= \sum_{i=1}^l \int_{\Omega} \int_{t_{i-1}}^{t_i} \BS^{l,n, m}(s, \bx)  \fint_{t_{i-1}}^{t_i} \BB(t, \bx) \, \d t \, \d \bx \, \d s \\
&= \ll \BS^{l,n, m}, \overline{\BB} \rr_{Q} \to \ll \BS^{m}, \BB \rr_{Q}, \quad \text{ as } l,n \to \infty,
\end{split}
\end{align}
where we have used \eqref{conv-l-11} and the fact that $\overline{\BB} \to \BB$ strongly in $\Leb^q(Q)^{d \times d}$, as $l \to \infty$, see \eqref{t-av-conv}. By \eqref{conv-l-12} and \eqref{conv-l-13}, the uniqueness of limits implies that $\overline{\BS}^m = \BS^m$ a.e. in $Q$.
 \end{proof}
 \end{Lemma}

For $m \in \mathbb{N}$, $t \in (0,T)$, $\bu \in \Leb^{2q'}(Q)^d$ and $\bv \in X^q(\Omega)$ let us introduce
\begin{align}\begin{split}\label{Lm-t}
\mathfrak{L}^{m}[\bu; \bv](t)
\coloneqq&
- b(\bu(t, \cdot),\bu(t, \cdot), \bv) - \ll \BS^{m}(t,\cdot), \BD\bv\rr_{\Omega}\\
&- \frac{1}{m}\ll \abs{\bu(t,\cdot)}^{2q'-2} \bu(t,\cdot), \bv \rr_{\Omega}
+ \ll \bf(t, \cdot), \bv\rr_{\Omega},
\end{split}
\end{align}
where $\BS^m$ is given by Lemma~\ref{Lem:ln-conv} and  $b(\cdot,\cdot,\cdot)$ is defined in \eqref{def:b}.

 \begin{Lemma}[Identification of the PDE as $l,n \to \infty$]\label{Lem:ln-id} 
  The limiting functions $\bu^m \in \Leb^{\infty}(0,T;\l2d)\cap \Xqd$ given in Lemma~\ref{Lem:ln-conv} satisfy that $\partial_t \bu^m \in \Leb^{\tau}(0,T;(X_{\diver}^q(\Omega))')$, with $\tau$ defined in  \eqref{def:tau}, and $X_{\diver}^q(\Omega)$ defined in  \eqref{def:Xq-s}. (Up to a representative) we have that $\bu^m \in
\Cw([0,T],\l2d)$ for all $m \in \mathbb{N}$.
Furthermore, for each $m \in \mathbb{N}$ the functions $\bu^m$ and $\BS^m \in \Leb^{q'}(Q)^{d \times d}$ from Lemma~\ref{Lem:ln-conv} satisfy
\begin{alignat}{2}
\ll \partial_t \bu^m(t,\cdot), \bw\rr_{\Omega}
&=  \mathfrak{L}^m[\bu^m;\bw](t)  \quad &&\text{ for all } \bw \in C^{\infty}_{0, \diver}(\Omega)^d
\label{ap-m-1}
\text{ for a.e. } t \in (0,T),
\\
(\BD \bu^m(\bz), \BS^m(\bz)) &\in \mathcal{A}(\bz)   \quad &&\text{ for a.e. } \bz \in Q, \label{impl-m}\\
 \esslim_{t \to 0_+ }
  \norm{\bu^m(t, \cdot) - \bu_0}_{\Leb^2(\Omega)} &  \,=  0. && \label{ap-m-2}
\end{alignat}
\begin{proof}
Let $m \in \mathbb{N}$ be arbitrary but fixed. \hspace*{\fill} \smallskip

\noindent
\textit{Step 1: Identification of the limiting equation}. For $\L=(l,n,m) \in \mathbb{N}^3$ multiplying  \eqref{ap-l-1} by $\varphi \in C^{\infty}_0((-T,T))$ and integrating over $(0,T) $ yields
\begin{align}\label{ln-1v}
\ll \partial_t \wU^{\L}, \bW \varphi \rr_{Q}  \overset{\eqref{ap-l-1}}{=} \ll \mathfrak{L}^{\L}[\oU^{\L}, \bW],\varphi\rr_{(0,T)}, \quad \text{ for any } \bW \in \Vd.
\end{align}
By integration by parts and the fact that
$\wU^{\L}\in C([0,T];\Leb^2(\Omega)^d)$ we obtain that 
\begin{align}\begin{split}\label{ln-2v}
\ll \partial_t \wU^{\L}, \bW \varphi \rr_{Q}  &= \ll \partial_t(\wU^{\L}\varphi), \bW \rr_Q - \ll \wU^{\L}, \bW \partial_t \varphi \rr_{Q}\\
&=-\ll \wU^{\L}(0,\cdot), \varphi(0) \bW \rr_{\Omega} - \ll \wU^{\L}, \bW \partial_t \varphi \rr_{Q}.
\end{split}
\end{align} 
Taking both \eqref{ln-1v} and \eqref{ln-2v} together, we obtain that
\begin{align}\label{ap-l-1-w}
- \ll \wU^{\L}, \bW \partial_t \varphi \rr_{Q} =  \ll \wU^{\L}(0,\cdot), \varphi(0) \bW \rr_{\Omega}
+ \ll \mathfrak{L}^{\L}[\oU^{\L}, \bW],\varphi\rr_{(0,T)}
\end{align}
for all $\bW \in \Vd$ and all $\varphi \in C^{\infty}_0((-T,T))$ and $\L=(l,n,m) \in \mathbb{N}^3$.

Now let $\bw \in C^{\infty}_{0, \diver}(\Omega)^d$ and $\varphi \in C^{\infty}_0((-T,T))$ be arbitrary.
Recall that by Remark~\ref{Rem:proj} \ref{itm:approx-Pin} for $\bw \in C^{\infty}_{0, \diver}(\Omega)^d$ we have that
\begin{align}\label{Pn-conv}
 \Vd\ni \Pn \bw \to \bw \quad \text{ strongly in } \Sob^{1,s}_0(\Omega)^d,\quad \text{ as } n \to \infty,\quad \text{ for any } s\in[1,\infty).
 \end{align}
 In order to deduce the limiting equation for $\bu^m$ we consider \eqref{ap-l-1-w} term by term, as $l,n \to \infty$:
 let $s\in [1, \infty)$ be large enough that the embedding $\Sob^{1,s}(\Omega)^d\ctemb \Leb^2(\Omega)^d$ is continuous.
By the strong convergence of $\wU^{l,n,m}$ in $\Leb^{p}(0,T;\Leb^2(\Omega)^d)$, for $p \in [1, \infty)$ by \eqref{cla-0}, with \eqref{Pn-conv} we obtain that
 \begin{align}\label{ln-3v}
 - \ll  \wU^{l,n,m}, \Pn (\bw) \partial_t \varphi\rr_{Q}
\to
- \ll  \bu^{m}, \bw \partial_t \varphi\rr_{Q}, \quad \text { as } l,n \to \infty.
\end{align}
Similarly, the strong convergence of $\wU^{l,n,m}(0,\cdot)\to \bu_0$ in $\Leb^2(\Omega)^d$ in \eqref{cla-2} yields that
 \begin{align}\label{ln-4v}
 \ll \wU^{l,n,m}(0,\cdot), \varphi(0) \Pn \bw \rr_{\Omega}
\to
\ll \bu_0, \varphi(0) \bw \rr_{\Omega}, \quad \text { as } l,n \to \infty.
\end{align}
By the fact that $\oU^{l,n,m} \to \bu^m$ strongly in $\Leb^r(Q)^d$ for all $r \in [1, \eta)$, as $l,n \to \infty$, by \eqref{cla-3}, it follows that $\oU^{l,n,m} \otimes \oU^{l,n,m} \to \bu^m \otimes \bu^m$ strongly in $\Leb^p(Q)^{d \times d}$ for all $p \in [1, \tfrac{\eta}{2})$. Such a $p>1$ exists, since $\eta = \max\left(2q',\frac{q(d+2)}{d}\right)>2$ by the assumption that $q>\frac{2d}{d+2}$.
With \eqref{Pn-conv} applied for $s=p'<\infty$ we obtain that $\varphi \nabla \Pn \bw \to \varphi \nabla \bw$ strongly in $\Leb^{p'}(Q)^{d \times d}$. Together these imply that
\begin{align} \label{ln-5v}
\ll \oU^{l,n,m} \otimes \oU^{l,n,m}, \varphi \nabla \Pn \bw \rr_{Q}
\to
\ll \bu^m \otimes \bu^m, \varphi \nabla \bw \rr_{Q},  \quad \text{ as }\; l,n \to \infty.
\end{align}
For the modification of the convective term note first that we have weak convergence of $\nabla \oU^{l,n,m}\wconv \nabla \bu^m$ in $\Leb^{q}(Q)^{d\times d}$ by \eqref{cla-4-v}. By \eqref{cla-3} we have in particular that $\oU^{l,n,m} \to \bu^m$ strongly in $\Leb^{q'}(Q)^d$, as $l,n \to \infty$ since $q'<2q'\leq \eta$. For $s>d$, the embedding $\Sob^{1,s}(\Omega)\ctemb \Leb^{\infty}(\Omega)$ is continuous, and hence we have $\varphi \Pn \bw \to \varphi \bw$ strongly in $\Leb^{\infty}(Q)^d$. Together, this yields that
\begin{align}\ \label{ln-6v}
\ll \oU^{l,n,m} \otimes \varphi \Pn \bw,  \nabla \oU^{l,n,m} \rr_{Q}
\to
\ll \bu^m \otimes \varphi \bw, \nabla \bu^m \rr_{Q}\quad \text{ as }\, l,n \to \infty.
\end{align}
By \eqref{Pn-conv} we have in particular that $\varphi \BD \Pn \bw \to \varphi \BD \bw$ strongly in $\Leb^{q}(Q)^{d \times d}$ and by \eqref{cla-7} that $\overline{\BS}^{l,n,m} \wconv \BS^m$ weakly in $\Leb^{q'}(Q)^{d\times d}$. Thus, it follows that
\begin{align}\label{ln-7v}
\ll \overline{\BS}^{l,n,m}, \varphi \BD \Pn \bw \rr_{Q}
&\to
\ll \BS^m, \varphi \BD \bw \rr_{Q}, \quad \text{ as } l,n \to \infty.
\end{align}
Since $|\oU^{l,n,m}|^{2q'-2} \oU^{l,n,m} \rightharpoonup \abs{\bu^m}^{2q'-2}\bu^m$ weakly in $\Leb^{(2q')'}(Q)^d$ by \eqref{cla-6} and $\varphi \Pn \bw \to \varphi \bw$ in particular in $\Leb^{2q'}(Q)^d$, we obtain
\begin{align} \label{ln-8v}
\frac{1}{m}\ll \abs{\oU^{l,n,m}}^{2q'-2}\oU^{l,n,m} , \varphi  \Pn \bw \rr_{Q}
&\to
\frac{1}{m}\ll \abs{{\bu}^{m}}^{2q'-2}{\bU}^{m} , \varphi  \bw \rr_{Q}, \quad \text{ as } l,n \to \infty.
\end{align}
Finally, with the strong convergence $\overline{\bf}\to \bf$ in $\Leb^{q'}(0,T;\Sob^{-1,q'}(\Omega)^d)$ by  \eqref{of-conv} and with \eqref{Pn-conv} we have that
\begin{align} \label{ln-9v}
\ll \overline{\bf}, \varphi \Pn \bw \rr_{Q}
&\to
\ll \bf, \varphi \bw \rr_{Q}, \quad \text{ as } l,n \to \infty.
\end{align}
By the fact that $\bu^m$ is divergence-free, it follows that $\widetilde{b}(\bu^m, \bu^m, \varphi \bw) = b(\bu^m, \bu^m, \varphi \bw)$. So with $\mathfrak{L}^{l,n,m}$ and $\mathfrak{L}^{m}$ as defined in  \eqref{Llnm-t} and \eqref{Lm-t}, respectively, the convergence results \eqref{ln-5v}--\eqref{ln-9v} yield that
\begin{align} \label{ln-10v}
\ll \mathfrak{L}^{l,n,m}[\oU^{l,n,m},\Pn \bw], \varphi \rr_{(0,T)}
\to
\ll \mathfrak{L}^{m}[\bu^{m},\bw], \varphi \rr_{(0,T)}, \quad \text{ as } l,n \to \infty.
\end{align}
Now, from \eqref{ap-l-1-w}, using \eqref{ln-3v}, \eqref{ln-4v} and \eqref{ln-10v}, as $l,n \to \infty$, we have that
\begin{align}\label{ap-m-1-w}
- \ll \bu^m, \bw \partial_t \varphi \rr_{Q} =  \ll \bu_0, \varphi(0) \bw \rr_{\Omega}
+ \ll \mathfrak{L}^m[\bu^m, \bw],\varphi\rr_{(0,T)}
\end{align}
 for all $\bw \in C^{\infty}_{0,\diver}(\Omega)^d$ and all $\varphi \in C^{\infty}_0((-T,T))$.
\smallskip

\noindent
\textit{Step 2: Bound on the time-derivative}. The distributional derivative of $\bu^m$ satisfies, by definition and using \eqref{ap-m-1-w}, that
\begin{align}\label{ln-11v}
 \ll \partial_t \bu^m, \bw \varphi \rr_{Q} =
- \ll \bu^m, \bw \partial_t \varphi \rr_{Q} \overset{\eqref{ap-m-1-w}}{=}
 \ll \mathfrak{L}^m[\bu^m, \bw],\varphi\rr_{(0,T)}
\end{align}
 for all $\bw \in C^{\infty}_{0,\diver}(\Omega^d)$ and all $\varphi \in C^{\infty}_0((0,T))$, since $\supp \varphi \subset (0,T)$. Using this equation we wish to show that $\partial_t \bu^m \in \Leb^{\tau}(0,T;(\Xod)')$ (not uniformly in $m \in \mathbb{N}$), for $\tau$ as in \eqref{def:tau} and $\Xod$ as in  \eqref{def:Xq-s}.
For $\mathfrak{L}^m$ as defined in \eqref{Lm-t}, using the fact that $\bu^m \in \Leb^{2q'}(Q)^d$ and $\BS^m \in \Leb^{q'}(Q)^{d \times d}$ let us first estimate 
 \begin{align}\begin{split}\label{ln-12v}
\abs{\ll \mathfrak{L}^m[\bu^m, \bw],\varphi\rr_{(0,T)}} &\leq
\abs{\ll \bu^m \otimes \bu^m, \varphi \nabla \bw \rr_{Q}}
+\abs{\ll \BS^m, \varphi \BD \bw \rr_Q} \\
&\quad + \frac{1}{m} \abs{\ll \abs{\bu^m}^{2q'-2} \bu^m, \varphi \bw \rr_{Q}} + \abs{\ll \bf, \varphi \bw \rr_Q}\\
&\leq \norm{\bu^m}_{\Leb^{2q'}(Q)}^2 \norm{\varphi \nabla \bw}_{\Leb^{q}(Q)} + \norm{\BS^m}_{\Leb^{q'}(Q)} \norm{\varphi \BD \bw}_{\Leb^q(Q)} \\
&\quad + \frac{1}{m} \norm{\bu^m}_{\Leb^{2q'}(Q)}^{2q'-1} \norm{\varphi \bw}_{\Leb^{2q'}(Q)} \\
&\quad + \norm{\bf}_{\Leb^{q'}(0,T;\Sob^{-1,q'}(\Omega))} \norm{\varphi \bw}_{\Leb^{q}(0,T;\Sob^{1,q}(\Omega))}\\
&\leq c(m) \left(\norm{\varphi}_{\Leb^q(0,T)} + \norm{\varphi}_{\Leb^{2q'}(0,T)}\right) \left(\norm{\bw}_{\Sob^{1,q}(\Omega)} + \norm{\bw}_{\Leb^{2q'}(\Omega)}\right)\\
&\leq c(m)\norm{\varphi}_{\Leb^{\tau'}(0,T)} \norm{\bw}_{X^q(\Omega)}
\end{split}
 \end{align}
 for all $\varphi \in C^{\infty}_0((0,T))$ and all $\bw \in C^{\infty}_{0,\diver}(\Omega)^d$, since $\tau'= \max(2q',q)$.
By the density of the respective test function spaces, $\ll \mathfrak{L}^m[\bu^m,\cdot], \cdot\rr_{(0,T)}$ represents a bounded linear functional on $\Leb^{\tau'}(0,T;\Xod)$, and thus we have that $\partial_t \bu^m \in \Leb^{\tau}(0,T;(\Xod)')$ by \eqref{ln-11v} and by reflexivity of the function space.
Consequently, $\ll \partial_t \bu^m, \bw \rr_{\Omega}$ is integrable for $\bw \in C^{\infty}_{0,\diver}(\Omega)^d$, and thus, we can rephrase \eqref{ln-11v} by the fundamental lemma of calculus of variations in the pointwise sense in time, so \eqref{ap-m-1} is proved.
\smallskip

\noindent 
\textit{Step 3: Identification of the initial condition}. Let us first show that $\bu^m \in \Cw([0,T];\l2d)$. Recall that the embedding $\Xod\ctemb \l2d$ is dense and hence we have $( \l2d)'\ctemb (\Xod)'$.
 Furthermore, since the embedding $\l2d \ctemb \Leb^2(\Omega)^d$ is continuous and $\l2d$ is a closed subspace, one has (e.g., by \cite{R.1991}) that $(\l2d)' \simeq \Leb^2(\Omega)^d/(\l2d)^{\bot} $, and hence it follows that
\begin{align}\label{emb:L2}
\l2d \ctemb \Leb^2(\Omega)^d \ctemb (\l2d)'.
\end{align}
This implies with the above that the embedding
  $\l2d \ctemb (\Xod)'$ is continuous.

   Consequently, we have that $\bu^m \in \Leb^{\infty}(0,T;\l2d)\ctemb \Leb^1(0,T;(\Xod)')$. With this and the fact that in particular $\partial_t \bu \in \Leb^1(0,T;(\Xod)')$, Lemma~\ref{Lem:cont-1} implies that $\bu^m \in \Cw([0,T];(\Xod)')$. Furthermore, by this and the fact that $\bu^m \in \Leb^{\infty}(0,T;\l2d)$, again with $\l2d \ctemb (\Xod)'$, Lemma~\ref{Lem:cont-2} shows that $\bu^m \in \Cw([0,T];\l2d)$.

Next, we shall show that $\bu^m(0,\cdot)=\bu_0 \in \l2d$. Let $\varphi \in C^{\infty}_{0}((-T,T))$ such that $\varphi(0)=1$. Multiplying  \eqref{ap-m-1} by $\varphi$ and integrating over $(0,T)$, yields that
\begin{align}\label{ln-13v}
\ll \partial_t \bu^m,\bw \varphi \rr_{Q} = \ll \mathfrak{L}^m [\bu^m, \bw], \varphi \rr_{(0,T)}
\end{align}
for all $\bw \in C^{\infty}_{0,\diver}(\Omega)^d$.
On the other hand, by integration by parts, the fundamental theorem of calculus and the fact that $\bu^m \in \Cw([0,T];\l2d)$ and also applying  \eqref{ap-m-1-w}, we have
\begin{align}\begin{split}\label{ln-14v}
\ll \partial_t \bu^m,\bw \varphi \rr_{Q}
 &= \ll \partial_t ( \bu^m \varphi),\bw \rr_{Q}
- \ll \bu^m,\bw  \partial_t \varphi \rr_{Q} \\
&\overset{\eqref{ap-m-1-w}}{=} - \ll  \bu^m(0,\cdot) \varphi(0),\bw  \rr_{\Omega}
+ \ll \bu_0, \varphi(0) \bw \rr_{\Omega}
+ \ll \mathfrak{L}^m[\bu^m, \bw],\varphi\rr_{(0,T)}
\end{split}
\end{align}
for all $\bw \in C^{\infty}_{0,\diver}(\Omega)^d$.
Comparing with \eqref{ln-13v} and noting that $\varphi(0)=1$, we obtain
\begin{align}
 \ll  \bu^m(0,\cdot),\bw  \rr_{\Omega}
= \ll \bu_0, \bw \rr_{\Omega} \quad \text{ for all } \bw \in C^{\infty}_{0, \diver}(\Omega)^d.
\end{align}
Since $\bu_0,\bu^m(0,\cdot)\in \l2d$ are divergence-free, this suffices to conclude that $\bu_0 = \bu^m(0,\cdot)$.  

By \eqref{cla-1} we have strong convergence $\wU^{l,n,m}(s,\cdot) \to \bu^m(s, \cdot)$ in $\Leb^2(\Omega)^d$, as $l,n \to \infty$, for a.e. $s \in (0,T)$. Furthermore, by \eqref{cla-2} we have $\wU^{l,n,m}(0,\cdot) \to \bu_0$ strongly in $\Leb^2(\Omega)^d$, as $l,n \to \infty$, and consequently for a.e. $s \in (0,T)$ we obtain that
\begin{align}\label{ln-ic-1}
\norm{\bu^m(s, \cdot)-\bu_0}_{\Leb^2(\Omega)}^2
= \lim_{l,n \to \infty}
\norm{\wU^{l,n,m}(s, \cdot) - \wU^{l,n,m}(0, \cdot)}_{\Leb^2(\Omega)}^2.
\end{align}
By \eqref{en-ineq-ln-t}, for all $s \in (0,T)$ and $\L=(l,n,m) \in \mathbb{N}^3$ we have that
\begin{align}\label{ln-ic-2}
\norm{\wU^{\L}(s, \cdot)}^2 -  \norm{\wU^{\L}(0, \cdot)}^2 \leq 2 \ll \overline{\bf}, \oU^{\L} \rr_{Q_s}, \quad \text{ for } \L =(l,n,m) \in \mathbb{N}^3,
\end{align}
since the other terms can be shown to be nonnegative. Indeed, the nonnegativity of the term $\langle \overline{\BS}^{\L}, \BD \oU^{\L}\rangle_{Q_s}$ follows by the relation \eqref{SL-i}, the fact that $(\BD \oU^{\L}, \BS^{\L}) \in \mathcal{A}(\cdot)$ a.e. in $Q$ by  \eqref{impl-l}, and that $\mathcal{A}(\cdot)$ is monotone and $(\B0,\B0) \in \mathcal{A}(\cdot)$ a.e. in $Q$ by Assumption~\ref{assump-A}.
Expanding the norm on the right-hand side in  \eqref{ln-ic-1}, adding and subtracting twice the term $\|\wU^{l,n,m}(0,\cdot)\|_{\Leb^2(\Omega)}^2$ and applying \eqref{ln-ic-2}, for $\L=(l,n,m)$, we obtain that
 \begin{align}
\notag
 \norm{\wU^{\L}(s, \cdot)-\wU^{\L}(0, \cdot)}_{\Leb^2(\Omega)}^2
 &= \norm{\wU^{\L}(s, \cdot)}_{\Leb^2(\Omega)}^2
- 2 \ll \wU^{\L}(s, \cdot), \wU^{\L}(0, \cdot) \rr_{\Omega}
+ \norm{\wU^{\L}(0, \cdot)}_{\Leb^2(\Omega)}^2\\
\notag
&= \norm{\wU^{\L}(s, \cdot)}_{\Leb^2(\Omega)}^2
- \norm{\wU^{\L}(0, \cdot)}_{\Leb^2(\Omega)}^2\\
 \label{ln-ic-3}
&\quad + 2 \ll \wU^{\L}(0, \cdot) - \wU^{\L}(s, \cdot),  \wU^{\L}(0, \cdot)\rr_{\Omega}\\
\notag
&\overset{\eqref{ln-ic-2}}{\leq}
\ll \overline{\bf}, \oU^{\L} \rr_{Q_s}
+ 2 \ll \wU^{\L}(0, \cdot) - \wU^{\L}(s, \cdot),  \wU^{\L}(0, \cdot)\rr_{\Omega}.
 \end{align}
  Then, applying $\limsup_{l,n \to \infty}$ gives that, for a.e. $s \in (0,T)$,
\begin{align}
\notag
\limsup_{l,n\to \infty} \norm{\wU^{l,n,m}(s, \cdot)-\wU^{l,n,m}(0, \cdot) }
&\overset{\eqref{ln-ic-3}}{\leq}
\limsup_{l,n\to \infty} \ll \overline{\bf}, \oU^{l,n,m} \rr_{Q_s} \\ \notag
&\quad\quad + 2 \limsup_{l,n\to \infty}
 \ll \wU^{l,n,m}(0, \cdot) - \wU^{l,n,m}(s, \cdot),  \wU^{l,n,m}(0, \cdot)\rr_{\Omega}\\  \label{ln-ic-4}
&=  \ll \bf, \bu^{m} \rr_{Q_s} +
2 \ll \bu_0 - \bu^{m}(s, \cdot),  \bu_0 \rr_{\Omega}\\ \notag
&=  \ll \bf, \bu^{m} \rr_{Q_s} +
2 \ll \bu^m(0,\cdot) - \bu^{m}(s, \cdot),  \bu_0 \rr_{\Omega},
\end{align}
 since we have the convergence $\oU^{l,n,m}\wconv \bu^m$ weakly in $\Leb^{q}(0,T;\Sob^{1,q}_{0}(\Omega)^d)$ by
  \eqref{cla-4-v}, $\overline{\bf} \to \bf$ strongly in $\Leb^{q'}(0,T;\Sob^{-1,q'}(\Omega)^d)$ by \eqref{of-conv},
  $\wU^{l,n,m}(0,\cdot) \to \bu_0$ strongly in $\Leb^2(\Omega)^d$ by \eqref{cla-2} and $\wU^{l,n,m}(s, \cdot) \to \bu^m(s, \cdot)$ strongly in $\Leb^{2}(\Omega)^d$ for a.e. $s \in (0,T)$ by \eqref{cla-1}. In the final line we have used that $\bu^m(0,\cdot) = \bu_0$. Let us denote by $N^m \subset [0,T]$ the zero subset of times for which \eqref{ln-ic-4} does not hold.
 Applying \eqref{ln-ic-4} in \eqref{ln-ic-1} and taking $\liminf_{s \to 0_{+}}$ omitting the zero set $N^m$ we have that
 \begin{align} \begin{split}
 \label{ln-ic-5}
 0 &\leq \liminf_{ (0,T)\backslash N^m \ni s \to 0} \norm{\bu^m(s, \cdot)-\bu_0}_{\Leb^2(\Omega)}^2 \\
  &\overset{\eqref{ln-ic-1}}{\leq}
  \liminf_{ (0,T)\backslash N^m \ni s \to 0} \limsup_{l,n \to \infty}
\norm{\wU^{l,n,m}(s, \cdot) - \wU^{l,n,m}(0, \cdot)}_{\Leb^2(\Omega)}^2\\
& \overset{\eqref{ln-ic-4}}{\leq}
\liminf_{ (0,T)\backslash N^m \ni s \to 0}
 \ll \bf, \bu^{m} \rr_{Q_s} +
2 \ll \bu^m(0,\cdot) - \bu^{m}(s, \cdot),  \bu_0 \rr_{\Omega} = 0,
\end{split}
 \end{align}
 where the last equality follows from the absolute continuity of the integral and from the fact that $\bu^m \in \Cw([0,T];\l2d)$, $\bu_0 \in \l2d$ and $\l2d \ctemb (\l2d)'$. This shows \eqref{ap-m-2}.
\smallskip

\noindent
\textit{Step 4: Energy identity}. Recall that $\bu^m \in \Xqd\ctemb \Leb^{\min(q,2q')}(0,T;\Xod^d)$ and $\partial_t \bu^m \in \Leb^{\tau}(0,T;(X_{\diver}^q(\Omega))')$, where $\tau = \min(q',(2q')')$, and equation \eqref{ap-m-1} is satisfied.
Because of the lack of integrability in time an approximation procedure by means of mollification in time can be applied to show the  energy identity
\begin{align}\label{en-id-m}
\frac{1}{2} \norm{  \bu^m (s,\cdot) }_{\Leb^{2}(\Omega)}^2
+  \ll \BS^m, \BD \bu^{m} \rr_{Q_s}
+ \frac{1}{m} \norm{\bu^m}_{\Leb^{2q'}(Q_s)}^{2q'} =  \ll  \bf, \bu^{m} \rr_{Q_s} +  \frac{1}{2} \norm{\bu_0 }_{\Leb^{2}(\Omega)}^2
\end{align}
for a.e. $s \in (0,T)$, where also the attainment of the initial datum in the sense of \eqref{ap-m-2} is used. The proof follows by a standard procedure and we therefore omit the details; see, e.g., \cite[Ch.~2.5]{L.1969}.
\smallskip

\noindent
\textit{Step 5: Identification of the implicit relation}. Recall that we have by the assertion \eqref{impl-l} that the inclusion $(\BD \oU^{l,n,m}(\bz), \BS^{l,n,m}(\bz))\in \mathcal{A}(\bz)$ holds for a.e. $\bz \in Q$. Furthermore, by  \eqref{cla-4-v} we have that $\BD \oU^{l,n,m} \wconv \BD \bu^m$ weakly in $\Leb^{q}(Q)^{d \times d}$ and by  \eqref{cla-8} that $\BS^{l,n,m} \wconv \BS^m$  weakly in $\Leb^{q'}(Q)^{d \times d}$, as $l,n \to \infty$.
By Lemma~\ref{minty} it suffices to show that
\begin{align}\label{ln-impl-1}
\limsup_{l,n \to \infty} \ll \BS^{l,n,m}, \BD \oU^{l,n,m} \rr_{Q_s}
\leq  \ll \BS^m, \BD \bu^m \rr_{Q_s},
\end{align} in order to obtain $(\BD \bu^m(\bz),\BS^m(\bz))\in \mathcal{A}(\bz)$ for a.e. $\bz \in Q_s$. Then we can exhaust $Q$ by letting $s \to T$.
 We can only show \eqref{ln-impl-1} for a.e. $s \in (0,T)$ since the energy identity \eqref{en-id-m} is available only for a.e. $s \in (0,T)$ and some of the arguments used to show \eqref{ln-impl-1} are only available for a.e. $s \in (0,T)$.

Let us add and subtract the term $\langle \overline{\BS}^{l,n,m}, \BD \oU^{l,n,m}\rangle_{Q_s}$ to obtain
\begin{align}\begin{split}\label{ln-impl-2}
\ll \BS^{l,n,m}, \BD \oU^{l,n,m}\rr_{Q_s} &=  \ll \overline{\BS}^{l,n,m}, \BD \oU^{l,n,m}\rr_{Q_s}
+ \ll \BS^{l,n,m} - \overline{\BS}^{l,n,m}, \BD \oU^{l,n,m}\rr_{Q_s}
\\
& \eqqcolon \mathrm{I} + \mathrm{II} ,
\end{split}
 \end{align}
 where the first term appears in the equation  \eqref{ap-l-1} for the approximate solutions and the second term has to be shown to vanish.
The energy inequality \eqref{en-ineq-ln-t} yields that
\begin{align}\begin{split} \label{ln-impl-3}
\mathrm{I} = \ll \overline{\BS}^{l,n,m}, \BD \oU^{l,n,m} \rr_{Q_s} &\overset{\eqref{en-ineq-ln-t}}{\leq} - \frac{1}{2}\norm{\wU^{l,n,m}(s,\cdot)}_{\Leb^2(\Omega)}^2 + \frac{1}{2}\norm{\wU^{l,n,m}(0,\cdot)}_{\Leb^2(\Omega)}^2
\\
&\quad\quad\quad  
+
\ll \overline{\bf}, \oU^{l,n,m}\rr_{Q_s}
- \frac{1}{m} \norm{\oU^{l,n,m}}_{\Leb^{2q'}(Q_s)}^{2q'}.
\end{split}
\end{align}
For the second term in \eqref{ln-impl-2}, for $l \in \mathbb{N}$ let $j \in \{1, \ldots, l\}$ be such that $s \in (t_{j-1}, t_j]$, i.e.,  $j$ depends on $s$ and on $l$. By the relation \eqref{SL-i} we have that
\begin{align} \begin{split} \label{ln-impl-4}
\ll \BS^{l,n,m} - \overline{\BS}^{l,n,m}, \BD \oU^{l,n,m} \rr_{Q_{i-1}^i}
&= \ll \int_{t_{i-1}}^{t_i} \BS^{l,n,m}(t, \cdot) - \overline{\BS}^{l,n,m} \, \d t, \BD \oU^{l,n,m}\rr_{\Omega}  \\
&=
\ll \int_{t_{i-1}}^{t_i} \BS^{l,n,m}(t, \cdot) \, \d t - \delta_l \overline{\BS}^{l,n,m}, \BD \oU^{l,n,m}\rr_{\Omega} = 0,\end{split}
\end{align}
for any $i \in \{1, \ldots, l\}$. So for $\mathrm{II}$ we obtain that
\begin{align}\begin{split} \label{ln-impl-5}
\mathrm{II}
&= \ll \BS^{l,n,m} - \overline{\BS}^{l,n,m}, \BD \oU^{l,n,m}\rr_{Q_{t_{j}}}
-
\ll \BS^{l,n,m} - \overline{\BS}^{l,n,m}, \BD \oU^{l,n,m}\rr_{Q_s^{t_{j}}} \\
&\overset{\eqref{ln-impl-4}}{=} 0
-
\ll \BS^{l,n,m}, \BD \oU^{l,n,m}\rr_{Q_s^{t_{j}}}
+
\ll  \overline{\BS}^{l,n,m}, \BD \oU^{l,n,m}\rr_{Q_s^{t_{j}}}\\
&\leq \ll  \overline{\BS}^{l,n,m}, \BD \oU^{l,n,m}\rr_{Q_s^{t_{j}}},
\end{split}
\end{align}
where the inequality follows since $\mathcal{A}(\cdot)$ is monotone a.e. in $Q$, $(\B0,\B0) \in \mathcal{A}(\cdot)$ a.e. in $Q$ and by the fact that $(\BD \oU^{l,n,m}, \BS^{l,n,m}) \in \mathcal{A}(\cdot)$ a.e. in $Q$ by \eqref{impl-l}. For the remaining term first we use again \eqref{en-ineq-ln-t} on $(s, t_j)$, noting that the term involving $\frac{1}{m}$ is nonnegative, which yields
\begin{align}\begin{split} \label{ln-impl-6}
\ll  \overline{\BS}^{l,n,m}, \BD \oU^{l,n,m}\rr_{Q_s^{t_{j}}}
&\leq \ll \overline{\bf}, \oU^{l,n,m} \rr_{Q_s^{t_j}}  - \frac{1}{2} \norm{\wU^{l,n,m}(t_j,\cdot)}_{\Leb^2(\Omega)}^2 + \frac{1}{2} \norm{\wU^{l,n,m}(s,\cdot)}_{\Leb^2(\Omega)}^2.
\end{split}
\end{align}
By the duality of norms, the estimate \eqref{est:l-4} and by \eqref{of-bd} we obtain
\begin{align}\begin{split} \label{ln-impl-7}
 \ll \overline{\bf}, \oU^{l,n,m} \rr_{Q_s^{t_j}}
& \leq \norm{\overline{\bf}}_{\Leb^{q'}(s,t_j;\Sob^{-1,q'}(\Omega))} \norm{\oU^{l,n,m}}_{\Leb^q(0,T;\Sob^{1,q}(\Omega))} \\
 &\leq c(m) \norm{\bf}_{\Leb^{q'}(t_{j-1},t_j;\Sob^{-1,q'}(\Omega))}
 \leq c(m) \norm{\bf}_{\Leb^{q'}(s-\delta_l,s +\delta_l;\Sob^{-1,q'}(\Omega))}.
 \end{split}
 \end{align}
 Furthermore, we have $\wU^{l,n,m}(t_j,\cdot)=\oU^{l,n,m}(t_j,\cdot) = \oU^{l,n,m}(s,\cdot)$, since $s \in (t_{j-1},t_j]$, and hence
 \begin{align}\begin{split} \label{ln-impl-8}
 \mathrm{II} &\overset{\eqref{ln-impl-5},\eqref{ln-impl-6}}{\leq}
 \ll \overline{\bf}, \oU^{l,n,m} \rr_{Q_s^{t_j}}  - \frac{1}{2} \norm{\wU^{l,n,m}(t_j,\cdot)}_{\Leb^2(\Omega)}^2 + \frac{1}{2} \norm{\wU^{l,n,m}(s,\cdot)}_{\Leb^2(\Omega)}^2\\
 &\overset{\eqref{ln-impl-7}}{\leq}
 c(m) \norm{\bf}_{\Leb^{q'}(s-\delta_l,s +\delta_l;\Sob^{-1,q'}(\Omega))}
  - \frac{1}{2} \norm{\oU^{l,n,m}(s,\cdot)}_{\Leb^2(\Omega)}^2 + \frac{1}{2} \norm{\wU^{l,n,m}(s,\cdot)}_{\Leb^2(\Omega)}^2.
  \end{split}
\end{align}
Now applying $\limsup_{l,n \to \infty}$ to $(\mathrm{I}+\mathrm{II})$ with \eqref{ln-impl-3} and \eqref{ln-impl-8}, noting that the term involving $\wU^{l,n,m}(s,\cdot)$ drops out, we obtain
\begin{align}\begin{split}\label{ln-impl-9}
 \limsup_{l,n \to \infty} \; (\mathrm{I} + \mathrm{II})  &\leq
 - \lim_{l,n \to \infty} \frac{1}{2} \norm{\oU^{l,n,m}(s,\cdot)}_{\Leb^2(\Omega)}^2
 +  \lim_{l,n \to \infty} \frac{1}{2}\norm{\wU^{l,n,m}(0,\cdot)}_{\Leb^2(\Omega)}^2 \\
&\quad  -  \frac{1}{m} \liminf_{l,n \to \infty} \norm{\oU^{l,n,m}}_{\Leb^{2q'}(Q_s)}^{2q'}
+   \lim_{l,n \to \infty}
\ll \overline{\bf}, \oU^{l,n,m}\rr_{Q_s} \\
&\quad
+
 c(m)  \lim_{l \to \infty} \norm{\bf}_{\Leb^{q'}(s-\delta_l,s +\delta_l;\Sob^{-1,q'}(\Omega))}\\
 &\leq
 - \frac{1}{2} \norm{\bu^{m}(s,\cdot)}_{\Leb^2(\Omega)}^2  + \frac{1}{2}\norm{\bu_0}_{\Leb^2(\Omega)}^2
 - \frac{1}{m}\norm{\bu^m}_{\Leb^{2q'}(Q_s)}^{2q'} + \ll \bf, \bu^{m}\rr_{Q_s},
 \end{split}
\end{align}
where the last inequality is based on the following arguments.
By \eqref{cla-9} we have that $\oU^{l,n,m}(s,\cdot) \to \bu^m(s, \cdot)$ strongly in $\Leb^2(\Omega)^d$, as $l,n \to \infty$, for a.e. $s \in (0,T)$. The second term converges to $\frac{1}{2}\|\bu_0\|_{\Leb^2(\Omega)}^2$, since by  \eqref{cla-2} we have that $\wU^{l,n,m}(0,\cdot) \to \bu_0$ strongly in $\Leb^2(\Omega)^d$. For the third term we use weak lower-semicontinuity with respect to the weak convergence in $\Leb^{2q'}(Q_s)^d$ and  \eqref{cla-4-v}. For the forth term we have convergence, since $\oU^{l,n,m} \wconv \bu^m$ weakly in $\Leb^{q}(0,T;\Sob^{1,q}(\Omega)^d)$ by  \eqref{cla-4-v} and $\overline{\bf} \to \bf$ strongly in $\Leb^{q'}(0,T;\Sob^{-1,q'}(\Omega)^d)$ by  \eqref{of-conv}, as $l,n \to \infty$. The last term vanishes by the absolute continuity of the integral, as $l \to \infty$.
Finally returning to \eqref{ln-impl-2}, applying $\limsup_{l,n \to \infty}$ and the energy identity \eqref{en-id-m} for a.e. $s \in (0,T)$, yields
\begin{align}\begin{split}\label{ln-impl-10}
 \limsup_{l,n \to \infty}  \ll \BS^{l,n,m}, \BD \oU^{l,n,m}\rr_{Q_s}
 &\overset{\eqref{ln-impl-2}}{\leq}
 \limsup_{l,n \to \infty} \; (\mathrm{I} + \mathrm{II})  \\
 & \overset{\eqref{ln-impl-9}}{\leq }
  - \frac{1}{2} \norm{\bu^{m}(s,\cdot)}_{\Leb^2(\Omega)}^2  + \frac{1}{2}\norm{\bu_0}_{\Leb^2(\Omega)}^2 \\
 &\qquad\quad - \frac{1}{m}\norm{\bu^m}_{\Leb^{2q'}(Q_s)}^{2q'} + \ll \bf, \bu^{m}\rr_{Q_s}  \\
 &\overset{\eqref{en-id-m}}{=} \ll \BS^m, \BD \bu^m\rr_{Q_s}
 \end{split}
\end{align}
for a.e. $s \in (0,T)$. This proves the claim in  \eqref{ln-impl-1} and completes the proof.
\end{proof}
\end{Lemma}

\subsection*{Limit $m \to \infty$}  In this step we loose the admissibility of the solution as a test function, and we have to use Lipschitz truncation to identify the implicit relation. The availability of the solenoidal Lipschitz truncation allows to simplify the arguments in \cite{BGMS.2012}, since no pressure has to be reconstructed.

Let us denote
\begin{align}\label{def:mu}
\mu &\coloneqq \min \left(\frac{q(d+2)}{2d},q',(2q')' \right) = \min\left(\wq',\tau\right),
\end{align}
where $\tau$ defined in \eqref{def:tau}. Note that  since $q >\frac{2d}{d+2}$, we have that $\mu>1$.
	
\begin{Lemma}[Convergence $m \to \infty$]\label{Lem:m-conv}
Let $\bu^{m} \in \Leb^{\infty}(0,T;\l2d) \cap \Xqd
$ be such that  $\partial_t \bu^m \in \Leb^{\tau}(0,T;(\Xod)')$ and let $\BS^m\in L^{q'}(Q)^{d\times d}$ be a solution to \eqref{ap-m-1}--\eqref{ap-m-2}, for $m \in \mathbb{N}$. Further, let $\mu>1$ and $\wq$ be defined in \eqref{def:mu} and \eqref{def:hq}, respectively. Then, there exists a constant $c>0$ such that we have
 \begin{align}\begin{split}\label{est:m}
 \norm{\bu^{m}}^2_{\Leb^{\infty}(0,T;\Leb^2(\Omega))}
 &+ \norm{\bu^{m}}_{\Leb^q(0,T;\Sob^{1,q}(\Omega))}^q
 +
\norm{\BS^{m}}_{\Leb^{q'}(Q)}^{q'}\\
 &+ \frac{1}{m}\norm{\bu^{m}}_{\Leb^{2q'}(Q)}^{2q'} + \norm{\bu^{m}}_{\Leb^{\frac{q(d+2)}{d}}(Q)} \leq c \quad \text{ for all } m \in \mathbb{N}.
 \end{split}
 \end{align}
  Furthermore, there exists a function $\bu \in \Leb^{\infty}(0,T;\l2d) \cap \Leb^q(0,T;\Sob^{1,q}_{0, \diver}(\Omega)^d)
  $ such that $\partial_t \bu \in \Leb^{\mu}(0,T; (\Sob^{1,\wq}_{0,\diver}(\Omega)^d)')$, an $\BS \in \Leb^{q'}(Q)^{d \times d}$ and subsequences such that, as $m \to \infty$,
 \begin{align}  \label{cm-0}
\bu^{m} &\to \bu  \; && \text{ strongly in  }\Leb^q(0,T;\l2d)\cap  \Leb^{r}(Q)^d, \, \forall r\in[1,\tfrac{q(d+2)}{d}),\\
\label{cm-0-v}
\bu^{m}(s,\cdot) &\to \bu(s,\cdot)  \; && \text{ strongly in  }\l2d \, \text{ for a.e. } s \in (0,T),\\
 \label{cm-1}
 \bu^{m} &\rightharpoonup \bu
 \; && \text{ weakly in  }\Leb^q(0,T;\Sob^{1,q}_{0,\diver}(\Omega)^d)
 \cap \Leb^{\frac{q(d+2)}{d}}(Q)^d, \\
 \label{cm-2}
 \bu^{m} & \overset{*}{\rightharpoonup} \bu
 \; && \text{ weakly* in  } \Leb^{\infty}(0,T;\l2d),\\
 \label{cm-3}
 \partial_t \bu^{m} &\rightharpoonup \partial_t \bu \; && \text{ weakly in } \Leb^{\mu}(0,T; (\Sob^{1,\wq}_{0,\diver}(\Omega)^d)'),\\
 \label{cm-4}
  \BS^{m} &\rightharpoonup \BS \; && \text{ weakly in } \Leb^{q'}(Q)^{d\times d},\\
   \label{cm-5}
\frac{1}{m}\abs{\bu^{m}}^{2q'-2}\bu^{m} & \to \b0  &&  \text{ strongly in } \Leb^{p}(Q)^{d\times d}, \; \text{ for any } p \in [1, (2q')').
 \end{align}
 \begin{proof}\hspace*{\fill} \smallskip

\noindent
\textit{Step 1: Estimates}. Recall that by \eqref{en-id-m} we have the following energy identity for a.e. $t \in (0,T)$:
\begin{align}\label{en-id-m-1}
\frac{1}{2} \norm{  \bu^m (t,\cdot) }_{\Leb^{2}(\Omega)}^2
+  \ll \BS^m, \BD \bu^{m} \rr_{Q_t}
+ \frac{1}{m} \norm{\bu^m}_{\Leb^{2q'}(Q_t)}^{2q'} =  \ll  \bf, \bu^{m} \rr_{Q_t} +  \frac{1}{2} \norm{\bu_0 }_{\Leb^{2}(\Omega)}^2.
\end{align}
By the fact that $(\BD \bu^m, \BS^m)\in \mathcal{A}(\cdot)$ a.e. in $Q$ by \eqref{impl-m}, we can use Assumption~\ref{assump-A} \ref{itm:A4} to show that
\begin{align}\begin{split}\label{est:m-1}
\ll \BS^m, \BD \bu^m \rr_{Q_t} &\geq - \norm{g}_{\Leb^1(Q_t)} + c_*\left(\norm{\BD \bu^m}_{\Leb^q(Q_t)}^q +  \norm{ \BS^m}_{\Leb^{q'}(Q_t)}^{q'} \right) \\
&\geq - \norm{g}_{\Leb^1(Q_t)} + c \norm{ \bu^m}_{\Leb^q(0,t;\Sob^{1,q}(Q))}^q + c_* \norm{ \BS^m}_{\Leb^{q'}(Q_t)}^{q'},
\end{split}
\end{align}
where we have used Poincar\'e's and Korn's inequalities in the last line.
Similarly as before, we use duality of norms and Young's inequality with $\varepsilon>0$ to bound
\begin{align}\begin{split}\label{est:m-2}
\ll \bf,  \bu^m \rr_{Q_t}
\leq c(\varepsilon) \norm{\bf}_{\Leb^{q'}(0,T;\Sob^{-1,q'}(\Omega))}^{q'} + \varepsilon
\norm{\bu^m}_{\Leb^{q}(0,t;\Sob^{1,q}(\Omega))}^q.
\end{split}
\end{align}
Applying \eqref{est:m-1} and \eqref{est:m-2} in  \eqref{en-id-m-1}, rearranging and choosing $\varepsilon>0$ small enough yields
\begin{align}\begin{split}\label{est:m-3}
 \norm{  \bu^m (t,\cdot) }_{\Leb^{2}(\Omega)}^2
+ \norm{ \bu^{m}}_{\Leb^{q}(0,t;\Sob^{1,q}(\Omega))}^q + \norm{\BS^m}_{\Leb^{q'}(Q_t)}^{q'}
+ \frac{1}{m} \norm{\bu^m}_{\Leb^{2q'}(Q_t)}^{2q'}
\leq c
\end{split}
\end{align}
for a.e. $t\in (0,T)$ and all $m \in \mathbb{N}$. Taking the essential supremum over $t\in (0,T)$ and also applying the parabolic interpolation from Lemma~\ref{par-interpol} shows \eqref{est:m}.
\smallskip

\noindent
\textit{Step 2: Bound on the time derivative}. In order to derive a uniform bound on the time derivative let us estimate $\mathfrak{L}^m[\bu^m;\bv]$. Since no uniform bounds on  $\norm{\bu^m}_{\Leb^{2q'}(Q)}$ are available at this point,
we use the bound \eqref{est:b-3} on the convective term, with $\wq$ as defined in \eqref{def:hq}, to deduce that
\begin{align*}
\abs{\ll \bu(t, \cdot) \otimes \bu(t, \cdot), \nabla \bv \rr_{\Omega}}  \leq c \norm{\bu(t, \cdot)}_{\Leb^{\frac{q(d+2)}{d}}(\Omega)}^2 \norm{\nabla \bv}_{\Leb^{\wq}(\Omega)},
\end{align*}
 which holds, since $q \geq \frac{2d}{d+2}$.
 Note that the embedding $\Sob^{1,\wq}(\Omega) \ctemb \Sob^{1,q}(\Omega)\cap \Leb^{2q'}(\Omega)$ is continuous for $\wq$ as in \eqref{def:hq}.
Also we have that $\mu'=\max\left(2q',q,\left(\tfrac{q(d+2)}{2d}\right)'\right)=\max(2q',\wq)$ for $\mu$ as in  \eqref{def:mu}, i.e., the embedding $\Leb^{\mu'}(\Omega)\ctemb \Leb^{\wq}(\Omega)\cap \Leb^{q}(\Omega)\cap \Leb^{2q'}(\Omega)$ is continuous.
With this, similarly as in \eqref{ln-12v}
 applying the uniform estimates in \eqref{est:m}
  one has that
 \begin{align}\begin{split}\label{est:m-4}
\abs{\ll \mathfrak{L}^m[\bu^m, \bw],\varphi\rr_{(0,T)}}
&\leq \norm{\bu^m}_{\Leb^{\frac{q(d+2)}{d}}(Q)}^2 \norm{\varphi \nabla \bw}_{\Leb^{\wq}(Q)} + \norm{\BS^m}_{\Leb^{q'}(Q)} \norm{\varphi \BD \bw}_{\Leb^q(Q)} \\
&\quad + \frac{1}{m} \norm{\bu^m}_{\Leb^{2q'}(Q)}^{2q'-1} \norm{\varphi \bw}_{\Leb^{2q'}(Q)} \\
&\quad + \norm{\bf}_{\Leb^{q'}(0,T;\Sob^{-1,q'}(\Omega))} \norm{\varphi \bw}_{\Leb^{q}(0,T;\Sob^{1,q}(\Omega))}\\
&\leq c\norm{\varphi}_{\Leb^{\mu'}(0,T)} \norm{\bw}_{\Sob^{1,\wq}(\Omega)}
\end{split}
 \end{align}
for all $\varphi \in C^{\infty}_0((0,T))$ and all $\bw \in C^{\infty}_{0,\diver}(\Omega)^d$, and all $m \in \mathbb{N}$.
 With \eqref{ap-m-1} and using the fact that $\mu>1$ and that the space $\Leb^{\mu}(0,T;(\Sob^{1,\wq}_{0,\diver}(\Omega)^d)')$ is reflexive, this shows that
 $\{\partial_t \bu^m\}_{m \in \mathbb{N}}$ is bounded in $\Leb^{\mu}(0,T;(\Sob^{1,\wq}_{0,\diver}(\Omega)^d)')$.
\smallskip

\noindent
\textit{Step 3: Convergence as $m \to \infty$}. Since $q >\frac{2d}{d+2}$, the embedding $\Sob^{1,q}_{0,\diver}(\Omega)^d \cpemb \l2d$ is compact. Because  $\wq\geq q >\frac{2d}{d+2}$, the embedding $\Sob^{1,\wq}_{0,\diver}(\Omega)\ctemb \l2d$ is in particular continuous and dense, which implies that
$ (\l2d) \ctemb (\Sob^{1,\wq}_{0,\diver}(\Omega))'$. Combined with the embedding in \eqref{emb:L2}, this yields that the embedding $\l2d \ctemb (\Sob^{1,\wq}_{0,\diver}(\Omega))'$ is continuous.
Hence, the Aubin--Lions compactness lemma implies that the embedding
\begin{align*}
\{\bv \in \Leb^q(0,T; \Sob^{1,q}_{0,\diver}(\Omega)^d)\colon \, \partial_t \bv\in \Leb^{\mu}(0,T; (\Sob^{1,\wq}_{0,\diver}(\Omega))') \} \cpemb \Leb^{q}(0,T, \l2d)
\end{align*}
is compact, see for example \cite[Lem.~7.7]{Ro.2013}.
The fact that by \eqref{est:m} the sequence $\{\bu^{m}\}_{m \in \mathbb{N}}$ is bounded in $\Leb^q(0,T;\Sob^{1,q}_{0,\diver}(\Omega)^d)$ and
 that $\{\partial_t \bu^{m}\}_{m \in \mathbb{N}}$ is bounded in $\Leb^{\mu}(0,T;(\Sob^{1,\wq}_{0,\diver}(\Omega)^d)')$ then ensures the existence of a subsequence such that
\begin{alignat}{2}\label{conv-m-0}
 \bu^{m} & \to \bu \quad && \text{ strongly in  } \Leb^{q}(0,T;\l2d), \quad \text{ as } m \to \infty.
\end{alignat}
By the estimates in \eqref{est:m}, the uniform bound on $\{\partial_t \bu^m\}_{m \in \mathbb{N}}$  in $\Leb^{\mu}(0,T;(\Sob^{1,\wq}_{0,\diver}(\Omega))')$ and the (sequential) Banach--Alaoglu theorem, there exists a subsequence such that \eqref{cm-1}--\eqref{cm-4} holds, where the limits can be identified with the help of \eqref{conv-m-0}.

The strong convergence in $\Leb^r(Q)^d$ for all $r \in [1, \frac{q(d+2)}{d})$ asserted in \eqref{cm-0}
 follows from the strong convergence in $\Leb^{1}(Q)^d$ by \eqref{conv-m-0}, and the boundedness in $\Leb^{\frac{q(d+2)}{d}}(Q)^d$ by \eqref{est:m} by means of interpolation. The convergence \eqref{cm-0-v} is deduced analogously to the proof of \eqref{cla-1} by the arguments following \eqref{conv-l-4}. With H\"older's inequality and the estimate in \eqref{est:m} we find that
\begin{align}\begin{split}\label{conv-m-1}
\norm{\frac{1}{m} \abs{\bu^m}^{2q'-2}\bu^m}_{\Leb^1(Q)}
 \leq \frac{c}{m} \norm{\bu^m}_{\Leb^{2q'}(Q)}^{2q'-1}
 \overset{\eqref{est:m}}{\leq}
 c m^{-\sfrac{1}{2q'}} \to 0, \quad \text{ as } m \to \infty,
 \end{split}
\end{align}
so strong convergence to $\b0$ of the regularization term in $\Leb^1(Q)^d$ is proved.
One can show uniform boundedness of $\frac{1}{m}\abs{\bu^m}^{2q'-2}\bu^m$ in $\Leb^{(2q')'}(Q)^d$; interpolation between $\Leb^1(Q)^d$ and $\Leb^{(2q')'}(Q)^d$
 then gives strong convergence to $\mathbf{0}$ in $\Leb^p(Q)^d$, for any $p \in [1, (2q')')$. Hence \eqref{cm-5}  follows.
\end{proof}
\end{Lemma}

For $t \in (0,T)$, $\bu \in \Leb^{\frac{q(d+2)}{d}}(Q)^d$ and $\bv \in \Sob^{1,\wq}_0(\Omega)^d$ with $\wq$ defined in \eqref{def:hq}, let us introduce
\begin{align}\begin{split}\label{L-n}
\mathfrak{L}[\bu; \bv](t)
\coloneqq&
- b(\bu(t, \cdot),\bu(t, \cdot), \bv) - \ll \BS(t,\cdot), \BD\bv\rr_{\Omega}
+ \ll \bf(t, \cdot), \bv\rr_{\Omega},
\end{split}
\end{align}
where $\BS \in \Leb^{q'}(Q)^{d \times d}$ is the limiting function introduced in Lemma~\ref{Lem:m-conv}.
	
\begin{Lemma}[Identification of the PDE as $m \to \infty$]\label{Lem:m-id} 
The limiting function $\bu \in \Leb^{\infty}(0,T;\l2d)\cap \Leb^{q}(0,T;\Sob^{1,q}_{0,\diver}(\Omega)^d)$ from Lemma~\ref{Lem:m-conv} satisfies that $\partial_t \bu \in \Leb^{\wq'}(0,T;(\Sob^{1,\wq}_{0,\diver}(\Omega)^d)')$, with $\wq$ defined in \eqref{def:hq}. (Up to a representative) we have that $\bu \in \Cw([0,T],\l2d)$. Furthermore, the functions $\bu$ and $\BS \in \Leb^{q'}(Q)^{d \times d}$ from Lemma~\ref{Lem:m-conv} satisfy
\begin{align}\label{p-1}
\ll \partial_t \bu(t,\cdot), \bw\rr_{\Omega}
&=  \mathfrak{L}[\bu;\bw](t)  \; &&\text{ for all } \bw \in C^{\infty}_{0, \diver}(\Omega)^d \text{ for a.e. } t \in (0,T), \\
(\BD \bu(\bz), \BS(\bz)) &\in \mathcal{A}(\bz)   \; &&\text{ for a.e. } \bz \in Q, \label{impl}\\
 \esslim_{t \to 0_+ } \norm{\bu(t, \cdot) - \bu_0}_{\Leb^2(\Omega)} &  \,=  0, && \label{p-2}
\end{align}
i.e., $(\bu, \BS)$ is a weak solution according to Definition~\ref{def:w-sol}.
\begin{proof}\hspace*{\fill}
\smallskip

\noindent 
\textit{Step 1: Identification of the limiting equation}. Let $\bw \in C^{\infty}_{0,\diver}(\Omega)^d$ and $\varphi \in C^{\infty}_0((0,T))$ and let us consider each of the terms in \eqref{ap-m-1} and \eqref{p-1}. By the weak convergence results in \eqref{cm-3} and \eqref{cm-4} we have that
\begin{align}\label{lm-1}
\ll \partial_t \bu^m, \varphi \bw \rr_{Q} &\to \ll \partial_t \bu, \varphi \bw \rr_{Q},\\\label{lm-2}
\ll  \BS^m, \varphi \BD \bw \rr_{Q} &\to \ll \BS, \varphi \BD \bw \rr_{Q},
\end{align}
as $m \to \infty$. Since by \eqref{cm-0} we have that $\bu^m \to \bu$ in $\Leb^{r}(Q)^{d \times d}$ for all $r \in \left[1, \frac{q(d+2)}{d}\right)$ it follows that $\bu^m \otimes \bu^m \to \bu \otimes \bu$ in $\Leb^{r}(Q)^{d \times d}$ for all $r \in \left[1, \frac{q(d+2)}{2d}\right)$. Since $q>\frac{2d}{d+2}$, this set is nonempty and the convergence holds in particular in $\Leb^1(Q)$, so
\begin{align}\label{lm-3}
\ll \bu^m \otimes \bu^m, \varphi \nabla \bw \rr_{Q} &\to \ll \bu \otimes \bu, \varphi \nabla \bw \rr_{Q}, \quad \text{ as } m \to \infty.
\end{align}
Taking the results in \eqref{lm-1}--\eqref{lm-3} and  \eqref{cm-5} shows that \eqref{ap-m-1} implies \eqref{p-1}.
\smallskip

\noindent
\textit{Step 2: Identification of the initial condition}. With similar arguments as in the proof of Lemma~\ref{Lem:ln-id}, Step 3, it follows that $\bu \in \Cw([0,T];\l2d)$, that $\bu_0 = \bu(0,\cdot) \in \l2d$ and that the initial datum is attained in the sense of  \eqref{p-2}.
\smallskip

\noindent
\textit{Step 3: Higher integrability of the time derivative}. As in Step 2 in the proof of Lemma~\ref{Lem:ln-conv} we can improve the integrability of $\partial_t \bu$ using the fact that \eqref{p-1} is satisfied.
This yields that $\partial_t \bu \in \Leb^{\wq'}(0,T;(\Sob^{1,\wq}_{0,\diver}(\Omega)^d)')$, for $\wq$ as defined in \eqref{def:hq}.
\smallskip

\noindent
\textit{Step 4: Identification of the implicit relation (compare \cite{BGMS.2012} and \cite[Sec.~3]{BDS.2013})}.
Recall that $\BD \bu^m \rightharpoonup \BD \bu$ weakly in $\Leb^{q}(Q)^{d \times d}$ by \eqref{cm-1}, that $\BS^m \rightharpoonup \BS $ weakly in $\Leb^{q'}(Q)^{d \times d}$ by \eqref{cm-4} and that we have that $(\BD \bu^m(\bz), \BS^m(\bz)) \in \mathcal{A}(\bz)$ for a.e. $\bz \in Q$ by \eqref{impl-m}. Hence, by Lemma~\ref{minty}, it suffices to show that
\begin{align}\label{lm-impl-1}
\limsup_{m \to \infty} \ll \BS^{m}, \BD \bu^{m} \rr_{\widetilde{Q}}
\leq  \ll \BS, \BD \bu \rr_{\widetilde{Q}},
\end{align}
for a set $\widetilde{Q}\subset Q$,
 to identify the implicit relation $(\BD \bu, \BS) \in \mathcal{A}(\cdot)$ a.e. on $\widetilde{Q}$.

Since there is no energy identity available for $\bu$, in order to identify the implicit relation one has to truncate the elements of the approximating sequence of velocity fields suitably so as to be able to use them as test functions. In contrast with \cite{BGMS.2012} we will not use a parabolic Lipschitz truncation after locally reconstructing the approximations to the pressure, but work with the solenoidal Lipschitz truncation introduced subsequently in \cite{BDS.2013} and stated in Lemma~\ref{Lem:LA-par-div}, as the argument is then more direct.

We wish to truncate $\bv^m \coloneqq \bu^m - \bu$, which satisfies, for all $\bxi \in C^{\infty}_{0,\diver}(Q)^d$, the equality
\begin{align}\label{lm-impl-2}
\ll \partial_t \bv^m, \bxi \rr_{Q} = \ll \bu^m \otimes \bu^m - \bu \otimes \bu, \nabla \bxi \rr_{Q} - \ll \BS^m - \BS, \BD \bxi \rr_{Q} - \frac{1}{m}\ll \abs{\bu^m}^{2q'-2}\bu^m, \bxi \rr_{Q},
\end{align}
by \eqref{ap-m-1} and \eqref{p-1} and by the density of $C^{\infty}_0(0,T)\times C^{\infty}_{0,\diver}(\Omega)^d$ in $C^{\infty}_{0,\diver}(Q)^d$.
In order to rewrite the equation in divergence form as required for Lemma~\ref{Lem:LA-par-div}, we adapt the last term locally.

Let $B_0 \subset \subset \Omega$ be an arbitrary but fixed ball compactly contained in $\Omega$. For a.e. $t \in (0,T)$ we seek a weak solution $\bg^m_0(t, \cdot)$ to
\begin{align}\begin{split}\label{lm-impl-3}
\Delta \bg^m_0(t, \cdot) &= \frac{1}{m} \abs{\bu^m(t,\cdot)}^{2q'-2}\bu^m(t,\cdot) \quad \text{ in } B_0,\\
\bg^m_0(t, \cdot) |_{\partial B_0} &= \b0,\end{split}
\end{align}
i.e., for suitable $p \in (1, \infty)$ we wish to find a $\bg^m_0(t,\cdot)\in \Sob^{1,p}_0(B_0)^d$ such that
\begin{align}\label{lm-impl-4}
\ll \nabla \bg^m_0(t, \cdot) , \nabla \bv \rr_{B_0}&= \frac{1}{m} \ll \abs{\bu^m(t,\cdot)}^{2q'-2}\bu^m(t,\cdot), \bv \rr_{B_0} \quad \text{ for all } \bv \in C^{\infty}_0(B_0)^d.
\end{align}
As $\partial B_0$ is smooth, standard regularity theory for Poisson's equation (see, \cite[Thm.~2.4.2.5]{Gr.2011} and \cite[Lem.~9.17]{GT.2001}) guarantees the existence of a unique $\bg^m_0(t,\cdot)\in \Sob^{2,p}(B_0)^d\cap \Sob^{1,p}_0(B_0)^d$ such that
\begin{align}\label{lm-impl-5}
\norm{\bg^m_0(t,\cdot)}_{\Sob^{2,p}(B_0)} \leq c \norm{\frac{1}{m} \abs{\bu^m(t,\cdot)}^{2q'-2}\bu^m(t,\cdot)}_{\Leb^{p}(B_0)}, \quad \text{ for } p \in (1, \infty).
\end{align}
For an arbitrary but fixed interval $I_0 \subset \subset (0,T)$, viewing $\bg^m_0$ as a function defined on 
$Q_0 \coloneqq I_0 \times B_0$ by \eqref{lm-impl-5}  one has that
\begin{align}\label{lm-impl-6}
\norm{\bg^m_0}_{\Leb^p(I_0;\Sob^{2,p}(B_0))} \leq c \norm{\frac{1}{m} \abs{\bu^m}^{2q'-2}\bu^m}_{\Leb^{p}(Q_0)},
\end{align}
and by \eqref{cm-5} the right-hand side converges, up to a subsequence, to zero, as $m \to \infty$, for $p \in [1, (2q')')$. This implies in particular that
\begin{align}\label{lm-impl-7}
\bg^m_0 \to \b0 \quad \text{ strongly in }  \Leb^p(I_0;\Sob^{1,p}(B_0)^d), \quad \text{ as } m \to \infty
\end{align}
for all $p \in [1, (2q')')$.
Furthermore, using \eqref{lm-impl-4} in \eqref{lm-impl-2} we have for the cylinder $Q_0 \subset \subset Q$ that
\begin{align}\label{lm-impl-8}
\ll \partial_t \bv^m, \bxi \rr_{Q_0} = \ll \bu^m \otimes \bu^m - \bu \otimes \bu, \nabla \bxi \rr_{Q_0} - \ll \BS^m - \BS, \BD \bxi \rr_{Q_0} - \ll \nabla \bg^m_0, \nabla \bxi \rr_{Q_0}
\end{align}
for all $\bxi \in C^{\infty}_{0,\diver}(Q_0)^d$.

Now we wish to apply Lemma~\ref{Lem:LA-par-div} with $p=q \in (1,\infty)$ and $\sigma$ such that
\begin{align}\label{lm-impl-9}
1 < \sigma < \min\left(2,q,q',\frac{q(d+2)}{2d}, (2q')'\right) = \min\left(q',\frac{q(d+2)}{2d}, (2q')'\right).
\end{align}
Such a $\sigma$ exists, since we have by assumption that $q > \frac{2d}{d+2}$. First note that $\bu$ and $\bu^m$ are (weakly) divergence-free, and so is $\bv^m$, and $\bv^m \rightharpoonup \b0$ weakly in $\Leb^{q}(I_0;\Sob^{1,q}(B_0)^d)$, as $m \to \infty$
by \eqref{cm-1}.

Since $\bu^m \to \bu$ strongly in $\Leb^{p}(Q)$ for $p \in [1, \frac{q(d+2)}{d})$ by \eqref{cm-0} and $\sigma < \frac{q(d+2)}{d}$ we have that $\bv^m \to \b0$ strongly in $\Leb^{\sigma}(Q_0)^d$, as $m \to \infty$. Furthermore, since $\{\bu^m\}_{m \in \mathbb{N}}$ is bounded in $\Leb^{\infty}(0,T;\Leb^2(\Omega)^d)$ by \eqref{est:m}  we have with $\sigma < 2$ that $\{\bv^m\}_{m \in \mathbb{N}}$ is bounded in $\Leb^{\infty}(0,T;\Leb^{\sigma}(\Omega)^d)$.
Now we set
\begin{align*}
\BG^m_1 \coloneqq \BS- \BS^m,\quad \text{ and } \quad
\BG^m_2 \coloneqq \bu^m \otimes \bu^m - \bu \otimes \bu - \nabla \bg^m_0.
\end{align*}
Note that $\BG^m_1 \rightharpoonup \B0$ weakly in $\Leb^{q'}(Q_0)^{d \times d}$ by \eqref{cm-4}.

By \eqref{cm-0} we have that $\bu^m \to \bu$ in $\Leb^{r}(Q)^d$ for all $r \in \left[1, \frac{q(d+2)}{d}\right)$, and thus, $\bu^m \otimes \bu^m \to \bu \otimes \bu$ in $\Leb^{r}(Q)^{d\times d}$ for all $r \in \left[1, \frac{q(d+2)}{2d}\right)$. This holds in particular for $r=\sigma<\frac{q(d+2)}{2d}$.
Furthermore, by \eqref{lm-impl-7} we have that $\nabla \bg^m_0 \to \B0$ strongly in $\Leb^{p}(Q_0)^{d \times d}$, for any $p \in [1, (2q')' )$, and hence also for $p=\sigma < (2q')'$.
This means that all the assumptions of Lemma~\ref{Lem:LA-par-div} are satisfied.

With the aid of the parabolic solenoidal Lipschitz truncation we show that
\begin{align}\label{lm-impl-10}
\lim_{m \to \infty} \int_{\tfrac{1}{8}Q_0} \left[\left( \BS^{m}-\BS^{\star}(\cdot, \BD \bu\right) \colon  \left( \BD \bu^{m} - \BD \bu\right)\right]^{\sfrac{1}{2}} \, \d \bz = 0,
\end{align}
where the exponent $\sfrac{1}{2}$ is used to control the size of the set where $\bv^m = \bu^m - \bu$ and its truncation do not coincide.
By the monotonicity of $\mathcal{A}$ and the fact that $(\BD \bu, \BS^{\star}(\cdot, \BD \bu)) \in \mathcal{A}(\cdot)$ and $(\BD \bu^m, \BS^m) \in \mathcal{A}(\cdot)$ a.e. in $Q$ by \eqref{impl-m}, it follows that the $\liminf_{m \to \infty}$ of the above is nonnegative.
To show the other direction, denote $H^m \coloneqq\left( \BS^{m}-\BS^{\star}(\cdot, \BD \bu\right) \colon  \left( \BD \bu^{m} - \BD \bu\right)\geq 0$, and let $j \geq j_0$, $\bad_{m,j} \subset Q_0$ and $\bv^{m,j}$ be given by Lemma~\ref{Lem:LA-par-div} applied on $Q_0$, and by \ref{itm:LA-u-eq} we have that $\bv^{m}=\bv^{m,j}$ on $\frac{1}{8} Q_0 \backslash \bad_{m,j}$. Dividing the domain into $\tfrac{1}{8}Q_0 \cap \bad_{m,j}$ and $\tfrac{1}{8}Q_0 \backslash \bad_{m,j}$, by H\"older's inequality we obtain
\begin{align}\begin{split}\label{lm-impl-11}
\int_{\tfrac{1}{8}Q_0} (H^m)^{\frac{1}{2}} \, \d \bz
&= \int_{\tfrac{1}{8}Q_0 \cap \bad_{m,j}} (H^m)^{\frac{1}{2}} \, \d \bz +
\int_{\tfrac{1}{8}Q_0 \backslash \bad_{m,j}} (H^m)^{\frac{1}{2}} \, \d \bz\\
&\leq
\abs{\tfrac{1}{8}Q_0 \cap \bad_{m,j}}^{\frac{1}{2}}
\left( \int_{\tfrac{1}{8}Q_0 \cap \bad_{m,j}} H^m \, \d \bz\right)^{\frac{1}{2}}\\
&\quad  +
\abs{\tfrac{1}{8}Q_0 \backslash \bad_{m,j}}^{\frac{1}{2}}
\left( \int_{\tfrac{1}{8}Q_0 \backslash \bad_{m,j}} H^m \, \d \bz\right)^{\frac{1}{2}} \\
&
\leq
\abs{ \bad_{m,j}}^{\frac{1}{2}}
\left( \int_{Q} H^m \, \d \bz \right)^{\frac{1}{2}} +
\abs{Q}^{\frac{1}{2}}
\left( \int_{\tfrac{1}{8}Q_0 \backslash \bad_{m,j}} H^m \, \d \bz \right)^{\frac{1}{2}},
\end{split}
\end{align}
where we have used the nonnegativity of $H^m$ in the first term. Since $H^m$ is bounded in $\Leb^1(Q)$ by the a~priori estimate in \eqref{est:m}, one has that
\begin{align}\label{lm-impl-12}
\int_{\tfrac{1}{8}Q_0} (H^m)^{\frac{1}{2}} \, \d \bz
\leq c
\abs{ \bad_{m,j}}^{\frac{1}{2}}
 +
c
\left( \int_{\tfrac{1}{8}Q_0 \backslash \bad_{m,j}} H^m \, \d \bz \right)^{\frac{1}{2}}.
\end{align}
By Lemma~\ref{Lem:LA-par-div} \ref{itm:LA-u-bad} we have that
\begin{align}\label{lm-impl-13}
\limsup_{m \to \infty} \abs{\bad_{m,j}}^{\frac{1}{2}} \leq \limsup_{m \to \infty} (\lambda_{m,j}^{q}\abs{\bad_{m,j}})^{\frac{1}{2}} \leq c 2^{-\frac{j}{2}}.
\end{align}
Let $\zeta \in C^{\infty}_0(\tfrac{1}{6}B_0)$ be the nonnegative function given by Lemma~\ref{Lem:LA-par-div} such that $\zeta|_{\tfrac{1}{8}B_0} \equiv 1$. In the second term in \eqref{lm-impl-12} we can use the nonnegativity of $H^m$, the definition of $H^m$ and $\bv^m$ and finally the definition of $\BG^m_1$ in order to find that $\BS^m = \BS - \BG^m_1$, and we obtain
\begin{align}\begin{split}\label{lm-impl-14}
\int_{\tfrac{1}{8}Q_0 \backslash \bad_{m,j}} H^m \, \d \bz
&=
\int_{\tfrac{1}{8}Q_0 \backslash \bad_{m,j}} H^m \,\zeta \, \d \bz
=
\int_{\tfrac{1}{8}Q_0} H^m \, \zeta \, \mathds{1}_{\bad_{m,j}^c} \, \d \bz \\
&\leq
\int H^m \, \zeta \, \mathds{1}_{\bad_{m,j}^c} \, \d \bz
=
\int \left(\BS^m - \BS^{\star}(\cdot, \BD \bu)\right)\colon \BD \bv^m \, \zeta \, \mathds{1}_{\bad_{m,j}^c} \, \d \bz\\
&=
- \int \left(\BG^m_1 -\BS + \BS^{\star}(\cdot, \BD \bu)\right)\colon \nabla \bv^m \, \zeta \, \mathds{1}_{\bad_{m,j}^c} \, \d \bz.
\end{split}
\end{align}
Since $\BS - \BS^{\star}(\cdot, \BD \bu) \in \Leb^{q'}(Q)^{d \times d}$, we are in the position to use Lemma~\ref{Lem:LA-par-div} \ref{itm:LA-u-GH}. Applying $\limsup_{m \to \infty}$ we find that
\begin{align}\begin{split}\label{lm-impl-15}
\limsup_{m \to \infty} \int_{\tfrac{1}{8}Q_0 \backslash \bad_{m,j}} H^m \, \d \bz
&\overset{\eqref{lm-impl-14}}{\leq}  \limsup_{m \to \infty} \abs{\int \left(\BG^m_1 -\BS + \BS^{\star}(\cdot, \BD \bu)\right)\colon \nabla \bv^m \, \zeta \, \mathds{1}_{\bad_{m,j}^c} \, \d \bz} \\
&\leq  c 2^{-\sfrac{j}{q}}.
\end{split}
\end{align}
Using \eqref{lm-impl-13} and \eqref{lm-impl-15} in  \eqref{lm-impl-12} yields
\begin{align}\begin{split}\label{lm-impl-16}
&\limsup_{m \to \infty} \int_{\tfrac{1}{8}Q_0} \left[\left( \BS^{m}-\BS^{\star}(\cdot, \BD \bu\right) \colon  \left( \BD \bu^{m} - \BD \bu\right)\right]^{\sfrac{1}{2}} \, \d \bz
=
\limsup_{m \to \infty} \int_{\tfrac{1}{8}Q_0} (H^m)^{\sfrac{1}{2}} \, \d \bz\\
&\quad \quad \quad \overset{\eqref{lm-impl-12}}{\leq}
c \limsup_{m \to \infty} \abs{ \bad_{m,j}}^{\frac{1}{2}}
+
c \limsup_{m \to \infty} \left( \int_{\tfrac{1}{8}Q_0 \backslash \bad_{m,j}} H^m \, \d \bz \right)^{\frac{1}{2}}
\\
&\quad \quad \quad \overset{\eqref{lm-impl-13}, \eqref{lm-impl-15}}{\leq}  c (2^{-\frac{j}{2}} + 2^{-\frac{j}{2q}}).\end{split}
\end{align}
Then taking $j \to \infty$ gives the claim and \eqref{lm-impl-10} is proved.
This means that $[H^m]^{\frac{1}{2}}\to 0$ strongly in $\Leb^1(\tfrac{1}{8}Q_0)$, as $m \to \infty$.

 However, to show \eqref{lm-impl-1} we need $\Leb^1$-convergence of $H^m$ at least on suitable subdomains, which can be achieved by use of Chacon's biting lemma, as was done in \cite{BGMS.2012}:
from the strong convergence of $[H^m]^{\frac{1}{2}}\to 0$ in $\Leb^1(\tfrac{1}{8}Q_0)$, as $m \to \infty$, we have that there exists a subsequence such that $[H^m]^{\frac{1}{2}}\to 0$ a.e. in $\tfrac{1}{8}Q_0$ and hence also $H^m\to 0$ a.e. in $\tfrac{1}{8}Q_0$.
 Furthermore, by the above estimates $\{H^m\}_{m \in \mathbb{N}}$ is bounded in $\Leb^1(\tfrac{1}{8}Q_0)$. By Chacon's biting lemma (see \cite{BM.1989}) there exists a further subsequence, a function $H \in \Leb^1(\frac{1}{8}Q_0)$ and a nonincreasing sequence of measurable subsets $E_i \subset \frac{1}{8}Q_0$, $i \in \mathbb{N}$, with $\abs{E_i} \to 0$ as $i \to \infty$, such that
\begin{align*}
H^m \rightharpoonup H\quad \text{ weakly in } \Leb^{1}(\tfrac{1}{8}Q_0\backslash E_i), \quad \text{ as } m \to \infty,\;\;\; \text{ for each fixed } i \in \mathbb{N}.
\end{align*}
Now let $i \in \mathbb{N}$  be arbitrary but fixed. With the Dunford--Pettis compactness criterion (see \cite[Ch.~8, Thm.~1.3]{ET.1999}) it follows that the sequence $\{H^m\}_{m \in \mathbb{N}}$ is equi-integrable on $\tfrac{1}{8}Q_0 \backslash E_i$. By the a.e. convergence of $H^m$ to zero in particular on $\frac{1}{8}Q_0 \backslash E_i$, Vitali's convergence theorem implies that $H^m \to 0$ in $\Leb^1(\tfrac{1}{8}Q_0 \backslash E_i)$, as $m \to \infty$, i.e., we have that
\begin{align}\label{lm-impl-17}
\ll  \BS^{m}-\BS^{\star}(\cdot, \BD \bu),  \BD \bu^{m} - \BD \bu \rr_{\tfrac{1}{8}Q_0 \backslash E_i} \to 0, \quad \text{ as } m \to \infty,
\end{align}
for any fixed $i \in \mathbb{N}$.
With the weak convergence of $\BS^m \rightharpoonup \BS$ in $\Leb^{q'}(Q)^{d \times d}$ by \eqref{cm-4}  and the weak convergence of $\BD \bu^m \rightharpoonup \BD \bu$ in $\Leb^{q}(Q)^{d \times d}$ following from  \eqref{cm-2} we thus deduce that
\begin{align*}
\lim_{m \to \infty} \ll \BS^m, \BD \bu^m \rr_{\tfrac{1}{8}Q_0  \backslash E_i} = \ll \BS, \BD \bu \rr_{\tfrac{1}{8}Q_0  \backslash E_i} \quad \text{ for all } i \in \mathbb{N}.
\end{align*}
This shows \eqref{lm-impl-1} for $\widetilde{Q} = \tfrac{1}{8}Q_0  \backslash E_i$, and thus we find that $(\BD \bu(\bz), \BS(\bz)) \in \mathcal{A}(\bz)$ for a.e. $\bz \in \tfrac{1}{8}Q_0  \backslash E_i$.
Since $\abs{E_i} \to 0$, as $i \to \infty$, we have that $(\BD \bu(\bz), \BS(\bz)) \in \mathcal{A}(\bz)$ for a.e. $\bz \in \tfrac{1}{8}Q_0 $.

Finally let us consider a cover of $Q$ consisting of (open) cylinders $Q^j = I^j \times B^j$, $j \in J$, for an index set $J$ such that
$Q = \bigcup_{j \in J} \tfrac{1}{8}Q^j$.
This can be, for example, chosen as a Whitney type cover, compare, e.g., \cite{DRW.2010}.
Then we can identify the implicit relation a.e. on $\tfrac{1}{8}Q^j$ for all $j \in J$ by the above and thus, have that  $(\BD \bu(\bz), \BS(\bz)) \in \mathcal{A}(\bz)$ for a.e. $\bz \in Q$, which proves  \eqref{impl}.
\end{proof}
\end{Lemma}

\begin{Remark}[Steady Problem]
The same regularization and splitting approach can be applied to show convergence in the steady case to cover the range
$q \in (\frac{2d}{d+2}, \frac{2d}{d+1}]$ for discretely divergence-free finite element functions, which is missing in \cite{DKS.2013}.
In \cite{DKS.2013} no regularization term was used, and hence the restriction $q>\frac{2d}{d+1}$ was required in the case of discretely divergence-free finite element functions. Considering the two results together, for $q \in (\frac{2d}{d+2}, \frac{2d}{d+1}]$ one should either use exactly divergence-free finite element functions as was done in \cite{DKS.2013},
or introduce a regularization term and pass to the limit $m\rightarrow \infty$ with the regularization parameter, as we have done here.
\end{Remark}

\section{Appendix: Results about the Constitutive Laws}


\begin{Lemma}{(Properties of $\BS^k$)}\label{Lem:prop-Sk}
For each $k \in \mathbb{N}$ the function $\BS^k: Q \times \Rds \to \Rds$ defined in \eqref{def:Sk} is measurable with respect to its first argument and smooth with respect to the second argument. Furthermore, the sequence $\{\BS^k\}_{k \in \mathbb{N}}$ satisfies the Assumption~\ref{Assump-Sk} with $\widetilde{g}=g$, $\widetilde{c}_*=c_*$ and $q\in (1,\infty)$, as in Lemma~\ref{sel}.
\begin{proof} Smoothness and measurability follow by the definition of the convolution and Fubini's theorem.
To show the properties \ref{itm:al-2} and \ref{itm:al-3} is straightforward.
For \ref{itm:al-4} we will follow \cite[Sec.~3.2]{BGMS.2012}: 
 Let $\bz \in Q$ be arbitrary but fixed and not in any of the zero-sets for which the properties of $\BS^{\star}$ do not hold. Then we first use the definition of $\BS^k$ in \eqref{def:Sk} and the fact that $\rho^k$ integrates to $1$ on $\Rds$ and then the  monotonicity of  $\BS^{\star}$ in Lemma~\ref{sel} \ref{itm:a2}, which gives that
\begin{align}\label{Sk-a3-0}
(\BS^{\star}(\bz, \BA) -  \BS^{\star}(\bz, \BB))\colon (\BA - \BB )\geq 0 \quad \text{ for all } \BA, \BB \in \Rds,
\end{align}
and the nonnegativity of $\rho^k$ to show that, for any sequence $\{\BD^k\}_{k \in \mathbb{N}}$ bounded in $\Leb^{\infty}(Q)^{d \times d}$ and for arbitrary but fixed $\BB \in \Rds$ we have
 \begin{align}\notag
&\left(\BS^k(\bz, \BD^k(\bz)) - \BS^{\star}(\bz, \BB)\right)\colon \left( \BD^k(\bz) - \BB\right) \\ \notag
 &\quad\quad \overset{\eqref{def:Sk}}{=} \int_{\Rds} ( \BS^{\star}(\bz, \BA) - \BS^{\star}(\bz, \BB)) \colon  (\BD^k(\bz)-\BB) \rho^k(\BD^k(\bz) - \BA) \, \d \BA \\ \label{Sk-a3-1}
&\quad\quad \overset{\eqref{Sk-a3-0}}{\geq} \int_{\Rds} ( \BS^{\star}(\bz, \BA) - \BS^{\star}(\bz, \BB)) \colon  (\BD^k(\bz)-\BA) \rho^k(\BD^k(\bz) - \BA) \, \d \BA \\
&\quad\quad \geq - \int_{\Rds} \abs{ \BS^{\star}(\bz, \BA) - \BS^{\star}(\bz, \BB)} \abs{ \BD^k(\bz)-\BA} \rho^k(\BD^k(\bz) - \BA) \, \d \BA
\notag\\
&\quad\quad \geq - \frac{1}{k}\int_{B_{\sfrac{1}{k}}(\BD^k(\bz))} \abs{ \BS^{\star}(\bz, \BA) - \BS^{\star}(\bz, \BB)} \rho^k(\BD^k(\bz) - \BA) \, \d \BA \eqqcolon (\star),
\notag
\end{align}
where in the last step we have used that by the definition of $\rho^k$ we have
\begin{align*}
\supp(\rho^k(\BD^k(\bz)-\BA))\subset B_{\sfrac{1}{k}}(\BD^k(\bz)),
\end{align*} and hence $\abs{\BD^{k}(\bz)-\BA} \leq \frac{1}{k}$. Also we integrate only over $\BA \in B_{\sfrac{1}{k}}(\BD^k(\bz))$, and since $\BD^k$ is uniformly bounded in $\Leb^{\infty}(\Omega)^{d \times d}$, for those we have
\begin{align}
\abs{\BA}\leq \tfrac{1}{k} + \norm{\BD^k}_{\Leb^{\infty}(Q)} \leq c,
\end{align}
where the constant $c$ is independent of $\bz \in \Omega$ and of $k \in \mathbb{N}$.

Then, since $\BS^{\star}(\bz, \cdot)$ is locally bounded (compare Lemma~\ref{sel} \ref{itm:a5}), there exists a constant $c>0$ independent of $k$ and $\bz$ such that
\begin{align*}
\abs{\BS^{\star}(\bz, \BA)} \leq c \quad \text{ for all } \BA \in B_{\sfrac{1}{k}}(\BD^k(\bz)).
\end{align*}
Since $\BB\in \Rds$ is fixed, again by the local boundedness of $\BS^{\star}$ we can estimate
\begin{align*}
\abs{\BS^{\star}(\bz, \BB)} \leq c(\BB),
\end{align*}
where the constant may depend on $\BB$, but not on $\bz$.
This implies
\begin{align*}
\abs{\BS^{\star}(\bz, \BA) - \BS^{\star}(\bz, \BB)} \leq \abs{\BS^{\star}(\bz, \BA)} + \abs{\BS^{\star}(\bz, \BB)} \leq c(\BB),
\end{align*}
independently of $\bz \in Q$ and for all $\BA \in B_{\sfrac{1}{k}}(\BD^k(\bz))$. With this we can estimate $(\star )$ further by
\begin{align}\begin{split}\label{Sk-a3-2}
(\star) &= - \frac{1}{k}\int_{B_{\sfrac{1}{k}}(\BD^k(\bz))} \abs{ \BS^{\star}(\bz, \BA) - \BS^{\star}(\bz, \BB)} \rho^k(\BD^k(\bz) - \BA) \, \d \BA \\
&\geq - \frac{1}{k}c(\BB)\int_{B_{\sfrac{1}{k}}(\BD^k(\bz))} \rho^k(\BD^k(\bz) - \BA) \, \d \BA
 = - \frac{c(\BB)}{k},\end{split}
\end{align}
where we have again used that $\rho^k$ integrates to $1$.
Now we take \eqref{Sk-a3-1} and \eqref{Sk-a3-2} together, multiply by $\varphi \in C^{\infty}_0(Q)$ such that $\varphi \geq 0$, and then integrate in $Q$. Recalling that the constant is independent of $\bz \in Q$ this gives
\begin{align*}
\int_{Q} (\BS^k(\cdot, \BD^k) - \BS^{\star}(\cdot, \BB))\colon (\BD^k - \BB) \varphi\, \d \bz \geq - \frac{c(\BB)}{k}\norm{\varphi}_{\Leb^1(Q)}.
\end{align*}
Then, applying $\liminf_{k \to \infty}$ yields
 \begin{align*}
\liminf_{k \to \infty} \int_{Q} (\BS^k(\cdot, \BD^k) - \BS^{\star}(\cdot, \BB))\colon (\BD^k - \BB) \varphi\, \d \bz \geq 0,
 \end{align*}
which proves the claim.
\end{proof}
\end{Lemma}

\begin{Lemma}[Properties of $\{\cS^k\}_{k \geq k_0}$]\label{Lem:prop-Sk-2}
The family of functions $\cS^k: Q \times \mathbb{R} \to \mathbb{R}$, $k \geq k_0$, defined in  \eqref{def:ap-Sk} has the following properties: each $\cS^k$ is a Carath\'{e}odory function
and
\begin{itemize}
\item[($\alpha 1$')] for a.e. $\bz \in Q$ the function $\cS^k(\bz, \cdot)$ is monotone;
\item[($\alpha 2$')] there exists a constant $\widetilde{c}_{*}>0$ and a nonnegative $\widetilde{g} \in \Leb^1(Q)$ such that
\begin{align*}
  \cS^{k}(\bz, B) B \geq -\widetilde{g}(\bz)+ \widetilde{c}_*\left(\abs{B}^{q}+ \abs{\cS^{k}(\bz, B)}^{q'}\right),
\end{align*}
for any $B \in \mathbb{R}_{\geq 0}$ and for a.e. $\bz \in Q$ and all $k \geq k_0$;
\item[($\alpha 3$')]  for any sequence $\{D^k\}_{k \in \mathbb{N}}$ bounded in $\Leb^{\infty}(Q)$, for any $B \in \mathbb{R}_{\geq 0}$ and all $\varphi \in C^{\infty}_0(Q)$ such that $\varphi \geq 0$, we have
\begin{align*}
\liminf_{k \to \infty} \int_Q \left( \cS^k(\cdot, D^k)- \cS^{\star}(\cdot, B) \right)\colon \left(D^k - B \right)\varphi \, \d \bz \geq 0.
\end{align*}
\end{itemize}
\begin{proof}
Since $\cS^{\star}(\cdot, \cdot)$ is measurable in the first argument by Lemma~\ref{sel}, so is $\cS^k(\cdot, \cdot)$ and by construction $\cS^k(\cdot, \cdot)$ is continuous in the second argument. Hence, $\cS^k(\cdot, \cdot)$ is a Carath\'{e}odory function.
\begin{itemize}
\item[($\alpha 1$')] By Lemma~\ref{sel} \ref{itm:a2}, $\cS^{\star}(\bz, \cdot)$ is monotone for a.e. $\bz \in Q$. Then, by construction $\cS^{k}(\bz, \cdot)$ is piecewise monotone for a.e. $\bz \in Q$, and consequently monotone for a.e. $\bz \in Q$.
\item[($\alpha 2$')] Let $\bz \in Q$ be arbitrary but fixed such that all properties hold for $\cS^{\star}(\bz, \cdot)$, excluding zero sets where necessary. Also let $B \in \mathbb{R}_{\geq 0}$ and $k \geq k_0$ be arbitrary but fixed.\\
If $B \notin A^k$, then we have that $\cS^k(\bz, B)=\cS^{\star}(\bz, B)$ and hence
the claim holds
by the estimate for $\BS^{\star}(\bz, \cdot)$ in Lemma~\ref{sel} \ref{itm:a3}, which is equivalent to the one for $\cS^{\star}(\bz, \cdot)$, with the same $g$ and $c_*$.\\
If $B \in A^k$, we make use of the fact that $A$ and hence also $A^k$ is bounded. In particular we know that $B \leq a^k_{I,+} \leq a_I + 1$. Hence, again with the estimate for $\cS^{\star}(\bz, \cdot)$ corresponding to Lemma~\ref{sel} \ref{itm:a3} we obtain
\begin{align}\label{est:mon-Sk-1}
\cS^k(\bz, B) B &= (\cS^k(\bz, B)- \cS^{\star}(\bz, B)) B + \cS^{\star}(\bz, B) B\\
&\geq - \left( \abs{\cS^k(\bz, B) } + \abs{\cS^{\star}(\bz, B)} \right) \abs{B} - g(\bz) + c_* \left(\abs{B}^q + \abs{\cS^{\star}(\bz, B)}^{q'}\right).\notag
\end{align}
Now using that $B <a_I+1$, and monotonicity of both $\cS^{\star}(\bz, \cdot)$ and $\cS^{k}(\bz, \cdot)$ we find that
\begin{align}\begin{split}\label{est:mon-Sk-2}
 \abs{\cS^{\star}(\bz, B) } &\leq \abs{\cS^{\star}(\bz, a_I+1) },\\
\abs{\cS^k(\bz, B) } &\leq \abs{\cS^k(\bz, a^k_{I,+}) } = \abs{\cS^{\star}(\bz, a^k_{I,+}) } \leq \abs{\cS^{\star}(\bz, a_I+1) }.\end{split}
\end{align}
Hence in \eqref{est:mon-Sk-1} we obtain
\begin{align}\begin{split}\label{est:mon-Sk-3}
\cS^k(\bz, B) B
&\overset{\eqref{est:mon-Sk-1}}{\geq}  - \left( \abs{\cS^k(\bz, B) } + \abs{\cS^{\star}(\bz, B)} \right) \abs{B} - g(\bz) + c_* \left(\abs{B}^q + \abs{\cS^{\star}(\bz, B)}^{q'}\right)\\
&\overset{\eqref{est:mon-Sk-2}}{\geq}
 - 2 \abs{\cS^{\star}(\bz, a_I+1) } (a_I+1)  - g(\bz) + c_* \left(\abs{B}^q + \abs{\cS^{\star}(\bz, B)}^{q'}\right).\end{split}
\end{align}
In order to estimate $\abs{\cS^{\star}(\bz, B)}^{q'}$ further, we use the inequality
$$
(a+b)^s \leq 2^{s-1} (a^s + b^s), \quad \text{ for } a,b \in \mathbb{R}_{\geq 0}, \, s \in [1, \infty).
$$
By the triangle inequality and again applying  \eqref{est:mon-Sk-2}, this yields
\begin{align*}
\abs{\cS^k(\bz, B)}^{q'}&\leq \left(\abs{\cS^{k}(\bz, B)-\cS^{\star}(\bz, B)} + \abs{\cS^{\star}(\bz, B)}  \right)^{q'}\\
 &\leq 2^{q'-1} \left(\abs{\cS^{k}(\bz, B)-\cS^{\star}(\bz, B)}^{q'} +  \abs{\cS^{\star}(\bz, B)}^{q'} \right)\\
 &\leq 2^{q'-1} \left(\abs{\cS^{k}(\bz, B)}+ \abs{ \cS^{\star}(\bz, B)}\right)^{q'} +   2^{q'-1} \abs{\cS^{\star}(\bz, B)}^{q'} \\
 &\overset{\eqref{est:mon-Sk-2}}{\leq }
 2^{2q'-1}  \abs{\cS^{\star}(\bz, a_I+1) }^{q'}  +2^{q'-1} \abs{\cS^{\star}(\bz, B)}^{q'}.
\end{align*}
Rearranging and applying this in \eqref{est:mon-Sk-3} gives
\begin{align*}
\cS^k(\bz, B) B
\geq& - 2 \abs{\cS^{\star}(\bz, a_I+1) } (a_I+1)  - g(\bz) + c_* \left(\abs{B}^q + \abs{\cS^{\star}(\bz, B)}^{q'}\right)\\
\geq& - 2 \abs{\cS^{\star}(\bz, a_I+1) } (a_I+1)  - g(\bz) - c_* 2^{q'}\abs{\cS^{\star}(\bz, a_I+1)}^{q'}\\
& + c_* \abs{B}^q + 2^{-(q'-1)} c_*\abs{\BS^k(\bz, B)}^{q'}.
\end{align*}
We set
$$
\widetilde{g}(\bz) \coloneqq 2 \abs{\cS^{\star}(\bz, a_I+1) } (a_I+1) +  c_* 2^{q'}\abs{\cS^{\star}(\bz, a_I+1)}^{q'} + g(\bz) \geq g(\bz) \geq 0,
$$
which is in $\Leb^1(Q)$ thanks to the local boundedness of $\cS^{\star}$ and we choose $ \widetilde{c}_{*} \coloneqq 2^{-(q'-1)}c_*$.
\item[($\alpha 3$')] Let $\{D^k\}_{k \in \mathbb{N}}$ be bounded in $\Leb^{\infty}(Q)$ by $\widetilde{C}$, let $B \in \mathbb{R}_{\geq 0}$ and let $\bz \in Q$ be arbitrary but fixed such that all properties hold for $\cS^{\star}(\bz, \cdot)$, ($\alpha 1$') and ($\alpha 2$') hold for $\cS^k(\bz, \cdot)$ and such that $\abs{D^k(\bz)} \leq \widetilde{C}$, excluding possibly zero sets.
We wish to show that there exists a constant $c>0$ depending on $\BS^{\star}$, $a_I$, $B$ and $\widetilde{C}$ such that
\begin{align}\label{est:id-Sk-1}
(\star) \coloneqq \left( \cS^k(\bz, D^k(\bz))- \cS^{\star}(\bz, B) \right)\colon \left(D^k(\bz) - B \right) \geq - \tfrac{c}{k}.
\end{align}
To this end, we distinguish the following three cases.
\begin{itemize}
\item[$B  \notin A^k$:] In this case we have that $\cS^{k}(\bz, B) = \cS^{\star}(\bz, B)$ and by monotonicity of $\cS^k(\bz, \cdot)$, it follows that $(\star) \geq 0$.
\item[$D^k(\bz) \notin A^k$:] Similarly, in this case we have that $\cS^{k}(\bz, D^k(\bz)) = \cS^{\star}(\bz, D^k(\bz))$, so by monotonicity of $\cS^{\star}(\bz, \cdot)$ it follows that $(\star) \geq 0$.
\item[$B, D^k(\bz) \in A^k$:] Assume that $B \in A^k_i$ and $D^k(\bz)\in A^k_j$ and further distinguish the cases
$i=j$ and $i \neq j$.\\
$i=j$: In this case we make use of the fact that $\abs{D^k(\bz)-B}\leq \frac{2}{k}$ and that $A^k$ is bounded, i.e., $B, D^k(z) \leq a_I + 1$. Then, using again the estimates in \eqref{est:mon-Sk-2} we have
\begin{align*}
&\left( \cS^k(\bz, D^k(\bz))- \cS^{\star}(\bz, B) \right)\colon \left(D^k(\bz) - B \right) \\
&\quad \quad \quad \geq
- \left(\abs{ \cS^k(\bz, D^k(\bz))} + \abs{\cS^{\star}(\bz, B)} \right) \abs{D^k(\bz) - B} \\
&\quad \quad \quad \geq
- \left(\abs{ \cS^k(\bz, D^k(\bz))} + \abs{\cS^{\star}(\bz, B)} \right) \tfrac{2}{k} \\
&\quad \quad \quad \overset{\eqref{est:mon-Sk-2}}{\geq}
- 2 \cS^{\star}(\bz, a_I+1) \tfrac{2}{k}
 \geq
-\tfrac{c}{k},
\end{align*}
since $\cS^{\star}$ is locally bounded.\\
$i\neq j$: In this case we make use of the fact that $\cS^k(\bz,\cdot)$ and $\cS^{\star}(\bz, \cdot)$ agree on $a^k_{i,\pm}$, $a^k_{j,\pm}$ and are monotone. Let $i<j$, i.e., $B \leq D^k(\bz)$; it then follows that
\begin{align*}
\cS^{\star}(\bz, B) \leq \cS^{\star}(\bz, a^k_{i,+}) =
\cS^{k}(\bz, a^k_{i,+}) \leq \cS^{k}(\bz, a^k_{j,-}) \leq
\cS^{k}(\bz, D^k(\bz)),
\end{align*}
which implies that
\begin{align*}
\left( \cS^k(\bz, D^k(\bz))- \cS^{\star}(\bz, B) \right)\colon \left(D^k(\bz) - B \right)  \geq 0.
\end{align*}
The case $i>j$ is shown analogously.
\end{itemize}
\end{itemize}
Altogether, these imply \eqref{est:id-Sk-1}, where the right-hand side is independent of $\bz$. Multiplying with $\varphi \in C^{\infty}_0(Q)$ such that $\varphi \geq 0$, and integrating over $Q$ gives
\begin{align*}
\int_Q \left( \cS^k(\cdot, D^k)- \cS^{\star}(\cdot, B) \right)\colon \left(D^k- B \right) \varphi \, \d \bz \geq - \frac{c}{k}.
\end{align*}
Then applying $\liminf_{k \to \infty}$ yields the assertion.
\end{proof}
\end{Lemma}

\begin{Corollary}\label{Cor:prop-Sk-2}
The family of functions $\BS^k: Q \times \Rds \to \Rds $, $k \geq k_0$, defined by
\begin{align*}
    \BS^{k}(\bz, \BB) \coloneqq  \cS^{k}(\bz, \abs{\BB})\frac{\BB}{\abs{\BB}}, \quad \text{ for } \bz \in Q, \,\BB \in \Rds,
\end{align*}
where $\cS^k$ is defined in \eqref{def:ap-Sk}, satisfies Assumption~\ref{Assump-Sk}.
\begin{proof}
 It is straightforward to show that monotonicity of $\cS^k$ implies monotonicity of $\BS^k$ and that the growth and coercivity bounds are equivalent for $\cS^k$ and $\BS^k$.
\end{proof}
\end{Corollary}

\bibliographystyle{plain}

\bibliography{library}

\begin{thebibliography}{10}

\bibitem{AF.1988}
E.~Acerbi and N.~Fusco.
\newblock An approximation lemma for {$W^{1,p}$} functions.
\newblock In {\em Material instabilities in continuum mechanics ({E}dinburgh,
  1985--1986)}, Oxford Sci. Publ., pages 1--5. Oxford Univ. Press, New York,
  1988.

\bibitem{BM.1989}
J.M. Ball and F.~Murat.
\newblock Remarks on {C}hacon's {B}iting {L}emma.
\newblock {\em Proceedings of the American Mathematical Society},
  107(3):655--663, 1989.

\bibitem{BBDR.2012}
L.~Belenki, L.C. Berselli, L.~Diening, and M.~R{\r{u}}{\v{z}}i{\v{c}}ka.
\newblock On the finite element approximation of $p$-{S}tokes systems.
\newblock {\em SIAM Journal on Numerical Analysis}, 50(2):373--397, 2012.

\bibitem{BL.1976}
J.~Bergh and J.~L\"ofstr\"om.
\newblock {\em Interpolation spaces. {A}n introduction}.
\newblock Springer-Verlag, Berlin-New York, 1976.
\newblock Grundlehren der Mathematischen Wissenschaften, No. 223.

\bibitem{BBF.2013}
D.~Boffi, F.~Brezzi, and M.~Fortin.
\newblock {\em Mixed Finite Element Methods and Applications}, volume~44 of
  {\em Springer Series in Computational Mathematics}.
\newblock Springer, 2013.

\bibitem{B.1979}
M.E. Bogovski\u\i.
\newblock Solution of the first boundary value problem for an equation of
  continuity of an incompressible medium.
\newblock {\em Dokl. Akad. Nauk SSSR}, 248(5):1037--1040, 1979.

\bibitem{BDF.2012}
D.~Breit, L.~Diening, and M.~Fuchs.
\newblock Solenoidal {L}ipschitz truncation and applications in fluid
  mechanics.
\newblock {\em Journal of Differential Equations}, 253(6):1910--1942, 2012.

\bibitem{BDS.2013}
D.~Breit, L.~Diening, and S.~Schwarzacher.
\newblock Solenoidal {L}ipschitz truncation for parabolic {PDE}'s.
\newblock {\em Mathematical Models and Methods in Applied Sciences},
  23(14):2671--2700, 2013.

\bibitem{BM.2016}
M.~Bul{\'i}{\v{c}}ek and J.~M{\'a}lek.
\newblock {\em Recent Developments of Mathematical Fluid Mechanics}, chapter On
  Unsteady Internal Flows of Bingham Fluids Subject to Threshold Slip on the
  Impermeable Boundary, pages 135--156.
\newblock Springer Basel, Basel, 2016.

\bibitem{BGMS.2009}
M.~Bul\'\i{}\v{c}ek, P.~Gwiazda, J.~M\'alek, and A.~\'{S}wierczewska{-}Gwiazda.
\newblock On steady flows of incompressible fluids with implicit power-law-like
  rheology.
\newblock {\em Advances in Calculus of Variations}, 2:109--136, February 2009.

\bibitem{BGMS.2012}
M.~Bul\'\i{}\v{c}ek, P.~Gwiazda, J.~M\'alek, and A.~\'{S}wierczewska{-}Gwiazda.
\newblock On unsteady flows of implicitly constituted incompressible fluids.
\newblock {\em SIAM Journal on Mathematical Analysis}, 44:2756--2801, 2012.

\bibitem{CHP.2010}
E.~Carelli, J.~Haehnle, and A.~Prohl.
\newblock Convergence analysis for incompressible generalized {N}ewtonian fluid
  flows with nonstandard anisotropic growth conditions.
\newblock {\em SIAM Journal on Numerical Analysis}, 48(1):164--190, 2010.

\bibitem{C.2002}
P.~Ciarlet.
\newblock {\em The Finite Element Method for Elliptic Problems}.
\newblock Society for Industrial and Applied Mathematics, 2002.

\bibitem{CR.1973}
M.~Crouzeix and P.-A. Raviart.
\newblock Conforming and nonconforming finite element methods for solving the
  stationary {S}tokes equations {I}.
\newblock {\em ESAIM: Mathematical Modelling and Numerical Analysis -
  Modélisation Mathématique et Analyse Numérique}, 7(R3):33--75, 1973.

\bibitem{D.1993}
E.~DiBenedetto.
\newblock {\em Degenerate Parabolic Equations}.
\newblock Springer Verlag, 1993.

\bibitem{DKS.2013}
L.~Diening, C.~Kreuzer, and E.~S\"uli.
\newblock Finite element approximation of steady flows of incompressible fluids
  with implicit power-law-like rheology.
\newblock {\em SIAM J. Numer. Anal.}, 51(2):984--1015, 2013.

\bibitem{DMS.2008}
L.~Diening, J.~M\'alek, and M.~Steinhauer.
\newblock On {L}ipschitz truncations of {S}obolev functions (with variable
  exponent) and their selected applications.
\newblock {\em ESAIM: Control, Optimisation and Calculus of Variations},
  14(2):211--232, 3 2008.

\bibitem{DRW.2010}
L.~Diening, M.~R\r{u}\v{z}i\v{c}ka, and J.~Wolf.
\newblock Existence of weak solutions for unsteady motions of general
  {N}ewtonian fluids.
\newblock {\em Ann. Scuola Norm. Sup. Pisa Cl. Sci.}, 9(1):1--46, 2010.

\bibitem{ET.1999}
I.~Ekeland and R.~Temam.
\newblock {\em Convex analysis and variational problems}.
\newblock North-Holland Publishing Co., Amsterdam-Oxford; American Elsevier
  Publishing Co., Inc., New York, 1976.
\newblock Translated from the French, Studies in Mathematics and its
  Applications, Vol. 1.

\bibitem{EG.2004}
A.~Ern and J.-L. Guermond.
\newblock {\em Theory and Practice of Finite Elements}, volume 159 of {\em
  Applied Mathematical Sciences}.
\newblock Springer-Verlag New York, 1st edition, 2004.

\bibitem{FMS.2003}
J.~Frehse, J.~M\'{a}lek, and M.~Steinhauer.
\newblock On analysis of steady flows of fluids with shear-dependent viscosity
  based on the {L}ipschitz truncation method.
\newblock {\em SIAM Journal on Mathematical Analysis}, 34(5):1064--1083, 2003.

\bibitem{GT.2001}
D.~Gilbarg and N.S. Trudinger.
\newblock {\em Elliptic Partial Differential Equations of Second Order}.
\newblock Classics in Mathematics. Springer, 2001.
\newblock reprint, originally in Grundlehren der mathematischen Wissenschaften,
  Vol. 224.

\bibitem{GL.2001b}
V.~Girault and J.-L. Lions.
\newblock Two-grid finite-element schemes for the steady {N}avier-{S}tokes
  problem in polyhedra.
\newblock {\em Port. Math. (N.S.)}, 58(1):25--57, 2001.

\bibitem{GS.2003}
V.~Girault and L.R. Scott.
\newblock A quasi-local interpolation operator preserving the discrete
  divergence.
\newblock {\em Calcolo}, 40(1):1--19, 2003.

\bibitem{GD.2003}
A.~Granas and J.~Dugundji.
\newblock {\em Fixed point theory}.
\newblock Springer Monographs in Mathematics. Springer-Verlag, New York, 2003.

\bibitem{Gr.2011}
P.~Grisvard.
\newblock {\em Elliptic Problems in Nonsmooth Domains}, volume~69 of {\em
  Classics in Applied Mathematics}.
\newblock Society for Industrial and Applied Mathematics, 2011.

\bibitem{GN.2014b}
J.~Guzm\'{a}n and M.~Neilan.
\newblock Conforming and divergence-free {S}tokes elements in three dimensions.
\newblock {\em IMA Journal of Numerical Analysis}, 34(4):1489--1508, 2014.

\bibitem{GN.2014}
J.~Guzm\'an and M.~Neilan.
\newblock Conforming and divergence-free {S}tokes elements on general
  triangular meshes.
\newblock {\em Math. Comp.}, 83:15--36, 2014.

\bibitem{KL.2002}
J.~Kinnunen and J.L. Lewis.
\newblock Very weak solutions of parabolic systems of $p$-{L}aplacian type.
\newblock {\em Ark. Mat.}, 40:105--132, 2002.

\bibitem{KS.2016}
C.~Kreuzer and E.~S\"uli.
\newblock Adaptive finite element approximation of steady flows of
  incompressible fluids with implicit power-law-like rheology.
\newblock {\em ESAIM Math. Model. Numer. Anal.}, 50(5):1333--1369, 2016.

\bibitem{La.1969}
O.A. Lady\v{z}enskaja.
\newblock {\em The mathematical theory of viscous incompressible flow}.
\newblock Second English edition, revised and enlarged. Translated from the
  Russian by Richard A. Silverman and John Chu. Mathematics and its
  Applications, Vol. 2. Gordon and Breach, Science Publishers, New
  York-London-Paris, 1969.

\bibitem{L.1969}
J.-L. Lions.
\newblock {\em Quelques M\'{e}thodes De R\'{e}solution Des Probl\`{e}mes Aux
  Limites Non Lin\'{e}aires}.
\newblock Dunod, Paris, 1969.

\bibitem{R.2003}
K.R. Rajagopal.
\newblock On implicit constitutive theories.
\newblock {\em Applications of Mathematics}, 48(4):279--319, 2003.

\bibitem{R.2006}
K.R. Rajagopal.
\newblock On implicit constitutive theories for fluids.
\newblock {\em Journal of Fluid Mechanics}, 550:243--249, 3 2006.

\bibitem{Ro.2013}
T.~Roub\'{\i}\v{c}ek.
\newblock {\em Nonlinear partial differential equations with applications},
  volume 153 of {\em International Series of Numerical Mathematics}.
\newblock Birkh\"auser/Springer Basel AG, Basel, 2nd edition, 2013.

\bibitem{R.1991}
W.~Rudin.
\newblock {\em Functional analysis}.
\newblock International Series in Pure and Applied Mathematics. McGraw-Hill,
  Inc., New York, 2nd edition, 1991.

\bibitem{S.1987}
J.~Simon.
\newblock Compact sets in the space {$L^p(0,T;B)$}.
\newblock {\em Ann. Mat. Pura Appl. (4)}, 146:65--96, 1987.

\bibitem{S.1999}
J.~Simon.
\newblock On the existence of the pressure for solutions of the variational
  {N}avier-{S}tokes equations.
\newblock {\em J. Math. Fluid Mech.}, 1(3):225--234, 1999.

\bibitem{T.1984}
R.~Temam.
\newblock {\em Navier--Stokes Equations: Theory \& Numerical Analysis}.
\newblock Studies in Mathematics and Its Applications. North-Holland, 1984.

\end{thebibliography}

\end{document}